\pdfoutput=1
\documentclass[journal ]{new-aiaa}
\usepackage[utf8]{inputenc}
\usepackage{textcomp}

\usepackage{graphicx}
\usepackage{gensymb}
\usepackage{amsmath}
\usepackage[version=4]{mhchem}
\usepackage{siunitx}
\usepackage{longtable,tabularx}
\graphicspath{{../Figures/}}

\usepackage{setspace}
\usepackage{tocloft}
\usepackage{caption}
\usepackage{subcaption}
\usepackage{tikz}
\usepackage{tkz-euclide}
\usepackage{flafter}
\usepackage{placeins}

\def\dxo{\delta\mathbf{x}(t_0)}
\def\dxom#1{\delta\mathbf{x}^{#1}(t_0)}

\def\ph#1#2#3{\boldsymbol{\Phi}(#1;t_{#2},t_{#3})}
\def\ps#1#2#3{\boldsymbol{\Psi}(#1;t_{#2},t_{#3})}
\def\psm#1#2#3#4{\boldsymbol{\Psi}^{(#1)}(#2;t_{#3},t_{#4})}

\def\xo{\mathbf{x}(t_0)}

\usepackage[normalem]{ulem}

\title{Applications of Induced Tensor Norms to Guidance Navigation and Control}

\author{Jackson Kulik \footnote{PhD Candidate, Center for Applied Mathematics, Email: jpk258@cornell.edu}, Cedric Orton-Urbina \footnote{Undergraduate Researcher, Physics}, Maximilian Ruth \footnote{PhD Candidate, Center for Applied Mathematics}, and Dmitry Savransky\footnote{Associate Professor, Mechanical and Aerospace Engineering}}
\affil{Cornell University, Ithaca, NY, 14850.}

\begin{document}

\maketitle

\begin{abstract}
    Linear methods are ubiquitous for control and estimation problems. In this work, we present a number of tensor operator norms as a means to approximately bound the error associated with linear methods and determine the situations in which that maximum error is encountered. An emphasis is placed on induced norms that can be computed in terms of matrix or tensor eigenvalues associated with coefficient tensors from higher-order Taylor series. These operator norms can be used to understand the performance and range of applicability of an algorithm exploiting linear approximations in different sets of coordinates. %or in some instances to provide heuristics for extensions of linear methods such as Gaussian mixture models. 
    We examine uses of tensor operator norms in the context of linear and higher-order rendezvous guidance, coordinate selection for a filtering measurement model, and to present a unified treatment of nonlinearity indices for dynamical systems.  %We find error bounds for the linear guidance methods around circular orbits in the 2-body problem as well as around the proposed NASA Gateway Near Rectilinear Halo Orbit. Next, we use the tensor norm methodology to compare measurement error associated with extended Kalman filter measurement updates when using azimuth/elevation angles over a unit vector for optical measurements. Finally, we present a number of novel nonlinearity indices as well as streamlined computational methods for existing nonlinearity indices. Another sentence
    Tensor norm computations can offer insights into these problems in one to two orders of magnitude less time than similarly accurate sampling methods while providing more general understanding of the error performance of linear or higher-order approximations.
\end{abstract}

\section{Introduction}
Methods relying on linearization are widespread in the fields of control and estimation due to their efficiency as well as ease of implementation and analysis. In particular, resource constrained embedded systems such as those on satellites often necessitate the efficiency of linear algorithms for onboard guidance and navigation tasks. %Examples of linear methods in dynamics, control, and estimation are too numerous and fundamental to list; so i
In this work, we focus on two applications: calculating rendezvous or stationkeeping guidance using state transition matrices\cite{clohessy1960terminal,mullins1992initial}, as well as the measurement update of the extended Kalman filter with an optical measurement \cite{schmidt1981kalman}. We present a method for approximating the maximum possible error when using some of these linear algorithms when the true function or dynamical system is nonlinear. Error estimates are expressed as a function of the scale of deviations from some nominal or reference values used in the linearization. 
%When using linear methods or developing semi-analytical methods that improve upon existing linear methods, it is important to understand the error behavior of linear methods. It becomes important not only to understand the maximum error, but what situations or inputs lead to that maximum error. Thus 

%%%%%stretching
%Additionally, the theory of linear systems is also useful in designing guidance and stationkeeping strategies for satellite trajectories through the examination of the singular value decomposition of the linearized dynamics about a reference trajectory. This leads to local lyapunov exponent based methods and maximal or minimal stretching direction methods  \cite{anderson2003application,muralidharan2023stretching,rivera2024statistical}. These linear stretching directions are generalized by higher-order stretching behavior which examines constrained optimization similar to what arises in the calculation of the singular value decomposition, but in the context of some higher-order truncated Taylor series about a reference orbit \cite{jenson2024bounding}. We will review and extend this higher-order stretching approach as the computational methods for this problem bear great similarity to those used to treat the previous problem, and allow us to cast the tensor eigenvalue problems solvers used in the previous section as a type of minorization maximization algorithm.

In relation to this idea about the error performance of linear algorithms, there has been much interest in the last few decades in quantifying the level of nonlinearity of a function such as the flow of a dynamical system. Measuring the nonlinearity of the flow of a dynamical system makes it possible to compare coordinate systems and choose the coordinate system with the best linear approximation before applying a linear control or estimation algorithm \cite{junkins1997karman, junkins2004nonlinear, jenson2022semianalytical}. Similarly, measures of nonlinearity and other related information about directions of maximal nonlinearity can aid development of semi-analytical algorithms that split the domain or employ Gaussian mixtures to improve linearization behavior for each split \cite{losacco2024low, tuggle2020model}. We present a number of nonlinearity indices that can be computed in terms of tensor norms, including a scale-free generalization of the original nonlinearity index that is easy to compute using the induced 2-norm of a matrix. Additionally, we develop an easier to implement approach for computing the recently developed nonlinearity index Tensor Eigenpair Measure of Nonlinearity (TEMoN)  \cite{jenson2022semianalytical} that requires lower-order tensors as compared with the existing approach and only a very simple power iteration like algorithm. Finally, we develop a novel measure of nonlinearity that we call D-Eigenvalue Measure of Nonlinearity (DEMoN).% which has a favorable interpretation as approximating the norm of an error as opposed to an error of a norm in the case of TEMoN.

We approach the problems described above by extending the formalism of induced or operator norms from matrices to multilinear operators or tensors. This allows us to treat a number of existing algorithms in the literature as well as other possible algorithms in terms of a single unified framework. In particular, this includes the tensor eigenvalue approach \cite{jenson2022semianalytical, jenson2024bounding}, the polynomial bounding scheme on a hypercube \cite{losacco2024low}, and the maximum change in the Frobenius norm of the first-order derivative approach \cite{tuggle2020model}. We define an induced tensor norm and present a number of examples in addition to these three. Much of the focus of this paper is on the induced tensor 2-norm, which relies on tensor Z-eigenvalues for computation, as well as the induced tensor $(2,\mathbf{D})$-norm which relies on a generalization known as a tensor D-eigenvalue \cite{qi2008d}. We begin by summarizing a number of tensors that arise in guidance, navigation, and control in Sec. \ref{sec:tensors}, build up from vector to tensor norms, describing their computation in Sec. \ref{sec:norms}, and finally apply tensor norms to the bounding of error in guidance (Sec. \ref{sec:guidance}) and estimation (Sec. \ref{sec:measurement}) contexts, as well as in nonlinearity index computations (Sec. \ref{sec:nonlin_index}).

\section{Tensors in Guidance Navigation and Control}
In this section, we summarize four varieties of tensors that play important roles in guidance, navigation, and control.
\label{sec:tensors}
\subsection{Local Dynamics Tensors}
Given an autonomous dynamical system in $\mathbb{R}^n$, the state vector $\mathbf{x}\in\mathbb{R}^n$ evolves according to the system of ordinary differential equations

\begin{equation}
    \frac{d}{dt}\mathbf{x}=\mathbf{F(x)}, \quad \mathbf{x}(0)=\mathbf{x}_0
    \label{eqn:dyn}
\end{equation}
Local dynamics tensors are defined as the partial derivative tensors of the vector field $\mathbf{F(x)}$ evaluated at some point~$\mathbf{x}^*$: 
\begin{equation}
    \left[\left(\frac{\partial^m \mathbf{F}}{\partial \mathbf{x}^m}\right)\bigg\rvert_{\mathbf{x}^*}\right]^i_{j_1...j_m}=\left(\frac{\partial F^i}{\partial x^{j_1}...\partial x^{j_m}}\right)\bigg\rvert_{\mathbf{x}^*}
\end{equation}
The $m$th order local dynamics tensor is a $(1,m)$-tensor, where an $(l,m)$-tensor has $m$ covariant indices and $l$ contravariant indices. An interpretation of this is that the tensor operates on $m$ vectors or $m$ copies of a single vector to produce one vector as a result.
Since a dynamical system evolves so that $\mathbf{x}^*$ changes over time, local dynamics tensors are most meaningfully interpreted at equilibria of that dynamical system where they describe the dynamical behavior around the equilibrium point. Analysis of the eigenvalues of these tensors has been discussed in previous literature  \cite{jenson2022semianalytical}, and is not expanded upon here. However, local dynamics tensors are necessary for calculating partial derivative tensors of the flow of the dynamical system.

\subsection{State Transition Tensors}
State transition tensors are the partial derivative tensors of the flow of a dynamical system evaluated at some reference trajectory (associated with some initial condition and some time of propagation). The flow map associated with the dynamical system in Eq.~\ref{eqn:dyn} is defined such that \begin{equation}
\frac{d}{dt}\varphi_t(\mathbf{x}_0)=\mathbf{F}(\varphi_t(\mathbf{x}_0)), \quad \varphi_0(\mathbf{x}_0)=\mathbf{x}_0
\end{equation}
The Jacobian of the flow map yields the state transition matrix (STM) $\boldsymbol{\Phi}(\mathbf{x}_0; t_f,t_0)$ associated with a given flow starting at the reference state $\mathbf{x}_0$ from time $t_0$ to time $t_f$ where $\Delta t=t_f-t_0$. We adopt indexing for the state transition matrix rows $i$ and columns $j$
\begin{equation}
    \Phi^i_j(\mathbf{x}_0;t_f,t_0)=\frac{\partial \varphi_{\Delta t}^i(\mathbf{x}_0)}{\partial x_0^j}
\end{equation}
where the upper index $i$ of the flow map refers to the $i$th component of the output. When the reference orbit initial state $\mathbf{x}_0$ is understood or the time span $t_0$ to $t_f$ is understood, the state transition matrix is be abbreviated as $\boldsymbol{\Phi}(t_f,t_0)$ or simply $\boldsymbol{\Phi}$.
Exchanging the order of temporal and spatial derivatives (assuming $\mathbf{F}$ has continuous spatial derivatives) and applying the chain rule yields the $n^2$ first-order variational equations 
\begin{equation}
    \frac{d\boldsymbol{\Phi}(t_f,t_0)}{dt_f}=\frac{\partial \mathbf{F(x)}}{\partial \mathbf{x}}\boldsymbol{\Phi}(t_f,t_0), \quad \boldsymbol{\Phi}(t_0,t_0)=\mathbf{I}_n
    \label{eqn:first_var}
\end{equation}
where $\mathbf{I}_n$ is the $n$ by $n$ identity matrix. Alternatively, in components
\begin{equation}
    \frac{d\Phi^i_j(t_f,t_0)}{dt_f}=\frac{\partial F_i(\mathbf{x})}{\partial x_l} \Phi^l_j(t_f,t_0), \quad \Phi^i_j(t_0,t_0)=\delta^i_{j}
    \label{eqn:phi1}
\end{equation}
where summation with respect to $l$ is understood and $\delta^i_{j}$ is the Kronecker delta.
Moving on to the second-order partial derivatives, we define the second-order $(1,2)$-state transition tensor (STT) $\boldsymbol{\Psi}(t_f,t_0)$
\begin{equation}
    \Psi^i_{j,k}(\mathbf{x}_0;t_f,t_0)=\frac{\partial^2 \varphi_t^i(\mathbf{x}_0)}{\partial x_0^j\partial x_0^k}
\end{equation}
Again, when the reference orbit initial state $\mathbf{x}_0$ is understood, the notation $\boldsymbol{\Psi}(t_f,t_0)$ is used to denote the second-order state transition tensor. Applying the product rule to Eq.~\ref{eqn:phi1}, we find that the $n^3$ equations known as the second-order variational equations depend not only on the values of $\boldsymbol{\Psi}$ but also $\boldsymbol{\Phi}$. Assuming that $\mathbf{F}$ has continuous second-order spatial partial derivatives, the second-order variational equations are
\begin{equation}
    \frac{d\Psi^i_{j,k}(t_f,t_0)}{dt}=\frac{\partial^2 F_i(\mathbf{x})}{\partial x_l \partial x_q} \Phi^l_j(t_f,t_0)\Phi^q_k(t_f,t_0)+\frac{\partial F_i(\mathbf{x})}{\partial x_l} \Psi^l_{j,k}(t_f,t_0),\quad \Psi^i_{j,k}(t_0,t_0)=0
    \label{eqn:second_var}
\end{equation}
With the STM and second-order STT, a second-order approximation of a perturbation to the flow map is given by a truncated Taylor series:
\begin{equation}
	\varphi_t^i(\mathbf{x}_0 + \mathbf{\delta x}_0) \approx \varphi_t^i(\mathbf{x}_0) + \Phi_j^i(\mathbf{x}_0;t_f,t_0) \delta x_0^j + \frac{1}{2} \Psi^i_{jk}(\mathbf{x}_0;t_f,t_0)\delta x_0^j \delta x_0^k
	\label{eqn:taylor}
\end{equation}
or, in shorthand to represent the single and double contractions:
\begin{equation}
	\varphi_t(\mathbf{x}_0 + \mathbf{\delta x}_0) \approx \varphi_t(\mathbf{x}_0) + \mathbf{\Phi}(\mathbf{x}_0; t_f,t_0) \delta \mathbf{x}_0 + \frac{1}{2} \mathbf{\Psi}(\mathbf{x}_0; t_f,t_0)\delta \mathbf{x}_0^2
\end{equation}
Note here that the superscript 2 above a bolded vector placed next to a bolded tensor indicates double contraction between the tensor and two copies of the vector. Indexing of vectors with a superscript is denoted instead with a non-bolded scalar as seen in Eq. \ref{eqn:taylor}. At times, indexing may take place with a bolded vector or tensor expression surrounded first by parentheses. By looking at the context, it should be clear when Einstein notation is being used by the matching indices for contractions as well as the non-bolded notation. In bolded instances, it should be assumed that this multiple contraction definition is used for brevity and to facilitate understanding.

Higher-order state transition tensors and variational equations can be derived under additional assumptions about the smoothness of the vector field \cite{boone2022directional}. We represent higher-order state transition tensors with the notation $\mathbf{\Psi}^{(m)}$ or when in components with $m$ covariant indices $\Psi^i_{j_1...j_m}$. The $m$-th order Taylor series is written in similar shorthand as:
\begin{equation}
    \varphi_t(\mathbf{x}_0 + \mathbf{\delta x}_0) \approx \varphi_t(\mathbf{x}_0) + \sum_{p=1}^m \frac{1}{p!}\mathbf{\boldsymbol{\Psi}}^{(p)}\delta \mathbf{x}_0^p
\end{equation}
Some applications of state transition tensors are summarized in the applications section of this paper; however, more information on higher-order state transition tensors, their analytical computation \cite{kellyhigher}, and their applications for dynamics \cite{rein2016second,bani2019exact,cunningham2023interpolated}, control \cite{boone2021orbital, kulik2023state}, and navigation/uncertainty propagation \cite{park2006nonlinear,majji2008high, boodramefficient, kulik2024overdetermined} exists across the literature of the last few decades. Parallel to the state transition tensor approach, differential algebra (DA) \cite{rasotto2016differential} is a technique for computing the Taylor expansion of an arbitrary function whether it be the flow of a dynamical system or another arbitrary nonlinear function. While DA is not employed in this paper, computing state transition tensors using higher-order variational equations becomes increasingly complicated at each order. As such, the DA approach offers a more convenient and typically efficient approach for computing the Taylor series up to a given order, from which the coefficients can be trivially collected into the corresponding partial derivative tensors such as state transition tensors. DA has been employed for uncertainty propagation \cite{valli2013nonlinear,wittig2015propagation}, filtering \cite{servadio2022maximum}, orbit determination \cite{armellin2018probabilistic}, and control purposes among others \cite{di2014high,lizia2008application}. Many of the algorithms described in the DA literature leverage efficient computations for manipulating formal power series and approximating the inverse of a truncated series. These methods can be used to efficiently obtain tensors related to and arising from the state transition tensor series such as the inverse flow of the system or the higher-order Cauchy-Green strain tensors described below.

\subsection{Higher-Order Cauchy-Green Strain Tensors}
Higher-order Cauchy-Green strain tensors are the coefficient tensors arising in the series expansion for the squared 2-norm of the final perturbation of a dynamical system as a function of some initial perturbation. These tensors were introduced to examine the stretching behavior of a dynamical system around some reference trajectory as well as the associated nonlinearity index \cite{jenson2022semianalytical,jenson2024bounding}
\begin{align}
    \delta\mathbf{x}^T(t_f)\delta\mathbf{x}(t_f) = \sum_{m=2}^\infty\mathbf{C}^{(m)}\dxo^m=\sum_{ p,q\geq1}^\infty \frac{1}{p!q!}\left(\boldsymbol{\Psi}^{(p)}\delta\mathbf{x}^p(t_0)\right)^T\left(\boldsymbol{\Psi}^{(q)}\delta\mathbf{x}^q(t_0)\right)
\end{align}
where $\mathbf{C}^{(m)}$ is the $m$th order (right) Cauchy-Green strain tensor. This series converges assuming that the flow map is infinitely differentiable.
%since it arises from the Taylor series expansion of the squared norm of the flow map differenced twice with a constant term (the flow of the reference trajectory) multiplied by the non-constant terms of the Taylor series expansion of the flow (a convergent series itself). 
Note that the higher-order Cauchy-Green strain tensors are covariant and not generally supersymmetric except in the case of the first-order 2-tensor. These tensors can be calculated in component form in terms of the Euclidean metric tensor $g_{i,j}=\delta_{i,j}$ which is the Kronecker delta:
\begin{align}
    C_{i_1...i_m}=\sum_{p+q=m; p,q\geq1} \frac{1}{p!q!}\Psi^{j_1}_{i_1...i_p}\delta_{j_1,j_2}\Psi^{j_2}_{i_{p+1}...i_m}
    \label{eq:cgt-tensor}
\end{align}
The first Cauchy-Green strain tensor is given by the square of the state transition matrix $\mathbf{C}^{(2)}=\mathbf{\Phi}^T\mathbf{\Phi}$, where the matrix is understood to represent a covariant (0,2)-tensor. %, though in matrix notation the purely covariant nature of the tensor is lost without the metric in Eq.~\ref{eq:cgt-tensor} to ``juggle" the indices. 
Notably, the singular value decomposition of the state transition matrix yields right singular vectors that are the same as the eigenvectors of the Cauchy-Green strain tensor, and singular values that are the same as the square root of the eigenvalues of the Cauchy-Green strain tensor. In a similar fashion, tensor eigenvalues of the higher-order Cauchy-Green strain tensors describe the stretching behavior of the higher-order terms in the state transition tensor series. For study with tensor eigenvalues, the supersymmetric tensor denoted $\hat{\mathbf{C}}^{(m)}$ is often used in calculations and surrounding theory. This is a tensor which produces the same result as $\mathbf{C}^{(m)}$ when operating on the same $m$ vectors, but is supersymmetric. Details about symmetrization, as well as avoiding symmetrization in calculations can be found in Sec. \ref{sec:2-norm} and Appendix \ref{appendix:partial_sym}.

\subsection{Measurement Partial Derivative Tensors}
Given a $d$-dimensional measurement function $h(x):\mathbb{R}^n\rightarrow\mathbb{R}^{d}$, the partial derivative tensors of this function play a role in some higher-order estimation algorithms such as the $j$th-moment Kalman filter \cite{majji2008jth} as well as in the context of assessing the error in using linearizations of measurement function for measurement underweighting schemes in Kalman filters \cite{majji2008jth, zanetti2010underweighting}. The order $(1,m)$ partial derivative tensors arise in the Taylor series expansion of the measurement function about some prior estimated state $\mathbf{x}^-$
\begin{equation}
    \mathbf{h}(\mathbf{x}^-+\delta\mathbf{x})=\mathbf{h}(\mathbf{x}^-)+\sum_{m=1}^\infty \frac{1}{m!}\left( \frac{\partial^m\mathbf{h}}{\partial \mathbf{x}^m}\right)_{\mathbf{x}=\mathbf{x}^-} \delta \mathbf{x}^m
\end{equation}
In this work, we particularly look at the first- and second-order partial derivatives from the series above.

%\subsection{Distribution Moment Tensors}
%Beyond partial derivative tensors which play a role in Taylor series of vector fields, flows, and measurement functions, tensors also play a role in representing the statistical moments of a distribution. While the sample moments/empirical moments of a distribution can be calculated from data coming from measurement or Monte Carlo study, it is often also convenient to consider analytical propagation of statistical moments through dynamics and measurement functions via their partial derivative tensors. A summary of these techniques can be found in the literature \cite{park2006nonlinear, boone2023directional}. We will represent the purely covariant tensor $\mathbf{P^(m)}$ $C^(m)$
%TODO: moments vs cumulants

\section{Norms of Tensors}
\label{sec:norms}
A norm $\Vert\cdot\Vert$ on a vector space $\mathcal{X}$ is defined as a function from the vector space to the real numbers that satisfies three properties: the triangle inequality $\Vert \mathbf{x} +\mathbf{y}\Vert\leq \Vert\mathbf{x}\Vert + \Vert\mathbf{y}\Vert$ for all $\mathbf{x},\mathbf{y}\in\mathcal{X}$, absolute homogeneity with scalars $\Vert a\mathbf{x}\Vert=\vert a\vert \Vert\mathbf{x}\Vert$ for all scalars $a$ and all vectors $\mathbf{x}\in\mathcal{X}$, and finally positive definiteness $\Vert\mathbf{x}\Vert=0$ if and only if $\mathbf{x}=0$ is the zero element of the vector space $\mathcal{X}$. Note that tensors of a given order and dimension along each index form a vector space. This includes Euclidean vectors, matrices, as well as higher-order tensors. We present a variety of norms of Euclidean vector spaces, as well as on the vector spaces associated with matrices and higher-order tensors, building up from simpler to more complex.

\subsection{Vector Norms}
The norm of a vector quantifies some notion of the size of that vector. We review a number of useful vector norms here. First, the $p$-norms of a vector $\mathbf{v}\in\mathbb{R}^n$ are defined for $p\geq 1$ as
\begin{align}
    \lVert \mathbf{v}\rVert_p = \left(\sum_{i=1}^n \lvert v_i\rvert^p\right)^{1/p}
\end{align}
and the $\infty$-norm is defined as
\begin{equation}
    \lVert \mathbf{v}\rVert_\infty=\max_{i=1...n} \lvert v_i\rvert
\end{equation}
The 2-norm is independent of the choice of coordinates for a vector unlike all of the other $p$-norms. A vector norm can also be induced by a symmetric positive definite matrix $\mathbf{D}$ using a quadratic form:
\begin{equation}
    \lVert \mathbf{v}\rVert_\mathbf{D}=\mathbf{v}^T\mathbf{D}\mathbf{v}
\end{equation}
The 2-norm of a vector is one such norm where $\mathbf{D}=\mathbf{I}_n$.
\subsection{Frobenius Norm}
The Frobenius norm of a tensor is easily computable and is equivalent to any other finite dimensional norm (bounded in both directions by potentially dimension dependent constant multiples of the other norm) \cite{golub2013matrix}. The Frobenius norm of an $(l,m)$-tensor under the Euclidean metric is the square root of the sum of the squared entries:
\begin{equation}
    %\left\Vert \mathbf{B} \right\Vert_F=\sqrt{\sum_{i,j,k...}\left(B_{i,j,k,...}\right)^2}
    \left\Vert \mathbf{B} \right\Vert_F^2=B^{i_1...i_{l}}_{j_1...j_m}\delta_{i_1,i'_1}...\delta_{i_l,i'_l}\delta^{j_1,j'_1}...\delta^{j_l,j'_l}B^{i'_1...i'_{l}}_{j'_1...j'_m}
\end{equation}
where $\delta^{i,j}$ with superscripts denotes the Euclidean inverse metric tensor given again by the Kronecker delta but expressed as a contravariant instead of covariant tensor. 
In general, the Frobenius norm of any tensor is the 2-norm of the tensor flattened into a vector in any order.

\subsection{Induced/Operator Norms}
These vector and Frobenius norms, among others, can be used to define other norms. The induced or operator norm subordinate to the $a$ and $b$ norms, or $(a, b)$-norm of a matrix $\mathbf{A}\in\mathbb{R}^{m,n}$ is defined as
\begin{equation}
    \Vert\mathbf{A}\Vert^a_b = \max_{\mathbf{x}\neq \mathbf{0}}\frac{\Vert\mathbf{A} \mathbf{x}\Vert_a}{\Vert\mathbf{x}\Vert_b}=\max_{\Vert\mathbf{x}\Vert_b=1}\Vert\mathbf{A} \mathbf{x}\Vert_a
\end{equation}
where the $a$-norm and $b$-norm conform to the dimensions of $\mathbf{A}$ and $\mathbf{x}\in\mathbb{R}^n$ \cite{golub2013matrix}.
We may extend this definition of an operator norm to tensors. Let $\mathbf{B}$ be a $(1,m)$-tensor such that the tensor is partially symmetric along the covariant indices. Let the $a$-norm on the output space be defined for a $(1,m_o)$-tensor, and the $b$-norm be some vector norm on the input space. Then the induced norm of a tensor is:
\begin{equation}
    \Vert\mathbf{B}\Vert^a_b = \max_{\mathbf{x}\neq \mathbf{0}}\frac{\Vert\mathbf{B} \mathbf{x}^{m-m_o}\Vert_a}{\Vert\mathbf{x}\Vert_b^{m-m_o}}=\max_{\Vert\mathbf{x}\Vert_b=1}\Vert\mathbf{B} \mathbf{x}^{m-m_o}\Vert_a
\end{equation}
where 
\begin{equation}
    \left(\mathbf{B} \mathbf{x}^{m-m_o}\right)^i_{j_1,...j_{m_o}}=B^i_{j_1...j_m}x^{j_{m_o+1}}...x^{j_{m}}
\end{equation}
When $a,b$ are both vector norms as is usually the case, $m_o=0$ and the induced tensor norm is defined:
\begin{equation}
    \Vert\mathbf{B}\Vert^a_b = \max_{\mathbf{x}\neq \mathbf{0}}\frac{\Vert\mathbf{B} \mathbf{x}^{m}\Vert_a}{\Vert\mathbf{x}\Vert_b^{m}}=\max_{\Vert\mathbf{x}\Vert_b=1}\Vert\mathbf{B} \mathbf{x}^{m}\Vert_a
\end{equation}
An example where $m_o$ is not zero is the induced (Frobenius, 2)-norm, where the $a$ norm is the Frobenius norm on a $(1,1)$-tensor. See Sec. \ref{sec:frob2} for details.
%these are not norms and may not be computable with power iteration even when the square of another polynomial. unsure.
%In full generality, suppose we have an operator representing a multivariate nonhomogeneous polynomial of order $m$:
%\begin{equation}
%    \mathcal{P}(x)=\sum_{i=1}^m\mathbb{P}^(i)
%\end{equation}

\subsection{2-norm}
\label{sec:2-norm}
The induced 2-norm of a $(1,m)$-tensor $\mathbf{B}$ is perhaps the most natural to consider given that it relies only on a single basic vector norm for its definition and it is unitarily invariant so the norm is invariant under rotations of the coordinates. We can compute the 2-norm of a tensor by examining the eigenvalues of the "square" of the tensor. Again, given the metric tensor for Euclidean space $g_{i,j}=\delta_{i,j}$ which is the Kronecker delta, we define the $(0,2m)$-tensor "square," $\tilde{\mathbf{B}}$:

\begin{equation}
\Tilde{B}_{i_1,..,i_m,j_1,...,j_m}=B^\alpha_{i_1,..,i_m}\delta_{\alpha,\beta}B^\beta_{j_{1},..,j_m}
\end{equation}

This square of the tensor gives the squared 2-norm of the output of the original tensor when applied to a vector ($2m$ times):

\begin{equation}
\Tilde{\mathbf{B}}\mathbf{x}^{2m}=\Vert\mathbf{B}\mathbf{x}^m\Vert^2_2
\end{equation}

For every covariant tensor such as $\Tilde{\mathbf{B}}$, there exists a symmetric tensor $\hat{\mathbf{B}}$ with the same output under repeated contraction with a single vector. That is, $\hat{\mathbf{B}}$ exists such that for any $\mathbf{x}\in\mathbf{R}^n$,

\begin{equation}
\hat{\mathbf{B}}\mathbf{x}^{2m}= \Tilde{\mathbf{B}}\mathbf{x}^{2m}
\end{equation}

and 
\begin{equation}
\hat{B}_{i_1,...,i_m,j_1,...,j_m}=\hat{B}_{\sigma(i_1,...,i_m,j_1,...,j_m)}
\end{equation}
for any permutation $\sigma$ of the indices. %The induced 2-norm of the tensor $\mathbf{B}$ is defined as the square root of the maximum Z-eigenvalue of the corresponding symmetrized square tensor.
%\begin{equation}
%    \Vert\mathbf{B}\Vert_2=\sqrt{\lambda_1(\hat{\mathbf{B}})}
%\end{equation}
%where a Z-eigenvalue is one of many generalizations of eigenvalueds for tensors. In particular, a Z-eigenvalue, eigenvector pair are defined in terms of a unit-norm constraint for the eigenvector:
%\begin{equation}
%    \mathbf{B}\mathbf{x}^{(m-1)}=\lambda \mathbf{x} \quad \text{s.t.} \quad \lVert \mathbf{x}\rVert_2=1
%\end{equation}
%The maximum eigenvalue of an even order symmetric tensor can easily be calculated by a higher-order generalization of power iteration. This method is globally convergent to some (typically large) eigenvalue, but may require multiple initial guesses to converge to the largest eigenvalue, eigenvector pair if an educated guess is not available \cite{kolda2011shifted}.
The symmetrization is equal to the mean of the tensor under all permutations:
\begin{equation}
\hat{B}_{i_1,...,i_m,j_1,...,j_m}=\frac{1}{(2m)!}\sum_{\sigma\in S_{2m} }\tilde{B}_{\sigma(i_1,...,i_m,j_1,...,j_m)}.
\end{equation}
where $S_{2m}$ is the permutation group of order $2m$.
Another more computationally efficient method for symmetrization appears in the MATLAB Tensor Toolbox \cite{kolda2006matlab}, though as we will demonstrate, explicit calculation of the symmetrized tensor is rarely necessary or desirable.
Finally, having obtained a symmetric positive semi-definite tensor ($\mathbf{B}\mathbf{x}^{2m}\geq 0$ for all $\mathbf{x}\in\mathbb{R}^n$), the maximum of the following constrained optimization problem can be obtained as the largest Z-eigenvalue of the tensor \cite{lim2005singular}:
\begin{equation}
    \lambda_\text{max}(\hat{\mathbf{B}})=\max_{\Vert\mathbf{x}\Vert_2=1} \hat{\mathbf{B}}\mathbf{x}^{2m}
\end{equation}
where a Z-eigenvalue satisfies the equation
\begin{equation}
    \hat{\mathbf{B}}\mathbf{x}^{2m-1}=\lambda \mathbf{x}, \quad \Vert\mathbf{x}\Vert_2=1
\end{equation}
More generally, any of the Z-eigenvectors of $\hat{\mathbf{B}}$ satisfy the Karush–Kuhn–Tucker (KKT) conditions and are constrained stationary points for the objective function $\hat{\mathbf{B}}\mathbf{x}^{2m}$ when $\mathbf{x}$ is constrained to the unit sphere \cite{lim2005singular}. Thus, the induced 2-norm of the tensor $\mathbf{B}$ is defined as the square root of the maximum Z-eigenvalue of the corresponding symmetrized square tensor $\hat{\mathbf{B}}$
\begin{equation}
    \Vert\mathbf{B}\Vert_2=\sqrt{\lambda_\text{max}(\hat{\mathbf{B}})}.
\end{equation}
This is a generalization of the idea that the induced 2-norm of a matrix is given by the largest singular value (square root of the largest eigenvalue of the square of the matrix). The maximum eigenvalue of an even-order symmetric tensor that is convex over the unit ball can be calculated by a higher-order generalization of power iteration \cite{de2000best, kolda2011shifted, kofidis2002best}.
\begin{equation}
    \mathbf{x}_{n+1}=\frac{\hat{\mathbf{B}}\mathbf{x}^{2m-1}_n}{\Vert\hat{\mathbf{B}}\mathbf{x}^{2m-1}_n\Vert}_2
    \label{eqn:power_iter}
\end{equation}
%Under this globally convergent scheme, $\mathbf{x}_n$ will always converge to the Z-eigenvector corresponding to one of the (typically largest) Z-eigenvalues for any initial guess $\mathbf{x}_0$. 
%In practice, the eigenvector corresponding to the largest eigenvalue is almost always obtained for random initial conditions; however, it is prudent to perform higher-order symmetric power iteration with a number of initial conditions to ensure the eigenvalue associated with the largest eigenvalue is obtained.
Having obtained the desired eigenvector, the corresponding eigenvalue is given by the magnitude of the numerator of Eq.~\ref{eqn:power_iter}, or $\Vert\mathbf{B}\Vert_2$ can be computed directly by taking
\begin{equation}
    \Vert\mathbf{B}\Vert_2 = \Vert\mathbf{B}(\mathbf{x}^*)^m\Vert_2
\end{equation}
where $\mathbf{x}^*$ is the Z-eigenvector that power iteration has converged to.
Since $\hat{\mathbf{B}}$ is positive semi-definite, even order, and fully symmetric, the symmetric higher-order power iteration method is globally convergent meaning it is guaranteed to converge to some eigenvector regardless of the starting guess $\mathbf{x}_0$. Typically, the eigenvector that is converged to corresponds to a large eigenvalue \cite{kofidis2002best}. This means that generally, an educated guess must be supplied for the initial iterate (like the dominant right singular vector of the state transition matrix in the case of a higher-order state transition tensor \cite{boone2023directional}) or multiple random initial guesses should be employed when trying to find the eigenvector associated with the largest eigenvalue. The number of random guesses needed to find the global maximum depends on the gap between the largest and second largest eigenvalue since the size of each basin of convergence varies with the size of the eigenvalue. Thus, when two eigenvalues are close to being the dominant eigenvalue, one might need tens to hundreds instead of a couple initial random guesses to reliably find the dominant eigenvector \cite{kolda2011shifted} %, though in practice the authors have rarely if ever noticed convergence to anything but the largest eigenpair due to its typically large basin of convergence. 
In other contexts for more general tensors, symmetric higher-order power iteration does not have global convergence guarantees, and shifted symmetric higher-order power iteration must be employed \cite{kolda2011shifted}. While symmetric higher-order power iteration requires a fully symmetric tensor such as $\hat{\mathbf{B}}$ to converge properly to an eigenvector, the partially symmetric structure of $\mathbf{B}$ and $\tilde{\mathbf{B}}$ enable us to perform symmetric higher-order power iteration on $\hat{\mathbf{B}}$ in Eq.~\ref{eqn:power_iter} without ever needing to perform the symmetrization to form $\hat{\mathbf{B}}$. Further performance optimizations are also possible by directly using the $\mathbf{B}$ tensor directly in calculations rather than its square. We refer to Appendix \ref{appendix:partial_sym} for further details.

\subsection{$(2,\mathbf{D})$-norm}
\label{sec:d-norm}
The notion of a D-eigenvalue was introduced as a generalization to the normalization in the definition of Z-eigenvalues for use in higher-order statistics \cite{qi2008d}. A D-eigenvalue of the tensor $\hat{\mathbf{B}}$ satisfies the equation
\begin{equation}
    \hat{\mathbf{B}}\mathbf{x}^{2m-1}=\lambda \mathbf{D}\mathbf{x}, \quad \Vert\mathbf{x}\Vert_\mathbf{D}=1
    \label{eqn:d-eig}
\end{equation}
where $\mathbf{D}$ is a symmetric positive definite matrix and the constraint enforces that eigenvectors lie on the ellipsoid given by the quadratic form induced by $\mathbf{D}$ rather than on the unit sphere as is the case with a Z-eigenvalue. All D-eigenvectors of $\hat{\mathbf{B}}$ satisfy the Karush–Kuhn–Tucker (KKT) conditions for the objective function $\hat{\mathbf{B}}\mathbf{x}^{2m}$ when $\mathbf{x}$ is constrained to the $\mathbf{D}$-ellipsoid of points satisfying the constraint in Eq.~\ref{eqn:d-eig}.
Generalized eigenproblem adaptive power method (GEAP) may be used to solve D-eigenvalue problems 
\cite{kolda2014adaptive}. However, since we are only interested in the largest eigenpair, we propose another simple and fast extension of symmetric higher-order power iteration which is usable because the tensors involved are convex multilinear operators over the unit ball. We are interested in the largest D-eigenpair in order to solve the problem
\begin{equation}
    \max_{\mathbf{x}^T\mathbf{D}\mathbf{x}=1}\hat{\mathbf{B}}\mathbf{x}^{2m}
\end{equation}
This problem can be changed into a Z-eigenvalue problem with the change of variables
\begin{equation}
    \mathbf{y}=\mathbf{D}^{1/2}\mathbf{x}
\end{equation}
where $\mathbf{D}^{1/2}$ is a matrix square root of the matrix $\mathbf{D}$ such as the Cholesky factor. That is
\begin{equation}
    \mathbf{D}=\left(\mathbf{D}^{1/2}\right)^T\mathbf{D}^{1/2}
\end{equation}
Making this substitution, the optimization problem becomes
\begin{equation}
     \max_{\mathbf{y}^T\mathbf{y}=1}\hat{\mathbf{B}}(\mathbf{D}^{-1/2}\mathbf{y})^{2m}  
\end{equation}
where $\mathbf{D}^{-1/2}$ denotes the inverse of the matrix square root of $\mathbf{D}$.
This constrained optimization leads to the Z-eigenvalue problem associated with the tensor
\begin{equation}
    \left(\hat{B}_\mathbf{D}\right)_{i_1...i_{2m}}=\hat{B}_{j_1...j_{2m}}\left(\mathbf{D}^{-1/2}\right)^{j_1}_{i_1}...\left(\mathbf{D}^{-1/2}\right)^{j_{2m}}_{i_{2m}}
    \label{eqn:d-eig-tensor}
\end{equation}
If $(\mathbf{v}_\mathbf{D},\lambda_{\mathbf{D}})$ is a Z-eigenpair of $\hat{\mathbf{B}}_\mathbf{D}$, then the D-eigenpair associated with $(\hat{\mathbf{B}},\mathbf{D})$ from Eq.~\ref{eqn:d-eig} is given by the inverse coordinate change
%\begin{equation}
%    (\mathbf{v},\lambda)=\left(\mathbf{D}^{-1/2}\mathbf{v}_\mathbf{D}, \frac{\Vert\hat{\mathbf{B}}(\mathbf{D}^{-1/2}\mathbf{v}_\mathbf{D})^{2m-1}\Vert_2}{\Vert\mathbf{D}^{1/2}\mathbf{v}_\mathbf{D}\Vert_2}\right)
%\end{equation}
\begin{equation}
    (\mathbf{v},\lambda)=\left(\mathbf{D}^{-1/2}\mathbf{v}_\mathbf{D}, \lambda_\mathbf{D}\right)
\end{equation}
since
\begin{equation}
    \lambda=\lambda_\mathbf{D}=\frac{\hat{\mathbf{B}}\mathbf{v}^{2m}}{\Vert\mathbf{v}\Vert_\mathbf{D}^2}=\hat{\mathbf{B}}\mathbf{v}^{2m}
\end{equation}
Note that this correspondence implies that the largest Z-eigenvalue corresponds to the largest D-eienvalue, and to find the largest D-eigenvalue, one can find the largest Z-eigenvalue of the associated problem. The $(2,\mathbf{D})$-norm of $\mathbf{B}$ is given in terms of the D-eigenvalue of the square of the tensor $\mathbf{B}$ as
\begin{equation}
    \Vert\mathbf{B}\Vert^{[2]}_\mathbf{D}=\sqrt{\lambda_\mathrm{max}(\hat{\mathbf{B}},\mathbf{D})}=\sqrt{\lambda_\mathrm{max}(\hat{\mathbf{B}}_\mathbf{D})}
\end{equation}
where the superscript ``$[2]$" in the above expression is not a square of the norm, but indicates the type of norm (a $(2,\mathbf{D})$-norm).
Rather than explicitly forming the tensor $\hat{\mathbf{B}}_\mathbf{D}$ and finding the Z-eigenpairs using symmetric higher-order power iteration, a similar procedure to that proposed in Sec. \ref{appendix:partial_sym} can be used after precomputing the Cholesky decomposition of $\mathbf{D}$ a single time prior to iteration. This algorithm takes only a few more operations than solving for the tensor 2-norm, with one computation of the Cholesky decomposition of a matrix prior to iteration (an $\mathcal{O}(n^3)$ operation), and one additional forward substitution and another backwards substitution to solve a lower and upper triangular linear system at each iteration (each taking only $\mathcal{O}(n^2)$ operations where $n$ is the dimension of the vector $\mathbf{x}$). To put this in perspective, one iteration of higher order power iteration on a covariant order $(2m)$-tensor naively takes $\mathcal{O}(n^{2m})$ or $\mathcal{O}(n^{m+1})$ when the tensor arises as the square of some partially symmetric $(1,m)$-tensor and the calculation is performed according to Appendix \ref{appendix:partial_sym}. This implies that computing a D-eigenvalue is negligibly slower than computing a Z-eigenvalue. Details can be found in Appendix Sec. \ref{appendix:partial_sym_d}.

\subsection{$(\infty,\, 2)$-norm}

In addition to the 2-norm, the induced $(\infty,\, 2)$-norms of both the matrix $\mathbf{A}$ and the $(1,\, 2)$-tensor $\mathbf{B}$ have analytical forms that are readily computable with standard linear algebra tools. The intuition for this development comes from viewing a tensor that is partially symmetric in the covariant indices (such as the second-order state transition tensor) as a vector of symmetric bilinear forms. The $(\infty,\, 2)$-norm of a matrix is well-known. 
\begin{align}
        \nonumber \Vert\mathbf{A}\Vert^\infty_2 &= 
    \max_{\Vert\mathbf{x}\Vert_2=1} \Vert\mathbf{A} \mathbf{x}\Vert_\infty\\
    &=\max_i \max_{\Vert\mathbf{x}\Vert_2=1} A^i_j\mathbf{x}^j\\
    %&=\max_i \frac{(\mathbf{A}^T)^j_i\mathbf{A}^i_j}\\
    \nonumber &=\max_i \Vert \mathbf{A}^i\Vert_2
\end{align}
To our knowledge, the computation of the $(\infty,\, 2)$-norm for a $(1,2)$-tensor that is symmetric along the covariant indices is a novel presentation:
\begin{align}
    \nonumber \Vert\mathbf{B}\Vert^\infty_2 &= \max_{\Vert\mathbf{x}\Vert_2=1} \Vert\mathbf{B} \mathbf{x}^2\Vert_\infty\\
    %&= \max_{\Vert\mathbf{x}\Vert_2=1} \max_i B^i_{j,k} \mathbf{x}^j\mathbf{x}^k\\
   &= \max_i \max_{\Vert\mathbf{x}\Vert_2=1} \left| B^i_{j,k} \mathbf{x}^j\mathbf{x}^k \right|\\    
    &= \max_i \sigma_1(\mathbf{B}^i)
\end{align}
where $\mathbf{A}^i$ denotes the $i$th row vector of the matrix $\mathbf{A}$, and $\sigma_1(\mathbf{B}^i)$ denotes the largest singular value of the covariant matrix $\mathbf{B}^i$ with entries such that $(\mathbf{B}^i)_{j,k}=B^i_{j,k}$. This appears as a result of the well-known identity that the maximum value of a symmetric bilinear form on the 2-norm unit ball is the largest absolute eigenvalue (also the largest singular value) of the matrix defining the bilinear form \cite{golub2013matrix}. When computing the $(\infty,2)$-norm, a solver should be used that is specifically designed to take advantage of the structure of the symmetric eigenvalue problem. We note that the $(\infty,2)$-norm, along with other popular induced norms for a matrix, can be found throughout the linear algebra literature \cite{drakakis2009calculation, lewis2010top}.

\subsection{(Frobenius,\, 2)-norm}
\label{sec:frob2}
The induced (Frobenius,\,2)-norm comes from interpreting the (1,\,2)-tensor as operating on a vector to produce matrices. As we discuss in following sections about nonlinearity indices, this norm can be employed to quantify the change in the Jacobian as the reference point it is evaluated at changes. For example, this norm applied to the second-order state transition tensor quantifies the changes to the state transition matrix as the initial reference state varies. The (Frobenius,\,2)-norm is defined
\begin{equation}
    \left\Vert \mathbf{B} \right\Vert^F_2=\max_{\Vert\mathbf{x}\Vert_2=1}\Vert\mathbf{B} \mathbf{x}\Vert_F
\end{equation}
This norm may be computed by flattening the tensor to the $n^2$ by $n$ matrix
\begin{equation}
    \overline{B}^{ni+j}_k=B^i_{j,k}
\end{equation}
Then the desired norm is given by the induced 2-norm of the matrix resulting from this unfolding: 
\begin{equation}
    \left\Vert \mathbf{B} \right\Vert^F_2=\Vert \overline{\mathbf{B}} \Vert_2
\end{equation}
%Implicitly there is some index juggling taking place in the above expressions that is avoided for clarity in the case where the Euclidean metric is employed. 
This norm is also unitarily invariant like the 2-norm. Another equivalent method to compute the (Frobenius, 2)-norm is through an eigendecomposition in the following manner derived in the context of Gaussian mixture splitting \cite{tuggle2020model}. Note that
\begin{align}
    \Vert\mathbf{Bx}\Vert_F^2&=B^{i_1}_{j_1,k_1}x^{k_1} \delta_{i_1, i_2}\delta^{j_1,j_2}B^{i_2}_{j_2,k_2}x^{k_2}\\
    &=B'_{k_1,k_2}x^{k_1}x^{k_2}
\end{align}
where we define
\begin{equation}
    B'_{k_1,k_2}=B^{i_1}_{j_1,k_1} \delta_{i_1, i_2}\delta^{j_1,j_2}B^{i_2}_{j_2,k_2}
\end{equation}
Ignoring covariance and contravariance of the indices of the matrix $\mathbf{B}'$ we obtain an expression more similar to that derived in \cite{tuggle2020model} without the tensor formalism
\begin{equation}
        \mathbf{B}'=\sum_i (\mathbf{B}^i)^T\mathbf{B}^i
\end{equation}
where $\mathbf{B}^i$ is the matrix given by fixing the contravariant index $i$ and leaving the two covariant indices free.
Then we can express the (Frobenius, 2)-norm as the maximum eigenvalue associated with $\mathbf{B}'$:
\begin{equation}
    \left\Vert \mathbf{B} \right\Vert^F_2=\lambda_{\mathrm{max}}(\mathbf{B}')
\end{equation}
The eigenvectors associated with $\mathbf{B}'$ and equivalently the right singular vectors associated with $\bar{\mathbf{B}}$ are all orthogonal with respect to one another. These offer directions along which the Frobenius norm is at a constrained stationary point. The (Frobenius, $\mathbf{D}$)-norm generalizes the (Frobenius, 2)-norm along the same lines in which the $(2, \mathbf{D})$-norm generalizes the 2-norm, though we do not present details here for space.

\subsection{Upper Bound on the 2-norm}
\label{sec:bounds-2}
By unfolding the tensor along another axis to preserve the first dimension instead of the last dimension, we can arrive at an upper bound on the induced $2$-norm of a tensor 
\cite{li2020eigenvalue}. This bound on the induced 2-norm may be computed by flattening the tensor to an $n$ by $n^{m-1}$ matrix. For a $(1,2)$-tensor, for example, the flattened tensor becomes
\begin{equation}
    \underline{B}^{i}_{nj+k}=B^i_{j,k}
\end{equation}
And tensor contraction becomes equivalent to multiplication with a flattened $n^2$ dimensional vector coming from the outer product of $\mathbf{x}$ with itself
\begin{equation}
    \underline{B}\text{vec}(\mathbf{x}\mathbf{x}^T)=B \mathbf{x}^2
\end{equation}
where 
\begin{equation}
(\text{vec}(\mathbf{x}\mathbf{x}^T))^{(nj+k)}=(\mathbf{x}\mathbf{x}^T)^{j,k}=x^jx^k
\end{equation}
With this relationship, we can show that if the true unit vector satsifying the constrained optimization for the induced 2-norm of a $(1,2)$-tensor is denoted
\begin{equation}
    \mathbf{x}^*=\arg\max_{\lVert\mathbf{x}\rVert_2=1}\lVert \mathbf{B} \mathbf{x}^2\rVert_2
\end{equation}
then, the inequality
\begin{equation}
    \lVert\underline{\mathbf{B}}\rVert_2\geq\lVert\underline{\mathbf{B}}\text{vec}(\mathbf{x}^*(\mathbf{x}^*)^T)\rVert_2=\lVert\mathbf{B} (\mathbf{x}^*)^2\rVert_2
\end{equation}
follows from the definition of the induced 2-norm of a matrix (the norm of the matrix is at least as large as the norm of the matrix operating on any single unit vector) and the fact that the vector being operated on in the above inequality is a unit vector
\begin{equation}
    \lVert \mathbf{x}\rVert_2=1 \implies \lVert\text{vec}(\mathbf{x}\mathbf{x}^T)\rVert_2 = 1
    \label{eq:unit_vec_id}
\end{equation}
The implication in Eq. \ref{eq:unit_vec_id} can be proved from the following observation
\begin{equation}
    \lVert\mathbf{x}\mathbf{x}^T\rVert_F^2=\lVert\text{vec}(\mathbf{x}\mathbf{x}^T)\rVert_2^2=\sum_{i=1}^n\sum_{j=1}^n(x^i)^2(x^j)^2=\left(\sum_{i=1}^n(x^i)^2\right)\left(\sum_{j=1}^n(x^j)^2\right)=\lVert\mathbf{x}\rVert_2^4
\end{equation}
These arguments boil down to the observation that the unit ball/sphere for vectors in $n$ dimensions mapped into $n^2$ dimensions by $\mathbf{x}\rightarrow\text{vec}(\mathbf{x}\mathbf{x}^T)$ is a subset of the unit ball/sphere in $n^2$ dimensions. This upper bound on the 2-norm is also a norm itself. Each necessary property in the definition of a norm follows from the related property of the matrix 2-norm of the flattened tensor.

\subsection{Upper Bound on the (Frobenius, $\infty$)-norm}
\label{sec:frob-inf}
Another induced norm that comes from interpreting a $(1,2)$-tensor as operating on a vector to produce matrices is the (Frobenius, $\infty$)-norm. 
\begin{align}
    \left\Vert \mathbf{B} \right\Vert^F_\infty&=\max_{\Vert\mathbf{x}\Vert_\infty=1}\Vert\mathbf{B} \mathbf{x}\Vert_F\\
\end{align}
While a convenient method to calculate this norm is not immediately apparent, an upper bound can be obtained very efficiently. A nonlinearity index using this bound on the (Frobenius, $\infty$)-norm of the second-order state transition tensor was introduced for automatic domain splitting using differential algebra \cite{losacco2024low}. 
Given a $(1,2)$-tensor $\mathbf{B}$, define the matrix $\mathbf{B}^{\infty}$
\begin{align}
\label{eq:inf-mat}
    (\mathbf{B}^{\infty})^i_j=\sum_{k}\vert B^i_{j,k}\vert
\end{align}
Then the induced (Frobenius, $\infty$)-norm of $\mathbf{B}$ is bounded above by the Frobenius norm of the resulting matrix
\begin{align}
    \left(\sum_i\sum_j\left(\sum_{k}\vert B^i_{j,k}\vert\right)^2\right)^{1/2}=\Vert\mathbf{B}^\infty\Vert_F \geq \left\Vert \mathbf{B} \right\Vert^F_\infty
\end{align}
Note that the original reference uses the subscript 2 on a norm of a matrix to indicate the Frobenius norm and not the induced 2-norm of the matrix \cite{losacco2024low}. While this is an upper bound for the induced (Frobenius, $\infty$)-norm, it also satisfies all the properties to be a norm in its own right, though not an induced norm in terms of any commonly used norms. The ease of calculating this norm makes it attractive. However, neither norm discussed in this section is unitarily invariant, making them dependent on the specific choice of coordinates. This coordinate dependence is attractive in the context of automatic domain splitting along basis vectors, but can be less attractive in other contexts where the physics or engineering applications are less tied to a specific coordinate frame. 

%\subsection{Constrained Maximum 2-norm of a Nonhomogeneous Multilinear Operator}
%\label{sec:mm-opt}

%\cite{de2011majorizing} \cite{kolda2011shifted}\cite{hunter2004tutorial}

% \subsection{Matrix Norms}
% \subsection{Tensors Norms}
% \subsection{Nonhomogeneous Multilinear Operator Norms}

\section{Applications}
\label{sec:applications}
In the following applications, we consider Keplerian two-body dynamics as well as circular restricted three-body dynamics, which we define here.
\subsection{Dynamics}
\subsubsection{Two-Body Dynamics}
The equations of motion for the two-body problem in an inertial frame are
\begin{align}
    \mathbf{F(x)}=\begin{bmatrix} \dot{x}& \dot{y}& \dot{z}& \frac{-\mu x}{\Vert\mathbf{r}\Vert^{3}}& \frac{-\mu y}{\Vert\mathbf{r}\Vert^{3}}& \frac{-\mu z}{\Vert\mathbf{r}\Vert^{3}}\end{bmatrix}^T
    \label{eqn:two-body}
\end{align}
where the state vector is $\mathbf{x}=[x,y,z,\dot{x},\dot{y},\dot{z}]^T$, an overdot denotes a time derivative, $\mu$ represents the standard gravitational parameter for the central body, and $\mathbf{r}=[x, y, z]^T$ is the position vector from the central body to the satellite \cite{vallado2001fundamentals}.
\subsubsection{Circular Restricted Three Body Dynamics}
The equations of motion for the circular restricted three-body problem are given in the synodic frame as 
\begin{align}
    \mathbf{F(x)}=\begin{bmatrix} \dot{x}& \dot{y}& \dot{z}& 2\dot{y}+\frac{\partial \overline{U}}{\partial x}& -2\dot{x}+\frac{\partial \overline{U}}{\partial y}& \frac{\partial \overline{U}}{\partial z}\end{bmatrix}^T
    \label{eqn:cr3bp}
\end{align}
where $\overline{U}(x,y,z)=\dfrac{1-\mu^*}{||\mathbf{r}_1||}+\dfrac{\mu^*}{||\mathbf{r}_2||}+\dfrac{x^2+y^2}{2}$ is the effective potential given the reduced mass $\mu^*=\dfrac{m_2}{m_1+m_2}$ for the two primary bodies with mass $m_1,m_2$ such that $m_1\geq m_2$ located along the x-axis at $[-\mu^*,0,0]$ and $[1-\mu^*, 0, 0]$ with respect to their common barycenter at the origin. The position of the satellite of interest with respect to the two primary bodies is given by $\mathbf{r}_1$ and $\mathbf{r}_2$ respectively \cite{koon2000dynamical}. The reference orbit used in the following sections comes from the proposed NASA Gateway orbit \cite{NationalAA2019, cunningham2023interpolated}. Initial conditions for the orbit are 
\begin{equation}
    \mu=1.0/(81.30059 + 1.0), \quad x_0=1.022022, \quad z_0 = -0.182097, \quad \dot{y}_0 = -0.103256
\end{equation}
in nondimensional units with other initial coordinates equal to zero. This initial condition is at apolune of the orbit. The period of the orbit is 1.511111 nondimensional time units where $2\pi$ time units is equivalent to the revolution period of the Earth-Moon system.

\subsection{Approximate Error Bounds for State Transition Tensor Guidance and Control}
\label{sec:guidance}
Rendezvous and proximity operations as well as formation flight of satellites often relies on linear or higher-order approximations of the dynamics in the vicinity of some chief satellite or reference orbit to facilitate computations regarding guidance and control. The state transition matrix and higher-order state transition tensors can be used to approximately propagate the relative motion of some secondary deputy satellite relative to the chief satellite or other reference orbit. In addition to being useful for forward propagation of dynamics, state transition tensors can be used to solve inverse problems. Calculating the initial velocity required to transfer to some position relative to the original orbit or relative to another satellite is one such boundary value problem that can be solved in the linear case by inversion of some block of the state transition matrix, and in the higher-order case by using series reversion methods. Calculating the initial velocity for a rendezvous with the chief satellite falls into this class of inverse problems.

These linear and higher-order series reversions are attractive for onboard implementation of autonomous guidance, navigation, and control systems due to the efficiency of these methods. One drawback to solving linear or higher-order approximations of boundary value problems for guidance is that there is a loss of accuracy as compared with solving the full nonlinear problem. For systems engineering purposes, it is important to understand and bound the error in guidance calculations that can be encountered under expected mission circumstances. We show here that the norm of tensors related to the order $m+1$ state transition tensor can be used to approximately bound the error associated with using an order $m$ state transition tensor model to perform one of the basic guidance and control tasks described above. Though not outlined here, state transition tensors may be used to describe evolution of the states and costates associated with the optimal control of a satellite in the vicinity of some orbit \cite{kulik2023state}. Approximation error for the thrust, control costs, and terminal states under that policy can be computed in a very similar manner to the impulsive controls presented here. The use of a differential algebra system such as DACE to compute the coefficients of forward and inverse series associated with the flow of a dynamical system makes complex operations to arbitrary order possible while eliminating the need to hand code the complicated tensor expressions described above \cite{rasotto2016differential}. Restructuring these coefficients of the polynomials into a tensor form would then allow one to perform tensor norm/eigenvalue analysis though we do not implement these calculations here.

\subsubsection{Error Bounds for Propagation of Relative Motion}
\begin{figure}
    \centering
    \includegraphics[width=.5\textwidth]{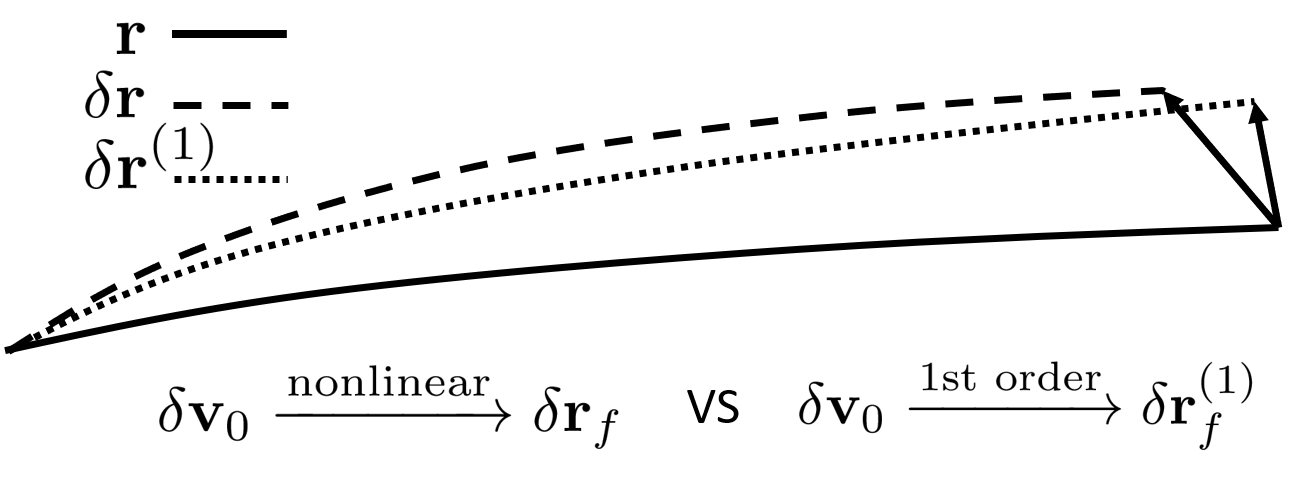}
    \caption{Notional depiction of a reference orbit, nonlinear propagation of relative position, and linear propagation of relative position.}
    \label{fig:prop_notional}
\end{figure}

The $m$-th order Taylor series approximation of relative motion of a satellite around some reference orbit is given by
\begin{equation}
	 \delta \mathbf{x}_f \approx \sum_{l=1}^m \frac{1}{l!}\mathbf{\Psi}^{(l)}(t_f,t_0)\delta \mathbf{x}_0^l
 %+ \mathbf{\Phi}^{(1)}(t_f,t_0) \delta \mathbf{x} + \frac{1}{2} \mathbf{\Phi}^{(2)}(t_f,t_0)\delta \mathbf{x}^2
\end{equation}
This $m$-th order approximation of relative motion has error
\begin{equation}
    \Vert \delta\mathbf{x}_f-\delta\mathbf{x}_f^{(m)}\Vert_2\leq \frac{1}{(m+1)!}\Vert \mathbf{\Psi}^{(m+1)}\Vert_2 \Vert \delta\mathbf{x}_0\Vert_2^{m+1}+\mathcal{O}(\Vert \delta\mathbf{x}_0\Vert_2^{m+2})
\end{equation}
where the notation $\cdot^{(m)}$ denotes the $m$th order Taylor approximation of the vector in question. The expression above makes sense to consider when the position and velocity are non-dimensionalized, allowing for them to be similarly scaled. In the case where position and velocity are not on similar scales, we can also consider errors in position and velocity separately by slicing the relevant tensors. For example, a first-order model position error based on some initial velocity alone can be described:
\begin{equation}
    \Vert \delta\mathbf{r}_f-\delta\mathbf{r}_f^{(1)}\Vert_2\leq \frac{1}{2}\Vert \mathbf{\Psi}^\mathbf{r}_\mathbf{vv}\Vert_2 \Vert \delta\mathbf{v}_0\Vert_2^{2}+\mathcal{O}(\Vert \delta\mathbf{v}_0\Vert_2^{3})
    \label{eqn:lin_prop_error}
\end{equation}
where $\mathbf{\Psi}^\mathbf{r}_{\mathbf{v}\mathbf{v}}(t_f, t_0)$ denotes a block of the second-order state transition tensor--the second-order partial derivative tensor for the final position with respect to the initial velocity:
\begin{equation}
    \mathbf{\Psi}^\mathbf{r}_{\mathbf{v}\mathbf{v}}(t_f, t_0)=\frac{\partial^2\mathbf{r}_f}{\partial\mathbf{v}_0^2}
\end{equation}
The inequality in Eq. \ref{eqn:lin_prop_error} can be employed to approximately bound the error between the linearly propagated and nonlinearly propagated position as depicted in Fig. \ref{fig:prop_notional}. To clarify the phrase ``approximately bound," we mean that this is a true upper bound on the second-order contribution to the error, which is dominant for small $\delta\mathbf{v}_0$ making this a good approximation of the true upper bound on the full nonlinear error when $\delta\mathbf{v}_0$ is sufficiently small. One potential heuristic to check whether $\delta\mathbf{v}_0$ is sufficiently small for the second-order term to well-approximate the full nonlinear error is to find the norm of the third-order term and check that this is small relative to the second-order term at the $\delta\mathbf{v}_0$ scale of interest. The second-order term is also likely to be a good approximation of the overall error when the second-order error bound is much smaller than the linear quantity under study at the $\delta\mathbf{v}_0$ scale of interest.
%\FloatBarrier
We present two examples utilizing the above error bound formulation. In both, we find the maximum error between nonlinear and linear models of the position $\Vert \delta\mathbf{r}_f-\delta\mathbf{r}_f^{(1)}\Vert_2$ given initial position at the reference orbit, and initial velocity perturbed within some ball of radius $R$ from the reference orbit initial velocity. This situation where the initial position is at the reference orbit, but the initial velocity perturbation is nonzero models the effects of an impulsive burn on a single satellite or the velocity imparted by an explosion on pieces of debris from a satellite. We wish to answer the question: how well does a linear model of the dynamics capture these behaviors for different scales of the velocity perturbation? To do so, analysis is conducted as the maximum scale $R$ of the velocity perturbation $\delta\mathbf{v}_0$ is varied from 0 to 200 meters per second. In the first case, two-body dynamics are employed with an International Space Station-like reference orbit having classical orbital elements semi-major axis, eccentricity, right-ascension of the ascending node, argument of perigee, and initial mean anomaly respectively are $(a,e,i,\Omega,\omega,M_0)=(6738 \text{km}, 0.000514, 51.6434^{\degree}, 0^{\degree}, 0^{\degree}, 0^{\degree})$. %Note that we set the last three elements to $0$ as these orbital elements change over time due to the $J_2$-perturbation as well as regular two-body motion, and make no impact on our analysis. 
In the second case, we study the restricted three-body dynamics with the initial conditions of the near-rectilinear halo orbit (NRHO) proposed for NASA's Gateway. We propagate each satellite one tenth of it's respective period in the below error approximations.

In order to validate the approximation in Eq. \ref{eqn:lin_prop_error}, we compare the results from that tensor norm bound with three other methods of characterizing the linearization error. First, to demonstrate a simple method to improve the estimate from Eq. \ref{eqn:lin_prop_error}, we make an educated guess of the perturbation to cause the maximum error between the linear and full nonlinear models. During the calculation of the tensor 2-norm, we find the maximal eigenvector arising in the calculation of the tensor 2-norm. This eigenvector is the unit-normed initial velocity perturbation that produces the largest error between the linear and quadratic models of the dynamics. We then take this unit-normed initial velocity perturbation, scale it by $R$, and propagate it under the full nonlinear dynamics and linear dynamics to find the final position error. We perform this calculation for the eigenvector calculated as well as the eigenvector in the opposite direction which also produces the same maximum error between linear and quadratic dynamics. The maximum of these two evaluations is taken. This tends to improve the approximation of the maximum error between the linear and full nonlinear models, with only two more evaluations of the nonlinear dynamics propagation. 

Next, to determine the true maximum error between linear and nonlinear models, we employ a local numerical optimization using the sequential least squares programming (SLSQP) method in SciPy \cite{2020SciPy-NMeth} on the negative of the objective function for maximization. Our objective function in this maximization is the comparison of the propagated true and linear dynamics as a function of the perturbation $\delta\mathbf{v}_0$. The constraint of this maximization is that the velocity perturbation must lie on the sphere defined by the radius $R$. The initial guess is the eigenvector with the sign that leads to the maximum error between linear and nonlinear evaluations. This method is the best approximation of the true upper bound of the full nonlinear error on the surface of the sphere under the assumption that the local optimization is finding the global maximum.

Finally, a sampling method is used to assure that a global maximum in the error is being approximated in the locally optimized ``true" maximum error above. We sample an initial velocity perturbation uniformly at random from the sphere of radius $R$. This perturbation is then propagated under the full nonlinear dynamics and linear dynamics to find the final position error between the two models. We repeat this sampling 5,000 times, reporting the norm of the largest error. This estimate will not be very precise, but has a good chance of finding a point near the global maximum rather than a local maximum. This process gives the maximum error via sampling, and we check that our maximum error from sampling is always below the ``true" maximum error we calculate above using local optimization to validate that we did not locally optimize to a point other than the global maximum.

\begin{figure}[hbt!]
     \centering
     \begin{subfigure}[b]{0.45\textwidth}
         \centering
         \includegraphics[width=\textwidth]{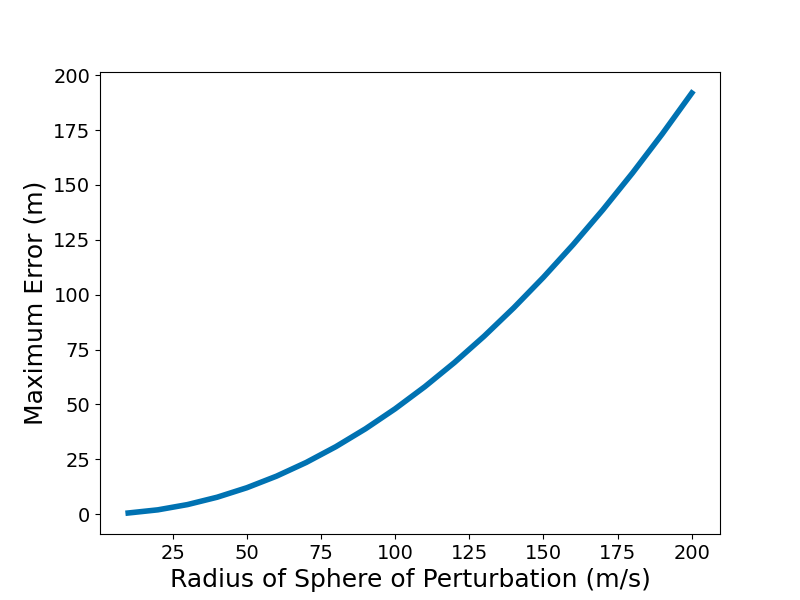}
         \caption{Perturbation vs. Maximum error in ISS example.}
     \end{subfigure}
     \hfill
     \begin{subfigure}[b]{0.45\textwidth}
         \centering
         \includegraphics[width=\textwidth]{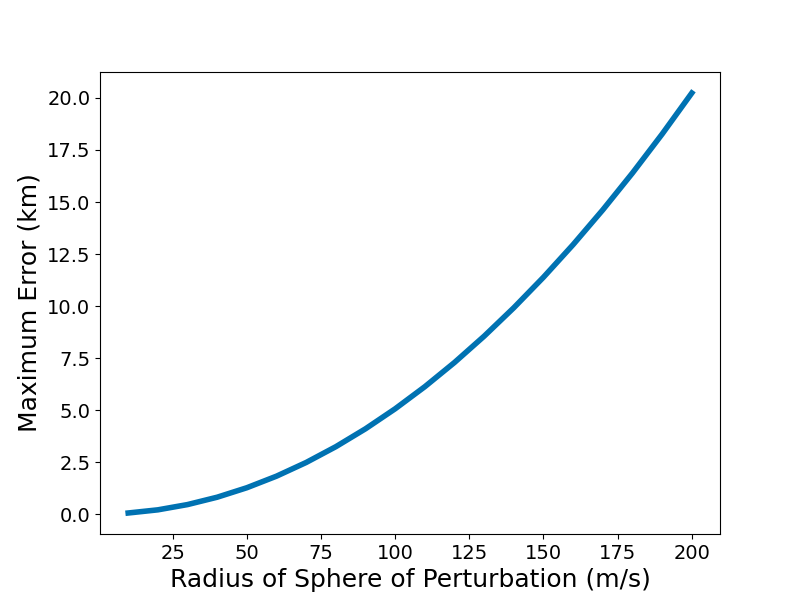}
         \caption{Perturbation vs. Maximum error in NRHO example.}
     \end{subfigure}\\
        \caption{Maximum error for propagation of relative motion in two-body and restricted three-body dynamics through SciPy optimization. }
        \label{fig:prop_max_error_plots}
\end{figure}

\begin{figure}[hbt!]
     \centering
     \begin{subfigure}[b]{0.45\textwidth}
         \centering
         \includegraphics[width=\textwidth]{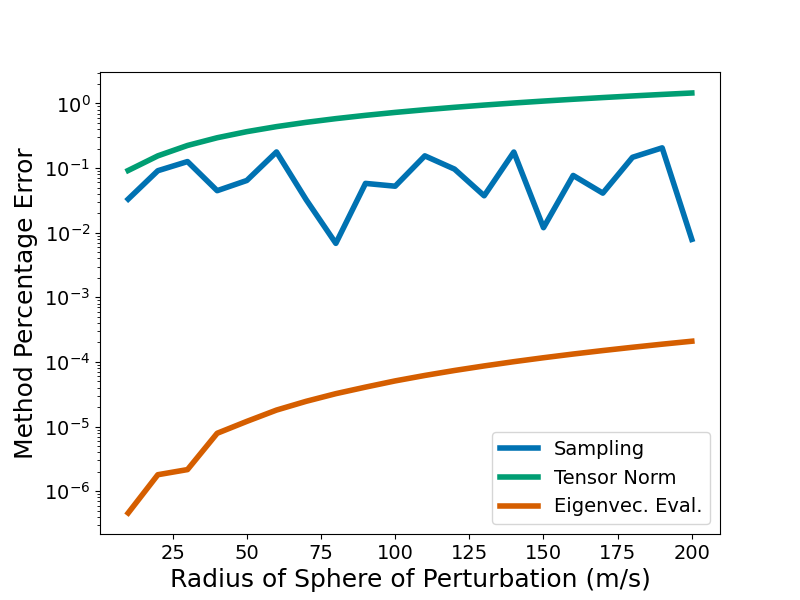}
         \caption{Percent error between methods in ISS example.}
     \end{subfigure}
     \hfill
     \begin{subfigure}[b]{0.45\textwidth}
         \centering
         \includegraphics[width=\textwidth]{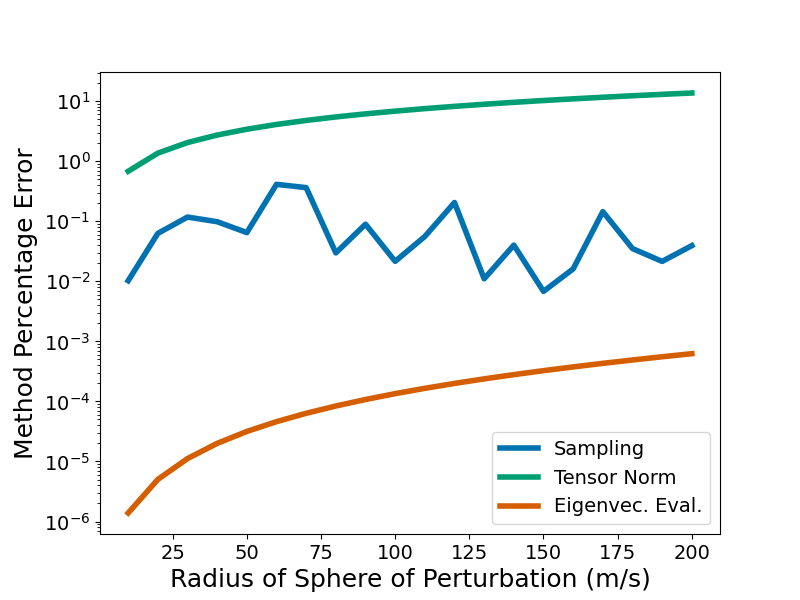}
         \caption{Percent error between methods in NRHO example.}
     \end{subfigure}\\
        \caption{Relative inaccuracy between methods of calculating maximum error and true upper bound for propagation of relative motion in two-body and restricted three-body dynamics. }
        \label{fig:prop_perc_error_plots}
\end{figure}

\begin{figure}[hbt!]
      \centering
      \begin{subfigure}[b]{0.45\textwidth}
          \centering
          \includegraphics[width=\textwidth]{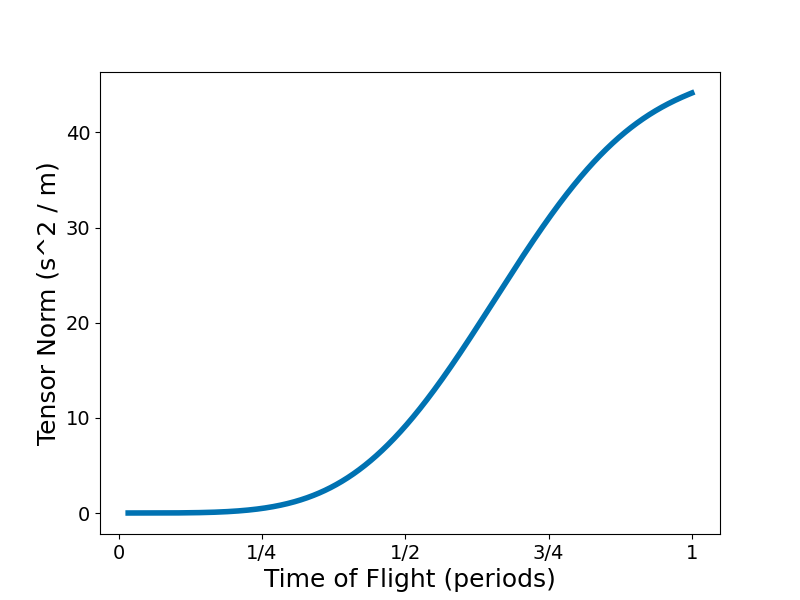}
          \caption{Propagation tensor norm vs. Time of flight in ISS example.}
      \end{subfigure}
      \hfill
      \begin{subfigure}[b]{0.45\textwidth}
          \centering
          \includegraphics[width=\textwidth]{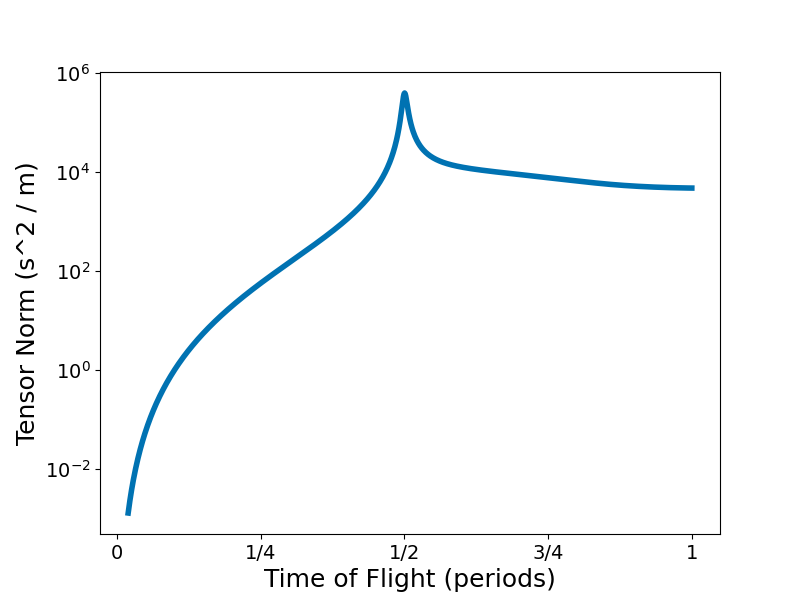}
          \caption{Propagation tensor norm vs. Time of flight in NRHO example (log).}
      \end{subfigure}\\
         \caption{Plots of velocity to position propagation tensor norm as time-of-flight varies.}
         \label{fig:prop_tnorm_plots}
\end{figure}

In Fig. \ref{fig:prop_max_error_plots}, we present the maximum errors for propagation of relative motion through the SLSQP method in SciPy. Results are presented in both the two- and restricted three-body cases as discussed above. In the two-body case, we see the magnitude of the maximum error reaches on the order of hundreds of meters with a quadratic relationship between position error and the perturbation magnitude. In the restricted three-body case, we see the magnitude of the maximum error is on the order of tens of kilometers, again with a quadratic relationship. We expect a larger error in the three-body case because this dynamical system is more sensitive to variations in initial conditions.
%\FloatBarrier
%\FloatBarrier
In Fig. \ref{fig:prop_perc_error_plots}, we present the percentage error between numerical optimization and the three other approaches for estimating the error: sampling, tensor norm (Eq. \ref{eqn:lin_prop_error}), and nonlinear evaluation of the tensor eigenvector. Results are presented in both the two- and restricted three-body cases discussed above. We see the error found by sampling to be on the scale of 1\%, which implies that our local SLSQP optimization is, in fact, finding the global maximum. To answer our proposed question regarding the accuracy of a linear model, in the two-body case, the difference between the tensor norm linear model and true upper error bound increases from $0.1\%$ to $10\%$ over the range of velocity perturbations. In the three body case, this error increases from $~3\%$ to $10\%$. Further, we find the error in our eigenvector evaluation to be on the scale of one one thousandth of one percent. This implies that the direction that maximizes the nonlinear error is very close to the direction which maximizes the quadratic error even though the quadratic and nonlinear errors may deviate by up to one tenth of the true nonlinear error.
%\FloatBarrier
In Fig. \ref{fig:prop_tnorm_plots}, we examine the tensor norm defined in Eq. (\ref{eqn:lin_prop_error}). Using the same ISS and Gateway NRHO examples, we instead propagate each satellite up to a single period. The tensor 2-norm associated with this flight is plotted against the time of flight, note that the NRHO example has a logarithmic y-scale. We see the tensor norm in the ISS example exhibits sub-exponential growth, while the tensor norm in the NRHO example exhibits exponential growth. The differences in growth can be attributed to each system's sensitivity to initial conditions and level of nonlinearity. %We quantify this value using the metric of the nonlinearity index, which can be calculated using tensor norms. This concept is explored in depth within section C. and D. of this paper.
We also can see a spikes in the growth of the tensor norm associated with the NRHO example at around half a period. This spike corresponds to the perilune of the NRHO orbit where the final position is particularly sensitive to initial conditions. The value of the tensor norm is not approaching infinity at this time, but reaching a high finite value. %These points are known to have high levels of uncertainty within the realm of linear dynamical approximations, and are referred to as relative transfer singularities. Explicitly, these are points where the STM slice $\mathbf{\Phi}^\mathbf{r}_\mathbf{v}$ is singular.
%\FloatBarrier
%\begin{equation}
%    \max_{\Vert\mathbf{x}\Vert_2\leq R}\Vert\mathbf{B}\mathbf{x}^m\Vert_2 = \Vert \mathbf{B}\Vert_2 R^m
%\end{equation}
%\subsubsection{Error Bounds for Impulsive Relative Motion Control}
\subsubsection{Final Position Error from an Impulsive Transfer}
\label{sec:miss-distance}

The impulsive relative transfer problem is a basic building block for station-keeping as well as rendezvous. In its simplest form, a satellite at the reference orbit incurs some $\Delta \mathbf{v}$ at the time $t_0$ in order to reach a desired final relative position $\delta\mathbf{r}_f^*$ with respect to the reference orbit at some time $t_f$. A linear approximation is often used to solve for the required initial relative velocity \cite{howell2005natural}
\begin{equation}
    \delta\mathbf{v}_0^{(1)}=(\mathbf{\Phi}^\mathbf{r}_\mathbf{v}(t_f, t_0))^{-1} \delta\mathbf{r}^*_f
\end{equation}
where $\mathbf{\Phi}^\mathbf{r}_\mathbf{v}$ is the upper right three by three block of $\Phi$, and the superscript $(1)$ indicates that the solution for the initial velocity is a first-order approximation. We might want to understand the accuracy of the linear approximation of the initial velocity $\delta\mathbf{v}_0^{(1)}$. A natural question is what level of error does implementing the linear guidance solution lead to in the final relative position achieved by the $\delta\mathbf{v}_0^{(1)}$ versus the desired $\delta\mathbf{r}_f^*$ when propagated by the true nonlinear dynamics? We are trying to find the maximum miss distance when implementing a linear guidance solution in a truly nonlinear dynamical system. To answer this question, we begin by substituting the linear solution for the initial impulse into a Taylor expansion for the final position
\begin{equation}
    \delta\mathbf{r}_f = \mathbf{\Phi}^\mathbf{r}_\mathbf{v}(t_f, t_0)\delta\mathbf{v}_0^{(1)} + \frac{1}{2}\mathbf{\Psi}^\mathbf{r}_{\mathbf{v}\mathbf{v}}(t_f, t_0)\left(\delta\mathbf{v}_0^{(1)}\right)^2+\mathcal{O}\left(\left(\delta\mathbf{v}_0^{(1)}\right)^3\right)
\end{equation}
If we define the position error tensor associated with the second-order order term from the above Taylor expansion
\begin{equation}
    (\mathbf{E}^{(1)})^i_{j,k}=\frac{1}{2}\left(\mathbf{\Psi}^\mathbf{r}_{\mathbf{v}\mathbf{v}}(t_f, t_0)(\mathbf{\Phi}^\mathbf{r}_\mathbf{v}(t_f, t_0))^{-2}\right)^i_{j,k}=\frac{1}{2}\frac{\partial^2 r_f^i}{\partial v_0^l\partial v_0^p}\left(\left(\frac{\partial\mathbf{r}_f}{\partial\mathbf{v}_0}\right)^{-1}\right)^l_j\left(\left(\frac{\partial\mathbf{r}_f}{\partial\mathbf{v}_0}\right)^{-1}\right)^p_k
    \label{eq:E1tens}
\end{equation}
then, when a linear approximation is used to determine the initial $\Delta \mathbf{v}$ for a transfer, the error in the final position is

    %\Vert\delta\mathbf{r}_f^*-\delta\mathbf{r}_f\Vert_2\leq  \left\Vert \frac{1}{2}\mathbf{\Psi}^\mathbf{r}_{\mathbf{v}\mathbf{v}}(t_f, t_0)(\mathbf{\Phi}^\mathbf{r}_\mathbf{v}(t_f, t_0))^{-2}\right\Vert_2 \Vert \delta\mathbf{r}_f^*\Vert_2^2 + \mathcal{O}\left(\Vert \delta\mathbf{r}_f^*\Vert_2^3\right)
    \begin{equation}
    \Vert\delta\mathbf{r}_f^*-\delta\mathbf{r}_f\Vert_2\leq  \left\Vert\mathbf{E}^{(1)}\right\Vert_2 \Vert \delta\mathbf{r}_f^*\Vert_2^2 + \mathcal{O}\left(\Vert \delta\mathbf{r}_f^*\Vert_2^3\right)
\end{equation}
where the 2-norm of the (1-2) tensor above can be calculated by the Z-eigenvector approach as described in Sec. \ref{sec:2-norm}.

Multivariate series reversion methods can be used to find higher-order approximations to inverse problems such as relative motion impulsive transfers and rendezvous \cite{griffith2004higher, turneraas, junkins2009generalizations}. For example, the second-order approximation of the initial relative velocity for the simple case of a transfer is
\begin{equation}
    \delta\mathbf{v}_0^{(2)}=(\mathbf{\Phi}^\mathbf{r}_\mathbf{v}(t_f, t_0))^{-1} \delta\mathbf{r}^*_f-\frac{1}{2}(\mathbf{\Phi}^\mathbf{r}_\mathbf{v}(t_f, t_0))^{-1} \mathbf{\Psi}^\mathbf{r}_\mathbf{vv}(t_f, t_0) (\mathbf{\Phi}^\mathbf{r}_\mathbf{v}(t_f, t_0))^{-2} \left(\delta\mathbf{r}^*_f\right)^2
\end{equation}
The tensor associated with the position error from the second-order initial velocity solution is 
\begin{equation}
    (\mathbf{E}^{(2)})^i_{j,k,l}=%((\mathbf{\Phi}^\mathbf{r}_\mathbf{v})^{-1})^i_{p_1}
    \left(-\left(\mathbf{\Psi}^\mathbf{r}_\mathbf{vv}\right)^{i}_{p_1,p_3}((\mathbf{\Phi}^\mathbf{r}_\mathbf{v})^{-1})^{p_1}_{p_2}\left(\mathbf{\Psi}^\mathbf{r}_\mathbf{vv}\right)^{p_2}_{p_4,p_5}+\frac{1}{6}\left(\mathbf{\Psi}^\mathbf{r}_\mathbf{vvv}\right)^{p_1}_{p_3,p_4,p_5}\right)((\mathbf{\Phi}^\mathbf{r}_\mathbf{v})^{-1})^{p_3}_{j}((\mathbf{\Phi}^\mathbf{r}_\mathbf{v})^{-1})^{p_4}_{k}((\mathbf{\Phi}^\mathbf{r}_\mathbf{v})^{-1})^{p_5}_{l}
\end{equation}
and the error in final position from a second-order approximation of the initial velocity is
\begin{equation}
    \Vert\delta \mathbf{r}_f^*-\delta\mathbf{r}_f\Vert_2\leq  \left\Vert \mathbf{E}^{(2)}\right\Vert_2 \Vert \delta\mathbf{r}_f^*\Vert_2^3 + \mathcal{O}\left(\Vert \delta\mathbf{r}_f^*\Vert_2^4\right)
\end{equation}
This process can be generalized to find the position error from using an $m$th-order approximation of the correct initial $\Delta \mathbf{v}$, but we only specify the exact form of the error tensor up to second-order due to the complexity of higher-order terms.
\begin{equation}
    \Vert\delta \mathbf{r}_f^*-\delta\mathbf{r}_f\Vert_2\leq  \left\Vert \mathbf{E}^{(m)}\right\Vert_2 \Vert \delta\mathbf{r}_f^*\Vert_2^{m+1} + \mathcal{O}\left(\Vert \delta\mathbf{r}_f^*\Vert_2^{m+2}\right)
\end{equation}
\begin{figure}
    \centering
    \includegraphics[width=.5\textwidth]{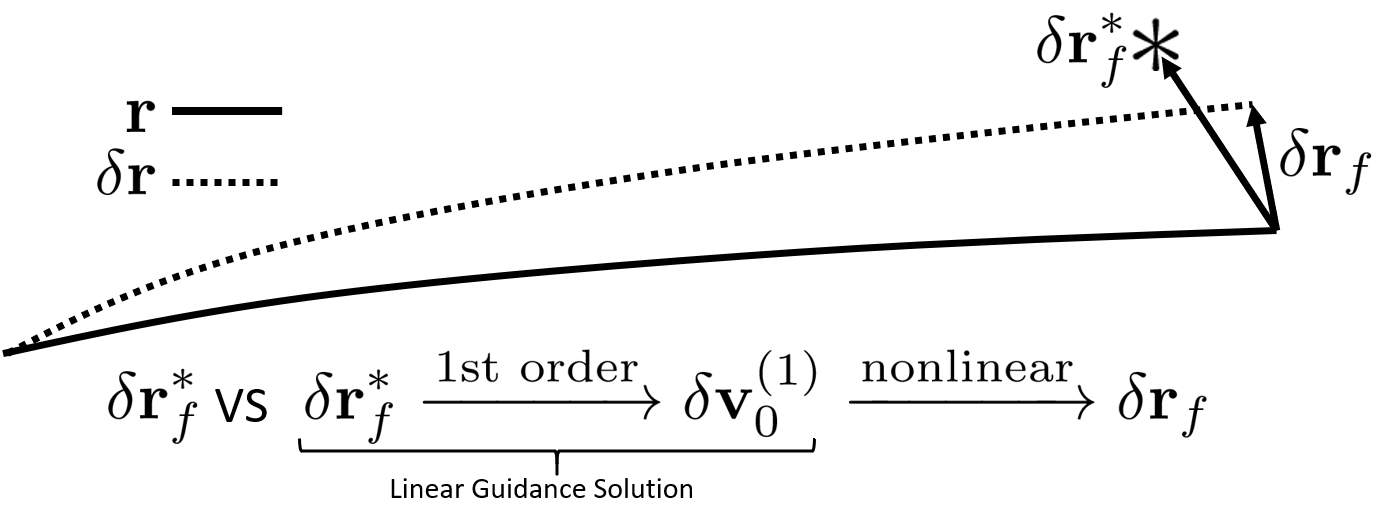}
    \caption{Notional depiction of transfer miss distance for linearized initial velocity calculation to reach $\delta\mathbf{r}_f^*$.}
    \label{fig:pos_notional}
\end{figure}
A notional depiction comparing the linearized guidance trajectory and showing the miss distance from the desired final is presented in Fig. \ref{fig:pos_notional}. We again present two examples utilizing the above error bound formulation. In both, we find the maximum error between nonlinear and linear models of the position $\|\delta \mathbf{r}_f^*-\delta\mathbf{r}_f\Vert_2$ given an initial position at the reference orbit, and some desired final position sampled from within a ball of radius $R$ from the reference orbit final position. %This situation in which a satellite attempts to reach some desired final position is a common problem throughout astrodynamics. 
We use a linear approximation to find the impulsive burn at $t_0$ that results in the satellite reaching this desired final position $\delta\mathbf{r}_f^*$. To find the maximum miss distance associated with implementing the linear guidance solution, analysis is conducted as the maximum scale $R$ of the final position $\delta\mathbf{r}_f^*$ is varied from 0 to 200 kilometers in the two-body orbit case and between 0 and 2,000 kilometers in the NRHO example. We use the same example cases as in the previous application. Further, we use the same methods of calculating the maximum error, with the adjustment of sampling $R$ as the magnitude of the final desired position deviation.

\begin{figure}[hbt!]
     \centering
     \begin{subfigure}[b]{0.45\textwidth}
         \centering
         \includegraphics[width=\textwidth]{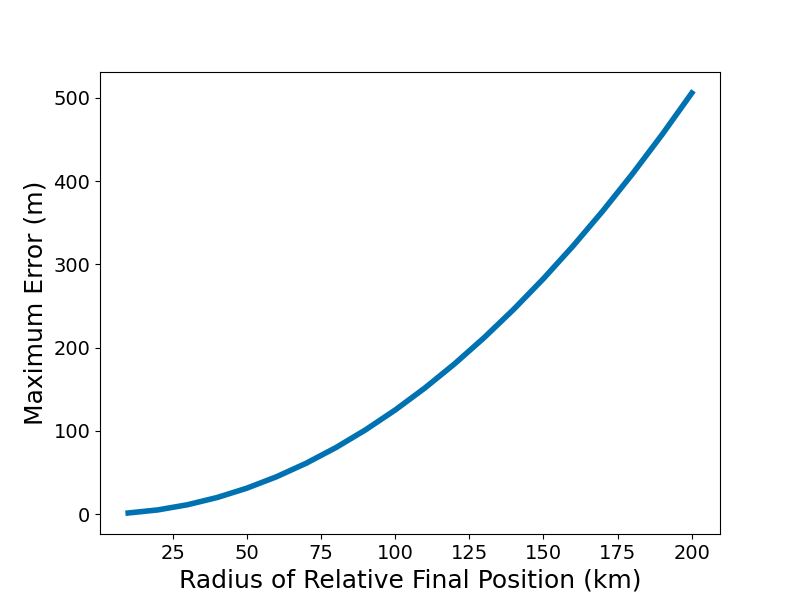}
         \caption{Final position vs. Maximum error in ISS example.}
     \end{subfigure}
     \hfill
     \begin{subfigure}[b]{0.45\textwidth}
         \centering
         \includegraphics[width=\textwidth]{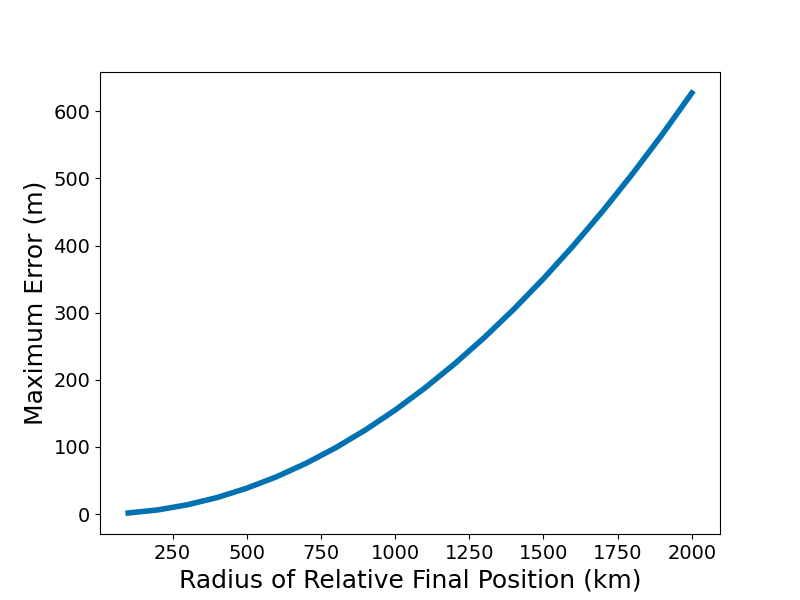}
         \caption{Final position vs. Maximum error in NRHO example.}
     \end{subfigure}\\
        \caption{Maximum error for propagation of an impulsive transfer in two-body and restricted three-body dynamics through SciPy optimization. }
        \label{fig:guide_max_error_plots}
\end{figure}

\begin{figure}[hbt!]
     \centering
     \begin{subfigure}[b]{0.45\textwidth}
         \centering
         \includegraphics[width=\textwidth]{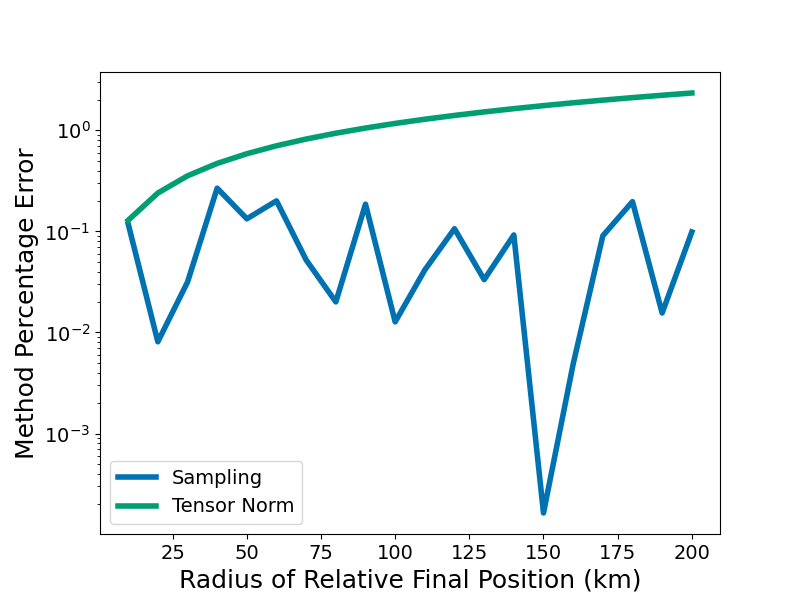}
         \caption{Percent error between methods in ISS example.}
     \end{subfigure}
     \hfill
     \begin{subfigure}[b]{0.45\textwidth}
         \centering
         \includegraphics[width=\textwidth]{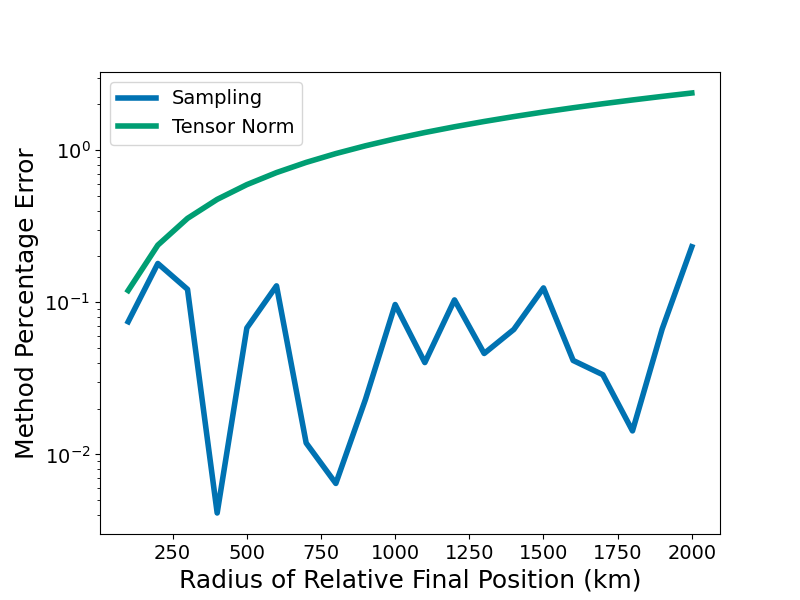}
         \caption{Percent error between methods in NRHO example.}
     \end{subfigure}\\
        \caption{Relative inaccuracy between methods of calculating maximum error and true upper bound for miss distance in an impulsive transfer in two-body and restricted three-body dynamics.}
        \label{fig:guide_perc_error_plots}
\end{figure}

In Fig. \ref{fig:guide_max_error_plots} we present the maximum errors for propagation of an impulsive transfer through the SLSQP method in SciPy. Similarly to the previous application, the results are presented in both the two- and restricted three-body cases. In the two-body case, we see the magnitude of the maximum error reaches on the order of hundreds of meters. For the order of magnitude larger distance traveled in the NRHO case, the magnitude of error is still on the order of hundreds of meters. This difference in error can be attributed to the restricted three-body system's sensitivity to initial conditions. In the opposite manner of propagation, this results in reduced error, as a given final position corresponds to a smaller initial velocity to reach the desired final position. Thus, a linear approximation performs better in terms of error for solving a transfer in the NRHO case when compared to the two-body case. 
%\FloatBarrier
In Fig. \ref{fig:guide_perc_error_plots}, we present the percentage error between numerical optimization and the three other approaches for estimating error. Although more erratic, we again find the sampling method-based error on the scale of tenths of a percent to just under ten percent of the true maximum nonlinear error. This again implies that our local SLSQP optimization is finding the global minimum. To answer our proposed question regarding the accuracy of this first-order linear model, we find the difference between the tensor
norm linear model and true upper error bound increases from 0.1\% to 10\% over the range of final positions. In
the three body case, this error increases from ~0.1\% to ~1\%. The sampling method giving around 0.1\% error relative to the optimized linearization error indicates that the global maximum is being found correctly by the local optimization method in this guidance problem as well. Further, we find the absolute relative error in our eigenvector evaluation (not shown) to be on the scale of less than $10^{-7}$ in both cases indicating very high accuracy in the agreement of the direction of maximal second-order error and nonlinear error. 

%\FloatBarrier
\begin{figure}[hbt!]
     \centering
     \begin{subfigure}[b]{0.45\textwidth}
         \centering
         \includegraphics[width=\textwidth]{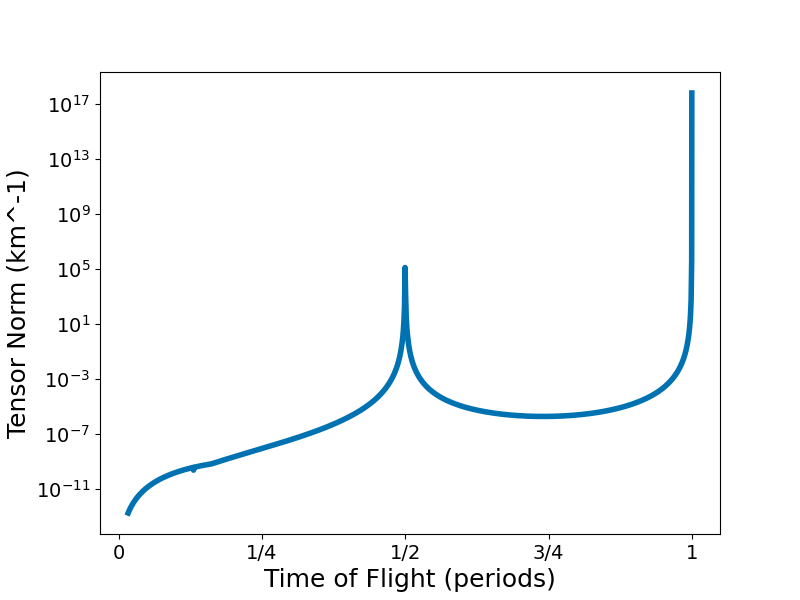}
         \caption{Miss distance tensor norm vs. Time of flight in ISS example.}
     \end{subfigure}
     \hfill
     \begin{subfigure}[b]{0.45\textwidth}
         \centering
         \includegraphics[width=\textwidth]{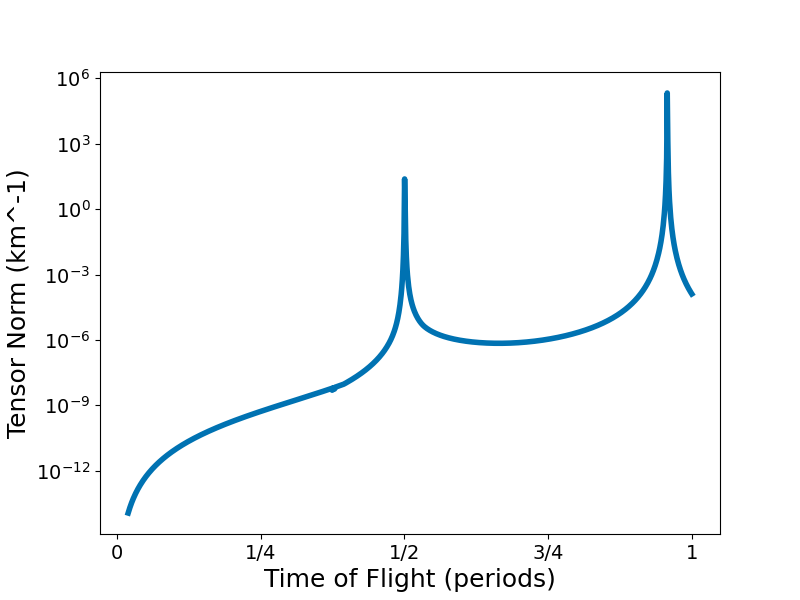}
         \caption{Miss distance tensor norm vs. Time of flight in NRHO example.}
     \end{subfigure}\\
        \caption{Plots of norm associated with miss distance calculation for linearized guidance over varying transfer durations.}
        \label{fig:guide_tnorm_plots}
\end{figure}

%\FloatBarrier
In Fig. \ref{fig:guide_tnorm_plots}, we examine the norm of the tensor defined in Eq. \ref{eq:E1tens}. We perform the same analysis as the previous application, propagating each satellite one period. In this case, both the ISS and NRHO examples are plotted using a logarithmic y-scale. The norm of the tensor spikes at around half and full periods. These represent relative transfer singularities and are times when the block of the STM $\mathbf{\Phi}^\mathbf{r}_\mathbf{v}$ becomes singular \cite{kulik2022relative,mullins1992initial,fitzgerald1998pinch}. At these points the tensor norm is approaching infinity unlike in the earlier previous propagation example where a high, but finite value is achieved. The same approach to infinity is true of spikes in subsequent guidance-related examples.
%\FloatBarrier

\subsubsection{Initial Velocity Error from an Impulsive Transfer}
\label{sec:vel-error}
\begin{figure}[ht]
    \centering
    \includegraphics[width=.5\textwidth]{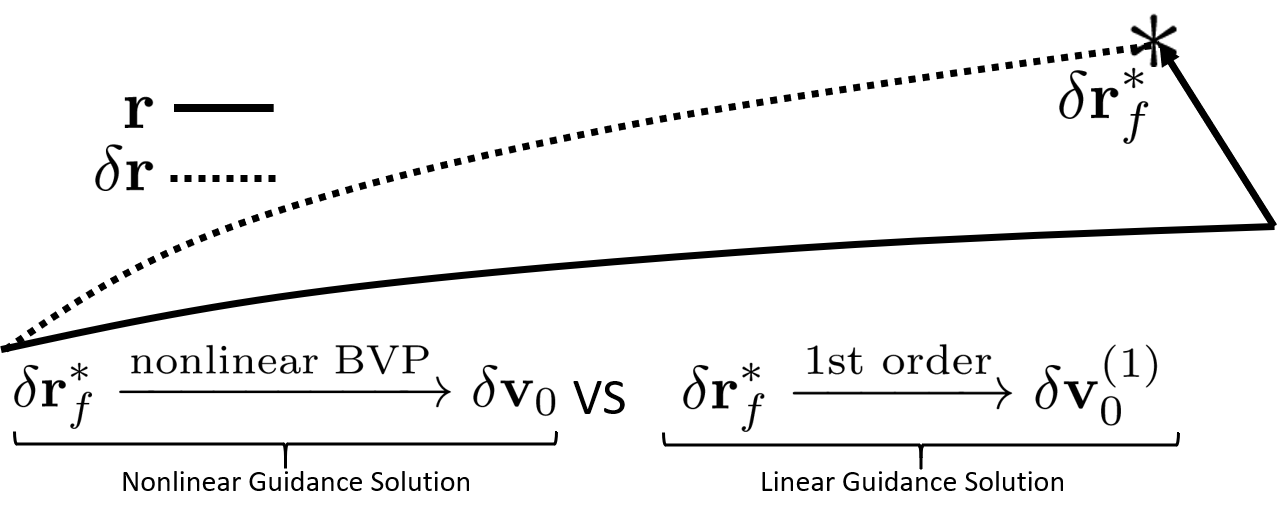}
    \caption{Notional depiction of error in linearized initial velocity calculation to reach $\delta\mathbf{r}_f^*$.}
    \label{fig:vel_notional}
\end{figure}

In Sec. \ref{sec:miss-distance}, we calculated bounds for the miss distance of a linear or higher-order series reversion-based transfer as a function of the desired distance of the transfer. Here we calculate bounds on the inaccuracy of the initial velocity calculated to reach the final target. The correction term to improve the approximation of the initial velocity of a transfer to one order higher is given in terms of the error tensor from the previous section, so that the $(m+1)$-th order solution is
\begin{equation}
    \delta \mathbf{v}_0^{(m+1)}=\delta \mathbf{v}_0^{(m)}-(\mathbf{\Phi}^\mathbf{r}_\mathbf{v})^{-1}\mathbf{E}^{(m)}(\delta\mathbf{r}_f^*)^{m+1}
    \label{eq:series-newton}
\end{equation}
where
\begin{equation}
    \left((\mathbf{\Phi}^\mathbf{r}_\mathbf{v})^{-1}\mathbf{E}^{(m)}\right)^i_{j_1...j_m}=\left((\mathbf{\Phi}^\mathbf{r}_\mathbf{v})^{-1}\right)^i_l\left(\mathbf{E}^{(m)}\right)^l_{j_1...j_m}
    \label{eq:E-tens-vel}
\end{equation}
Eq.~\ref{eq:series-newton} can be viewed as Newton's method applied to gain one more order of accuracy in the series reversion per iteration \cite{griffith2004higher, lizia2008application, zhang2010second}. %Though there exist more efficient methods for performing multivariate series reversion, 
The resulting final term in the inverse series for initial velocity in terms of final position affords the following observation about the error in the $m$-th order solution to the relative transfer problem: %Thus, the error in the $m$th order approximation of the initial velocity using series reversion for an impulsive transfer is
\begin{equation}
    \Vert\delta \mathbf{v}_0^{*}-\delta \mathbf{v}_0^{(m)}\Vert_2 \leq \Vert(\mathbf{\Phi}^\mathbf{r}_\mathbf{v})^{-1}\mathbf{E}^{(m)}\Vert_2\Vert\delta\mathbf{r}_f^*\Vert^{m+1}+\mathcal{O}(\Vert\delta\mathbf{r}_f^*\Vert^{m+2})
    \label{eq:bound-vel}
\end{equation}
Thus, we can approximately bound the error in the order $m$ series reversion solution for the initial velocity to achieve some desired transfer. In particular, this yields an approximate error bound for the initial velocity computed by linear guidance methods versus the true initial velocity that solves the full nonlinear boundary value problem. A notional depiction of the transfer trajectory is presented in Fig. \ref{fig:vel_notional}. The boundary conditions are the same for both the linear and nonlinear versions of the guidance problem, and in each self-consistent model of the dynamics, the solution leads to the final desired position, though the two initial velocity solutions are clearly different (and the linearized solution does not reach the desired position under the nonlinear dynamics as discussed above).
%%%change from here on

\begin{figure}[hbt!]
     \centering
     \begin{subfigure}[b]{0.45\textwidth}
         \centering
         \includegraphics[width=\textwidth]{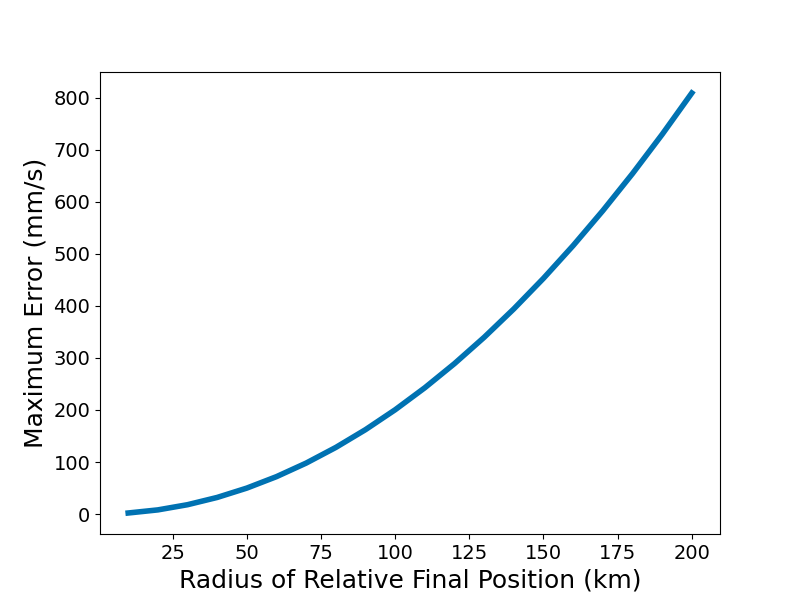}
         \caption{Final position vs. Maximum initial velocity error in ISS example.}
     \end{subfigure}
     \hfill
     \begin{subfigure}[b]{0.45\textwidth}
         \centering
         \includegraphics[width=\textwidth]{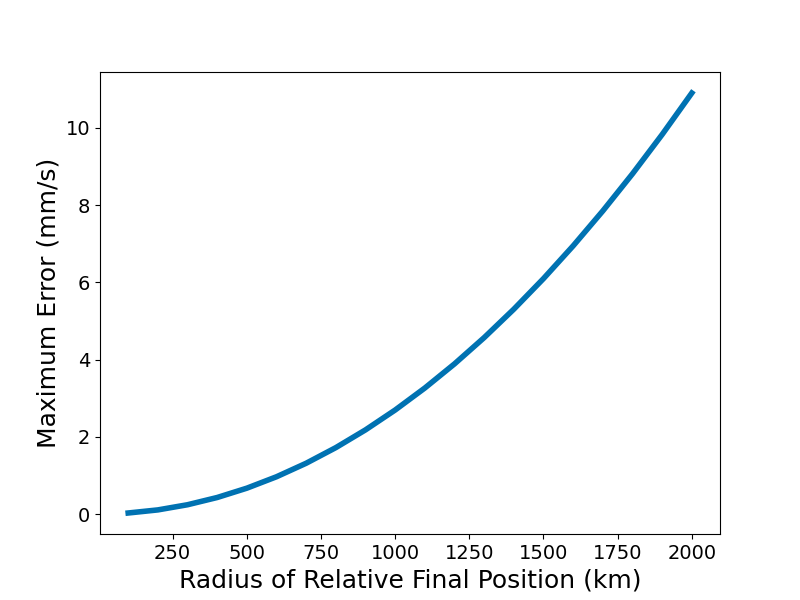}
         \caption{Final position vs. Maximum initial velocity error in NRHO example.}
     \end{subfigure}\\
        \caption{Maximum error for linear guidance calculation of initial velocity as a function of transfer distance.}
        \label{fig:vel_max_error_plots}
\end{figure}

\begin{figure}[hbt!]
     \centering
     \begin{subfigure}[b]{0.45\textwidth}
         \centering
         \includegraphics[width=\textwidth]{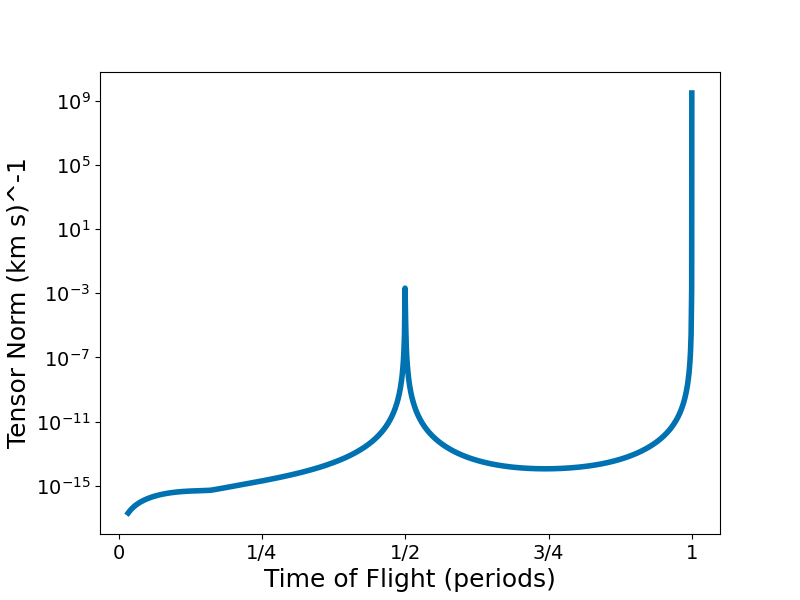}
         \caption{Tensor norm associated with velocity error vs. Time of flight in ISS example.}
     \end{subfigure}
     \hfill
     \begin{subfigure}[b]{0.45\textwidth}
         \centering
         \includegraphics[width=\textwidth]{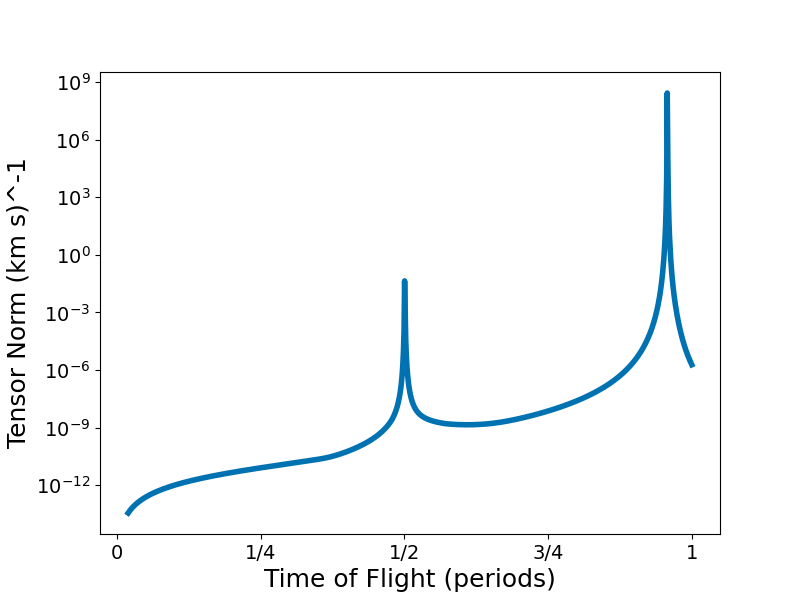}
         \caption{Tensor norm associated with velocity error vs. Time of flight in NRHO example.}
     \end{subfigure}\\
        \caption{Plots of the norm of the tensor that describes initial velocity error as a function of the desired transfer distance while varying time of flight.}
        \label{fig:vel_tnorm_plots}
\end{figure}

%\FloatBarrier
We present a set of examples under the same circumstances as in the miss-distance scenario for a linearized transfer. The true maximum errors for one tenth of a period transfer can be found in Fig. \ref{fig:vel_max_error_plots}. In the two-body case, hundred kilometer transfers calculated with linearized guidance can lead to nearly meter per second errors in the initial velocity, while thousand kilometer transfers in the three-body problem lead to errors on the order of centimeters per second. In some sense, this combined with the previous evidence on miss distance shows that while linear propagation has its limits, linear guidance can actually be quite accurate in the three-body problem. While this might be surprising from the perspective that the three-body problem is a chaotic dynamical system and control would seem more difficult, the large distance scales in this system come into play as well to make thousands of kilometers a relatively small distance as compared with the scale of the NRHO. The error of the tensor norm, sampling, and eigenvector evaluation based methods as compared relative to the numerical optimization around the eigenvector are on the same scale and exhibit the same trends as Fig. \ref{fig:guide_perc_error_plots} and are omitted.
%\FloatBarrier
% \begin{figure}[hbt!]
%      \centering
%      \begin{subfigure}[b]{0.45\textwidth}
%          \centering
%          \includegraphics[width=\textwidth]{Figures/Vel/twoBodyVelError.png}
%          \caption{Percent error between methods in ISS example.}
%      \end{subfigure}
%      \hfill
%      \begin{subfigure}[b]{0.45\textwidth}
%          \centering
%          \includegraphics[width=\textwidth]{Figures/Vel/threeBodyVelError.png}
%          \caption{Percent error between methods in NRHO example.}
%      \end{subfigure}\\
%         \caption{Relative inaccuracy between methods of calculating maximum error and true upper bound for propagation of an impulsive transfer in two-body and restricted three-body dynamics.}
%         \label{fig:vel_perc_error_plots}
% \end{figure}
%\FloatBarrier
In Fig. \ref{fig:vel_tnorm_plots}, we examine the norm of the tensor defined in Eq. \ref{eq:E-tens-vel} over the course of an orbital period. Similar spikes can be seen at relative transfer singularity times. Using this tensor norm from this plot and Eq. \ref{eq:bound-vel}, the initial velocity error from linearized guidance calculations can be bounded for small transfer distances. % Eq. \ref{eq:E1tens}. We perform the same analysis as the previous application, propagating each satellite two respective periods. In this case, both the ISS and NRHO examples are plotted using a logarithmic y-scale.
%\FloatBarrier
\subsubsection{Final Position Error During Rendezvous}
\begin{figure}
    \centering
    \includegraphics[width=.5\textwidth]{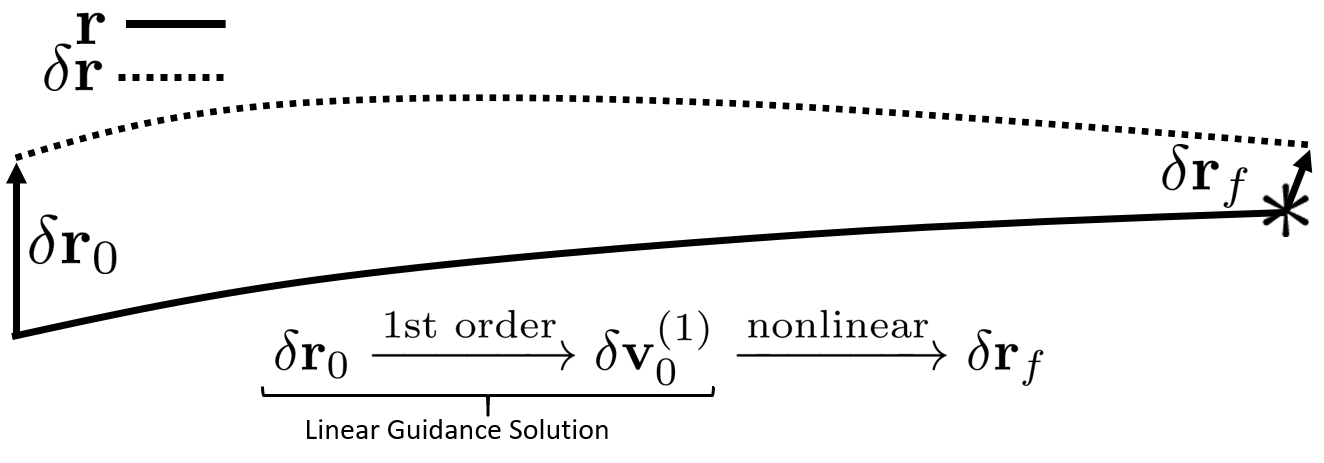}
    \caption{Notional depiction of miss distance in a rendezvous calculated with linearized dynamics.}
    \label{fig:rend_notional}
\end{figure}

\begin{figure}[hbt!]
     \centering
     \begin{subfigure}[b]{0.45\textwidth}
         \centering
         \includegraphics[width=\textwidth]{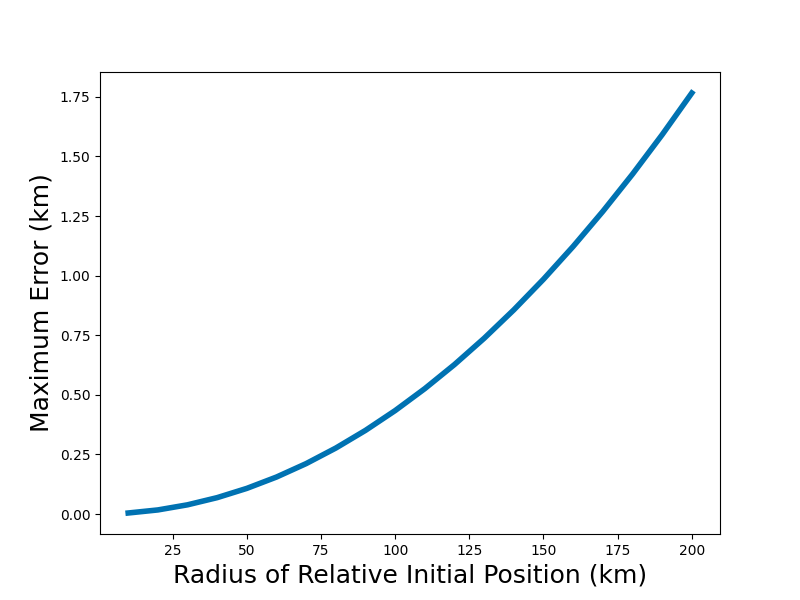}
         \caption{Initial position vs. Max miss distance in ISS example.}
     \end{subfigure}
     \hfill
     \begin{subfigure}[b]{0.45\textwidth}
         \centering
         \includegraphics[width=\textwidth]{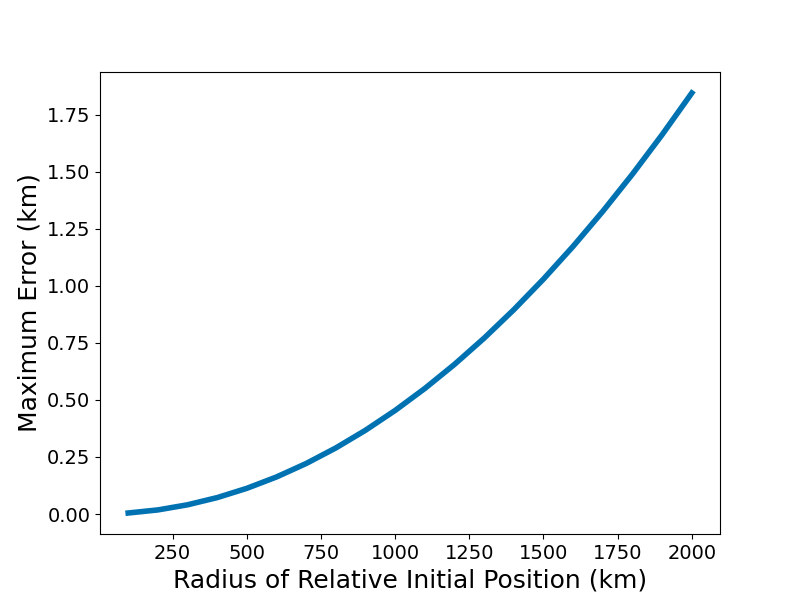}
         \caption{Initial position vs. Max miss distance in NRHO example.}
     \end{subfigure}\\
        \caption{Maximum miss distance of an impulsive rendezvous implemented using linear guidance.}
        \label{fig:rend_max_error_plots}
\end{figure}

\begin{figure}[hbt!]
     \centering
     \begin{subfigure}[b]{0.45\textwidth}
         \centering
         \includegraphics[width=\textwidth]{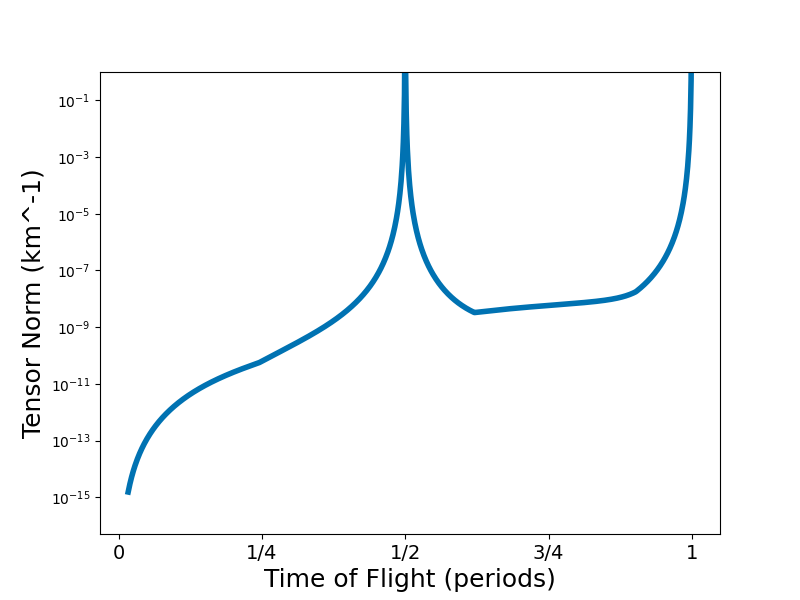}
         \caption{Rendezvous tensor norm vs. Time of flight in ISS example.}
     \end{subfigure}
     \hfill
     \begin{subfigure}[b]{0.45\textwidth}
         \centering
         \includegraphics[width=\textwidth]{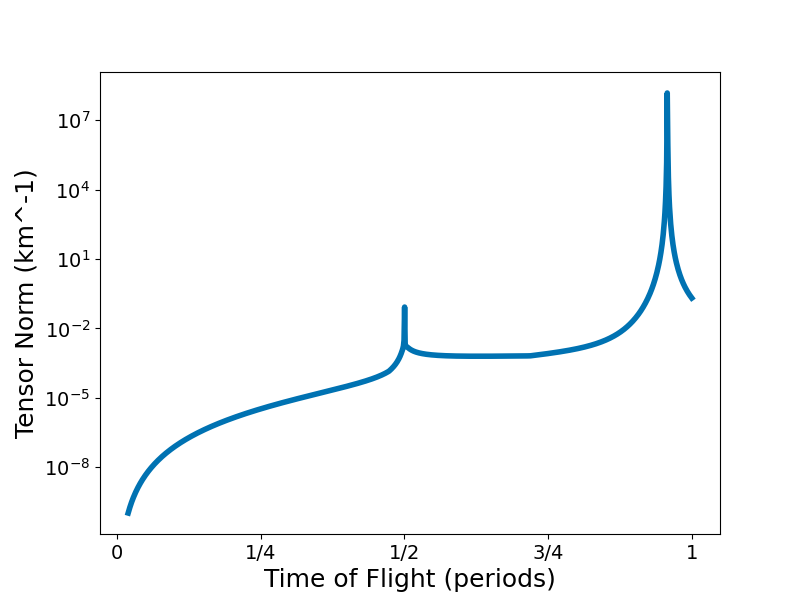}
         \caption{Rendezvous tensor norm vs. Time of flight in NRHO example.}
     \end{subfigure}\\
        \caption{Plots of tensor norm associated with rendezvous miss distance as transfer duration is varied.}
        \label{fig:rend_tnorm_plots}
\end{figure}

Before, we were considering the simplified transfer problem, where a satellite begins at the reference orbit and attempts to transfer to another position at a later time. Another related problem is the impulsive rendezvous problem. A satellite begins at some position $\delta\mathbf{r}_0$ relative to the reference orbit and attempts to reach the origin (reference orbit) at some time later. The linear solution to this problem is given:
\begin{equation}
    \delta\mathbf{v}_0^{(1)}=-(\mathbf{\Phi^\mathbf{r}_\mathbf{v}})^{-1}\mathbf{\Phi^\mathbf{r}_\mathbf{r}}\delta\mathbf{r}_0
    \label{eq:rend_1st}
\end{equation}
The final position error when the first-order solution for the initial velocity is used is
\begin{equation}
    \Vert \delta\mathbf{r}_f\Vert_2\leq  \Vert \mathbf{F}^{(1)}\Vert_2 \Vert\delta\mathbf{r}_0\Vert_2^2+\mathcal{O}(\Vert\delta\mathbf{r}_0\Vert_2^3)
    \label{eq:rend-bound}
\end{equation}
where the tensor comes from substituting the first-order transformation from $\delta\mathbf{r}_0$ to $\delta\mathbf{v}_0$ (Eq.~\ref{eq:rend_1st}) into any part of the second-order Taylor expansion that takes initial velocity as input
\begin{equation}
    (\mathbf{F}^{(1)})^i_{j,k}=\frac{1}{2}\left((\mathbf{\Psi}^\mathbf{r}_\mathbf{rr})^i_{j,k}-(\mathbf{\Psi}^\mathbf{r}_\mathbf{rv})^i_{j,l}((\mathbf{\Phi}^\mathbf{r}_\mathbf{v})^{-1}\mathbf{\Phi}^\mathbf{r}_\mathbf{r})^l_k-(\mathbf{\Psi}^\mathbf{r}_\mathbf{vr})^i_{l,k}((\mathbf{\Phi}^\mathbf{r}_\mathbf{v})^{-1}\mathbf{\Phi}^\mathbf{r}_\mathbf{r})^l_j+(\mathbf{\Psi}^\mathbf{r}_\mathbf{vv})^i_{l,p}((\mathbf{\Phi}^\mathbf{r}_\mathbf{v})^{-1}\mathbf{\Phi}^\mathbf{r}_\mathbf{r})^l_j((\mathbf{\Phi}^\mathbf{r}_\mathbf{v})^{-1}\mathbf{\Phi}^\mathbf{r}_\mathbf{r})^p_k\right)
    \label{eq:F-tens}
\end{equation}
%We do not present the more general form or the initial velocity error term due to the complexity in writing out this expression. 
What we have in this case is an approximate bound for the error in the final position during a rendezvous as a function of the initial distance from the chief satellite being targeted for rendezvous. A notional depiction of miss distance in a rendezvous using linear guidance is presented in Fig. \ref{fig:rend_notional}.
%\FloatBarrier
Examples of the miss distance from using a linearized guidance solution for rendezvous under the same circumstances as previously described are presented in Fig. \ref{fig:rend_max_error_plots} with the modification that the distance being varied is the initial distance from the reference orbit rather than the final distance. The miss distance is on the order of two to three times that of the transfer that begins at the reference orbit (between one and two kilometers in this particular case). This is because there are more coupled second-order terms coming into play in this scenario than there are in the simple transfer guidance miss distance scenario presented earlier. Again, the error performance of the three methods for approximating the maximum error exhibit the same trends and scales as in the previous two examples and are omitted.
%\FloatBarrier
% \begin{figure}[hbt!]
%      \centering
%      \begin{subfigure}[b]{0.45\textwidth}
%          \centering
%          \includegraphics[width=\textwidth]{Figures/Rend/twoBodyRendError.png}
%          \caption{Percent error between methods in ISS example.}
%      \end{subfigure}
%      \hfill
%      \begin{subfigure}[b]{0.45\textwidth}
%          \centering
%          \includegraphics[width=\textwidth]{Figures/Rend/threeBodyRendError.png}
%          \caption{Percent error between methods in NRHO example.}
%      \end{subfigure}\\
%         \caption{Relative inaccuracy between methods of calculating maximum error and true upper bound for propagation of an impulsive transfer in two-body and restricted three-body dynamics.}
%         \label{fig:rend_perc_error_plots}
% \end{figure}
%\FloatBarrier
In Fig. \ref{fig:rend_tnorm_plots}, we examine the norm of the tensor defined in Eq. \ref{eq:F-tens}, which can be used to bound the error of rendezvous miss distance using Eq. \ref{eq:rend-bound}. Unsurprisingly, relative transfer singularities manifest as spikes at certain transfer times in this example as well since each tensor examined in the guidance problems described here has relied on the invertibility of the relevant block of the state transition matrix which is known to become singular at these times. Interestingly, even between these spikes, the tensor norm is not a smooth function of time. Just as is the case with matrix eigenvalues, tensor eigenvalues, upon which this tensor norm is based, do not always vary in a differentiable fashion as a parameter is varied and the eigenvalues pass by one another. 

%\FloatBarrier
We have demonstrated that the error associated with linear propagation of motion around some reference orbit as well as linear guidance in the vicinity of that orbit can be quantified using a tensor norm approach. In the cases examined, the tensor norm approach for bounding error was very accurate at small distances (off by fractions of a percent of the actual error), but degrades as the distances from the reference orbit increase. The advantage of the tensor norm approach is that it can be used to predict the error as a function of the distance scale using a polynomial expression, and does not need to be recomputed for difference scale values. Additionally, the calculation is dominated by calculation of the second-order state transition tensor which requires on the order of 36 times as many equations to be integrated as a single integration of the original dynamical system. Compare this cost with sampling methods which rely on propagating thousands of trajectories to yield similar levels of accuracy and it is evident that the tensor norm method is one to two orders of magnitude faster than sampling, and even faster if the second-order state-transition tensor is already being calculated for another purpose.

\subsection{Measurement Nonlinearity}
\label{sec:measurement}
To quantify the nonlinearity of a measurement function, it is possible to directly study the norm of a coefficient tensor from the Taylor series approximation. However, another related metric can offer more insight into the error introduced by employing linearization of the measurement function in estimation problems. Suppose we have an initial estimate $\mathbf{x}^{-}$ for some state vector $\mathbf{x}^*$, and we want to update that estimate given a single measurement $\mathbf{h}(\mathbf{x}^*)=\mathbf{z}$ that is completely certain (i.e. no sensor noise) but potentially underdetermined so that the dimension of $\mathbf{z}$ may be less than the dimension of the state $\mathbf{x}$). Under an additional assumption on the homoscedasticity of the prior state covariance, the extended Kalman update reduces to the estimate 
\begin{equation}
    \mathbf{x}^+=\mathbf{x}^- + \left(\frac{\partial \mathbf{h}}{\partial \mathbf{x}}\right)^\dagger_{\mathbf{x}=\mathbf{x}^-}\left(\mathbf{z}-\mathbf{h}(\mathbf{x}^-)\right)
\end{equation}
where $(\cdot)^\dagger$ is the Moore-Penrose pseudoinverse. Expanding $\mathbf{h}(\mathbf{x}^*)=\mathbf{z}$ about $\mathbf{x}^-$ up to second order where we define $\delta\mathbf{x}=\mathbf{x}^*-\mathbf{x}^-$, we obtain
\begin{align}
    \mathbf{x}^+&=\mathbf{x}^- + \left(\frac{\partial \mathbf{h}}{\partial \mathbf{x}}\right)^\dagger_{\mathbf{x}=\mathbf{x}^-}\left(\mathbf{h}(\mathbf{x}^-)+\left(\frac{\partial \mathbf{h}}{\partial \mathbf{x}}\right)_{\mathbf{x}=\mathbf{x}^-}\delta\mathbf{x}+\frac{1}{2}\left(\frac{\partial^2 \mathbf{h}}{\partial \mathbf{x}^2}\right)_{\mathbf{x}=\mathbf{x}^-}\delta\mathbf{x}^2    -\mathbf{h}(\mathbf{x}^-)\right)+\mathcal{O}(\delta\mathbf{x}^3)\\
    &=\mathbf{x}^- + \boldsymbol{\Pi}_\mathbf{H}\delta\mathbf{x} +\frac{1}{2}\left(\frac{\partial \mathbf{h}}{\partial \mathbf{x}}\right)^\dagger_{\mathbf{x}=\mathbf{x}^-}\left(\frac{\partial^2 \mathbf{h}}{\partial \mathbf{x}^2}\right)_{\mathbf{x}=\mathbf{x}^-}\delta\mathbf{x}^2+\mathcal{O}(\delta\mathbf{x}^3)\\
    &=\mathbf{x}^*-\boldsymbol{\Pi}^\perp_\mathbf{H}\delta\mathbf{x}+\frac{1}{2}\left(\frac{\partial \mathbf{h}}{\partial \mathbf{x}}\right)^\dagger_{\mathbf{x}=\mathbf{x}^-}\left(\frac{\partial^2 \mathbf{h}}{\partial \mathbf{x}^2}\right)_{\mathbf{x}=\mathbf{x}^-}\delta\mathbf{x}^2+\mathcal{O}(\delta\mathbf{x}^3)
\end{align}
where $\boldsymbol{\Pi}_\mathbf{H}$ is the projection of the state into the linearly observable subspace under $\mathbf{h}$---a projection into the row space of the Jacobian of the measurement function
\begin{equation}
    \mathbf{H}=\frac{\partial \mathbf{h}}{\partial \mathbf{x}}
\end{equation}
and $\boldsymbol{\Pi}_\mathbf{H}^\perp$ is the projection onto the unobservable subspace of the state given by the right null space of $\mathbf{H}$.
Thus, the error in the observable subspace of the state using the linearization for the update is determined up to second-order by the tensor
\begin{equation}
    \bar{\mathbf{H}}=\left(\frac{\partial \mathbf{h}}{\partial \mathbf{x}}\right)^\dagger\frac{\partial^2 \mathbf{h}}{\partial \mathbf{x}^2}
\end{equation}
with components
\begin{equation}
    \bar{H}^i_{jk}=\left(\left(\frac{\partial \mathbf{h}}{\partial \mathbf{x}}\right)^\dagger\right)^i_l \left(\frac{\partial^2 \mathbf{h}}{\partial \mathbf{x}^2}\right)^l_{jk}
    \label{eq:hbar}
\end{equation}
It is possible to compute another tensor generalizing Eq. \ref{eq:hbar} that depends specifically on some prior state covariance and sensor noise, but the approach we take above allows us to compare measurement functions without getting into specifics of the sensor noise covariance or prior covariance.
Under the homoscedasticity and noise-free sensor assumption above, the 2-norm of $\bar{\mathbf{h}}$ evaluated at some estimate of the state gives the maximum level of error in the observable subspace of the state when using the linear update as a function of how far off that estimate is from the truth
\begin{equation}
    \Vert\boldsymbol{\Pi}_\mathbf{H}(\mathbf{x}^*-\mathbf{x}^+)\Vert_2\leq \frac{1}{2}\Vert\bar{\mathbf{H}}\Vert_2 \Vert\delta\mathbf{x}\Vert_2^2+\mathcal{O}(\delta\mathbf{x}^3)
\end{equation}
Using the norm of $\bar{\mathbf{H}}$, we can examine how well differing measurement models perform in terms of linearization error in the context of a linear estimation algorithm such as an extended Kalman filter. Note that the error in the unobservable subspace is the same no matter the choice of coordinates for an observation, because the unobservable subspace is the same regardless of the particular coordinate choice for the measurement model.

\begin{figure}
    \centering
    \includegraphics[width=.4\textwidth]{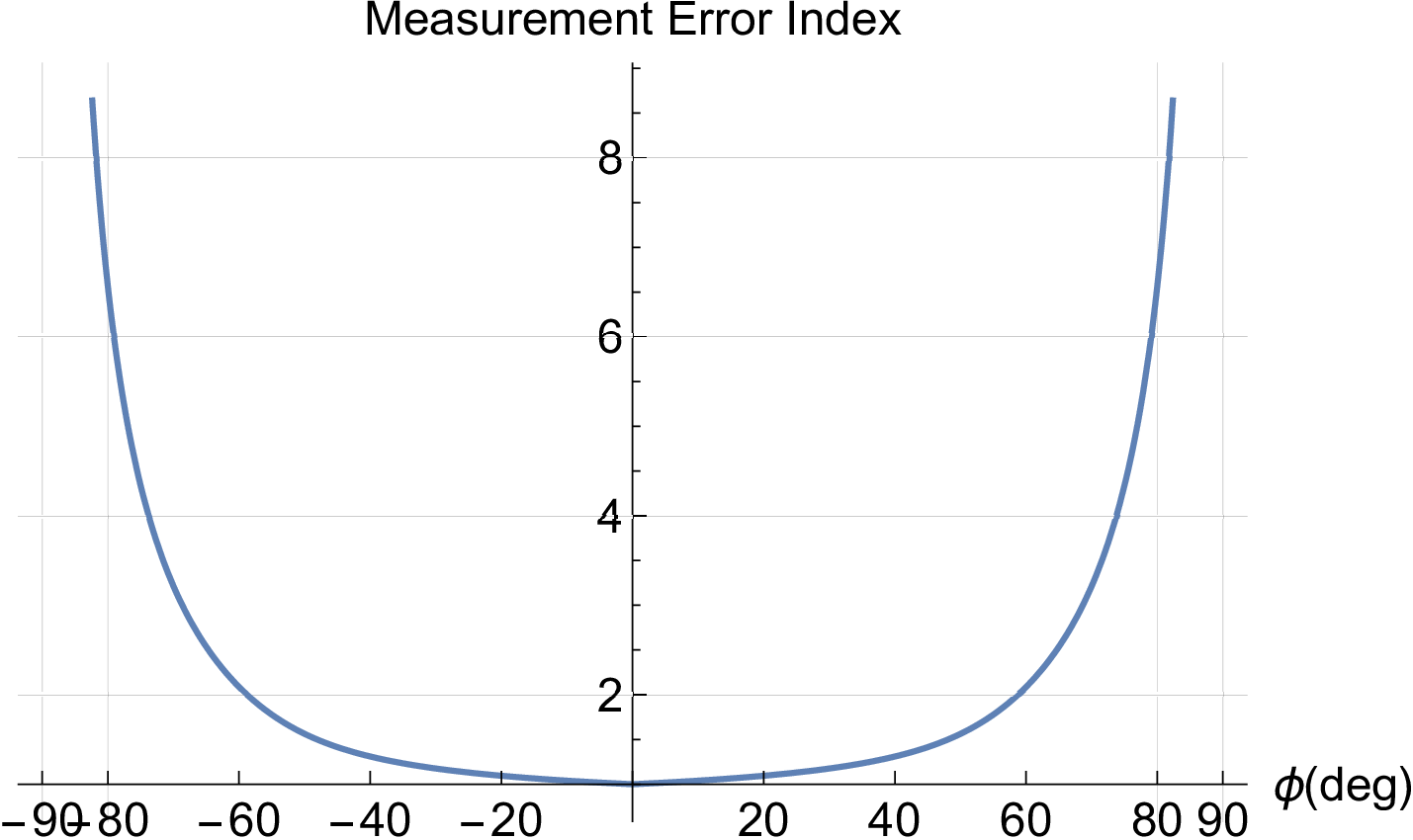}
    \caption{$\Vert \bar{\mathbf{H}}_a\Vert_2$ associated with angle measurements as a function of the estimate latitude $\phi$}
    \label{fig:angle-meas}
\end{figure}

\begin{figure}[hbt!]
     \centering
%\begin{picture}(220,180)
\begin{picture}(210,160)
\centering
\put(0,0){ \includegraphics[height=160pt]{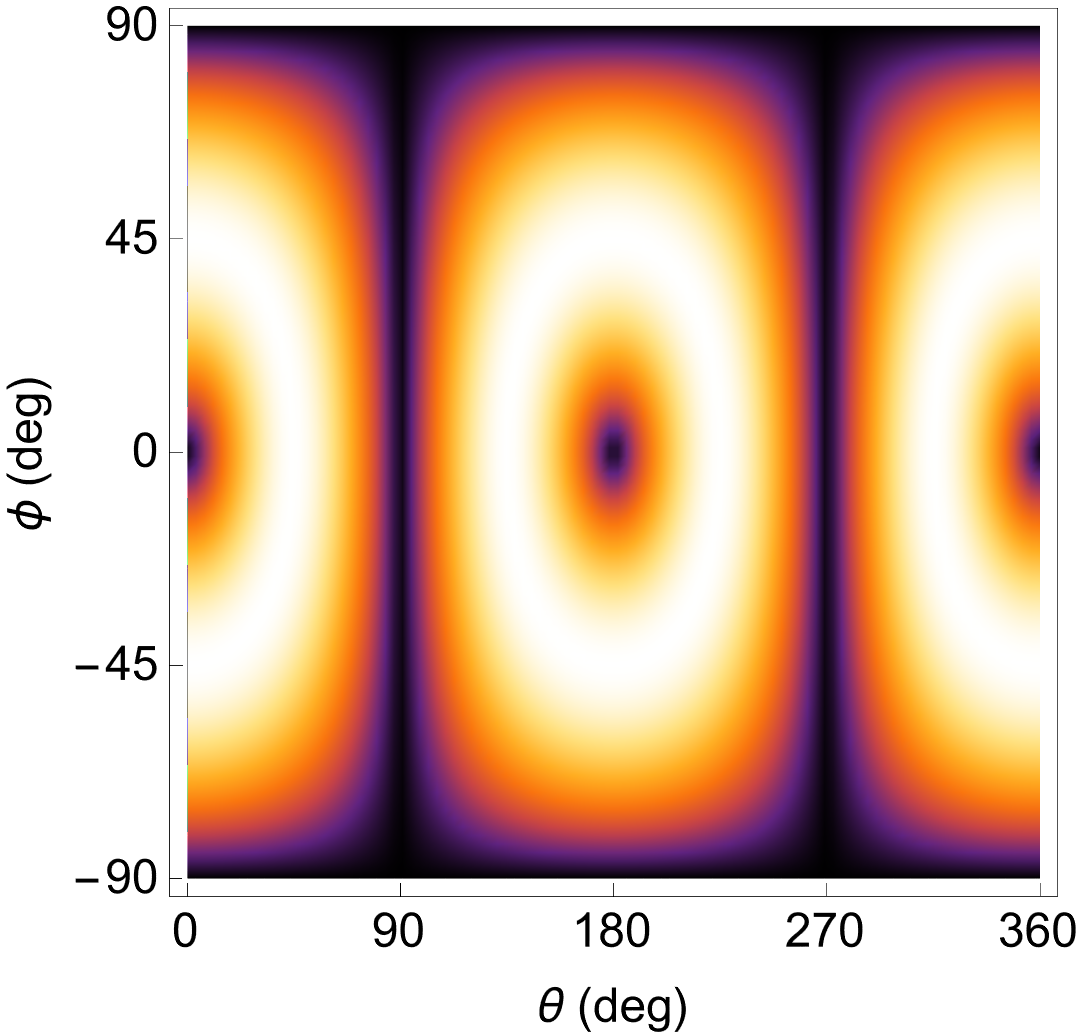}}
\put(165,20){ \includegraphics[height=120pt]{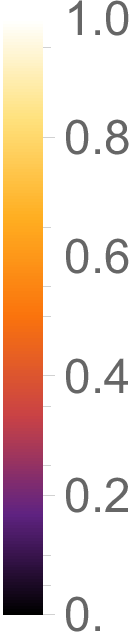}}
\put(160,170){\begin{minipage}[t]{50pt}
\begin{align*}
\Vert\bar{\mathbf{H}}_u\delta\mathbf{r}^2\Vert_2
\end{align*}
 \end{minipage}}
\end{picture}

\caption{The value of $\Vert\bar{\mathbf{H}}_u\delta\mathbf{r}^2\Vert_2$ from Eq.~\ref{eq:unit_error} is plotted for initial estimate errors on the unit sphere.}
\label{fig:unit_meas}
\end{figure}

For example, two common models of an optical measurement are given by 
\begin{align}
    \mathbf{h}_a(\mathbf{r})&=[\theta, \phi]^T = \left[\tan^{-1}\left(\frac{y}{x}\right),\sin^{-1}\left(\frac{z}{\Vert\mathbf{r}\Vert_2}\right)\right]^T\\
    \mathbf{h}_u(\mathbf{r})&=\frac{\mathbf{r}}{\Vert\mathbf{r}\Vert_2} = \frac{1}{\sqrt{x^2+y^2+z^2}}[x,y,z]^T
\end{align}
where $\mathbf{r}=[x,y,z]^T$, $\theta$ is the azimuthal angle, and $\phi$ is the elevation angle. We can evaluate the angle and unit vector models for their error performance in a linear estimation algorithm by evaluating $\Vert \bar{\mathbf{H}}\Vert_2$ at various estimates for $\mathbf{r}$. We use position vectors confined to the unit sphere for comparison. The angle measurement function should exhibit the same error behavior at any point with the same value of $\phi$ due to symmetry. In Fig.~\ref{fig:angle-meas}, we plot the value of $\Vert \bar{\mathbf{H}}\Vert_2$ associated with the measurement function $\mathbf{h}_a$ as a function of $\phi$ for an estimate on the unit sphere with that given latitude angle. The value of $\Vert \bar{\mathbf{H}}_a\Vert_2$ varies from 1 to $\infty$, with a minimum at $\phi=0$ and singularity at $\phi=\pm90\deg$.
On the other hand, the unit vector measurement model possesses spherical symmetry, and so it can be expected that the value of $\Vert \bar{\mathbf{H}}_u\Vert_2$ is constant when evaluated at estimates anywhere along the unit sphere. For an estimate at $\mathbf{r}=[1,0,0]^T$, the value of the norm of $\bar{\mathbf{H}}_u\delta\mathbf{r}^2$ is plotted in Fig.~\ref{fig:unit_meas} for $\mathbf{r}=[\cos\theta\cos\phi,\sin\theta\cos\phi,\sin\phi]^T$. One can show that the value is given by the expression
\begin{equation}
    \Vert\bar{\mathbf{H}}_u\delta\mathbf{r}^2\Vert_2=4 \cos ^2(\theta ) \cos ^2(\phi ) \left(\sin ^2(\theta ) \cos ^2(\phi )+\sin ^2(\phi )\right)
    \label{eq:unit_error}
\end{equation}
which has a maximum of $1$ along a 1-dimensional subspace including the point $\delta\mathbf{r}=\sqrt{2}/2[1,1,0]^T$. Thus, $\Vert \bar{\mathbf{H}}_u\Vert_2$ is 1 everywhere on the unit ball. At zero degree latitudes, both models perform equally well in terms of linearization error. However, at latitudes greater than around $80\deg$, the effects of linearization error on an extended Kalman update are around an order of magnitude greater for the angle measurement model as compared to the unit vector model which has uniformly superior linearization error performance. This, of course, comes at the cost of an increase in the dimension of the measurement model for the unit vector model and a minor additional computational burden when dealing with larger matrices. When minor increases in computational cost are not important, the unit vector approach provides superior linearization error performance for linear estimation algorithms. 
While, in practice no measurement is perfectly certain, this stylized problem bears a close resemblance to the stressing situation of an extended Kalman filter in which the measurement noise covariance is very low when compared with the state error covariance and a filter can become smug. In fact, a principled method for underweighting measurements in an extended Kalman filter makes use of the bounds on the two norm (as presented in Sec. \ref{sec:bounds-2}) of the second-order measurement partial derivative tensor \cite{zanetti2010underweighting}. While this bound could likely be made tighter using the 2-norm of the measurement partial derivatives, it is unlikely that the actual 2-norm would be as easy and quick to compute in closed form for general measurement models as the bound on the 2-norm tends to be. As such, this approach is not pursued further here, and this method is highlighted mostly as a means to compare the efficacy of different measurement models for an estimation problem in much the same way the nonlinearity index is employed to compare different state representations for a given dynamical system.

\subsection{Nonlinearity Index}
\label{sec:nonlin_index}
The nonlinearity index was initially introduced to quantify the level of nonlinearity of a dynamical system as described in various coordinates \cite{junkins1997karman,junkins2004nonlinear}. It has also been employed in the study of distribution propagation under nonlinear dynamics \cite{park2006nonlinear}. Additionally, it has been used to justify the use of calibration and decalibration processes to improve accuracy in settings that employ linearized dynamics \cite{sinclair2014calibration}. To name a few other applications, the nonlinearity index has been calculated to compare coordinate system choices for optimal control \cite{kellyindex}, relative motion \cite{alfriend2005evaluation}, and aerodynamics \cite{abdallah2013measuring, tapolcai2017aircraft}. Originally, the nonlinearity index was defined in a manner that depends on the scale of the perturbation to the reference trajectory, and also relies on numerical integration of a large number of randomly sampled perturbed trajectories. 

The nonlinearity index $\nu(t_f,t_0)$ associated with a trajectory was originally presented in terms of the Frobenius norm ($\Vert\cdot\Vert_F$) of the difference between the STM about some reference orbit and the STM about $N$ other nearby orbits each with an initial deviation of their state vector $\dxo_i$ from the reference initial state $\xo$ \cite{junkins1997karman, junkins2004nonlinear}. That is, it was defined as
\begin{equation}
    \nu(t_f,t_0)=\sup_{i=1..N}\left(\frac{\left\Vert \ph{\xo+\dxo_i}{f}{0}-\ph{\xo}{f}{0}\right\Vert_F}{\left\Vert\ph{\xo}{f}{0}\right\Vert_F}\right)
    \label{eqn:adhoc}
\end{equation}
Typically, the chosen initial deviations lie uniformly at random on a sphere centered at zero, but the choice of size for the sphere and the sampling are ad hoc. We present a similar alternative formulation that is scale-free and does not rely on sampling trajectories with initial deviations from the reference. First, we introduce new notation to clarify the dependence of the nonlinearity index on the scale of perturbation $r$ and the choice of output norm $a$ and input norm $b$. We define
\begin{equation}
    \nu^a_{b,r}(t_f,t_0)=\sup_{\left\Vert\dxo\right\Vert_b=r}\left(\frac{\left\Vert \ph{\xo+\dxo}{f}{0}-\ph{\xo}{f}{0}\right\Vert^a_b}{\left\Vert\ph{\xo}{f}{0}\right\Vert_b^a\Vert\dxo\Vert_b}\right)
    \label{eqn:with_scale}
\end{equation}
To excise the scale parameter $r$ from equation \eqref{eqn:with_scale}, we can take the limit as $r$ approaches zero. Then
\begin{equation}
    \nu_{b}^{a}(t_f,t_0)=\lim_{r\rightarrow 0}\nu^a_{b,r}(t_f,t_0)=\frac{\Vert\ps{\xo}{f}{0}\Vert^a_b}{\Vert\ph{\xo}{f}{0}\Vert^a_b}
    \label{eqn:psi_index}
\end{equation}
A scale-free version of the original nonlinearity index due to Junkins\cite{junkins1997karman} does not match this form where the two norms in the quotient are the same as one another, but takes the similar form
\begin{equation}
    \nu^*(t_f,t_0)=\frac{\Vert\ps{\xo}{f}{0}\Vert^F_2}{\Vert\ph{\xo}{f}{0}\Vert_F}
    \label{eqn:junkins-like}
\end{equation}
As was discussed in Sec. \ref{sec:frob-inf}, a similar nonlinearity index was presented in terms of an easily computed upper bound on the (Frobenius, $\infty$)-norm of the second-order state transition tensor \cite{losacco2024low}
\begin{equation}
    \nu^\Box(t_f,t_0)= \frac{\Vert\boldsymbol{\Psi}^\infty(\mathbf{x}(t_0);t_f,t_0)\Vert_F}{\Vert\ph{\xo}{f}{0}\Vert_F}\geq\frac{\Vert\ps{\xo}{f}{0}\Vert^F_\infty}{\Vert\ph{\xo}{f}{0}\Vert_F}
    \label{eqn:armellin-index}
\end{equation}
where the matrix $(.)^\infty$ was defined in Eq. \ref{eq:inf-mat}. This nonlinearity index is not unitarily invariant and varies based on rotations of the coordinate system. We denote this index with a superscript box because of its relation to the infinity norm and the need to disambiguate from a norm/index induced by the infinity norm. Note that while the expression above is identical to that presented in \cite{losacco2024low}, the language of differential algebra was used and the connection to the induced norm bound, while implied in the original work, is our own formalized observation.

In any of the above cases, by taking a quotient of a tensor norm and a matrix norm, we have split up two quantities that are each obtained as the solutions to constrained optimization problems. %Splitting the nonlinearity index into two separate optimizations tends to make it easier to compute
%The original index in Eq.~\ref{eqn:adhoc} and its scale-free counterpart in Eq.~\ref{eqn:junkins-like} are quantifying how much the Frobenius norm of the state transition matrix can change as the reference orbit about which it is calculated is perturbed and normalizing this with respect to the Frobenius norm of the state transition matrix at the reference orbit. This answers the question: how much does the linearization change percentage-wise (if we multiply by one hundred) as we change the reference orbit? This is a perfectly valid and interesting question to ask, though there are other equally valid and interesting questions.
Alternatively, Jenson and Scheeres presented a formulation of the nonlinearity index called Tensor Eigenpair Measure of Nonlinearity (TEMoN) in which the nonlinearity index is framed in terms of a single constrained optimization problem, rather than two independent optimizations \cite{jenson2022semianalytical}. This innovative work was the first to introduce the study of tensor eigenvalues to the astrodynamics as well as guidance, navigation, and control communities. In their paper, two optimization problems appear. One nonlinearity measure of interest is given by the optimization
\begin{equation}
    \aleph_{m,R}=\max_{\dxo\in\mathcal{B}_R}\left\vert \frac{\mathbf{C}^{(3)}\dxom{3}+...+\mathbf{C}^{(m)}\dxom{m}}{\mathbf{C}^{(2)}\dxom{2}}\right\vert
    \label{eqn:aleph}
\end{equation}
On the other hand, TEMoN is given by the optimization
\begin{equation}
    \tau_{m,R}=\max_{\dxo\in\mathcal{B}_R}\left\vert \frac{\mathbf{C}^{(m)}\dxom{m}}{\mathbf{C}^{(2)}\dxom{2}}\right\vert
    \label{eqn:TEMON}
\end{equation}
and the previous measure of nonlinearity from Eq.~\ref{eqn:aleph} is upper bounded by the sum of individual TEMoN optimizations for each order up to $m$
\begin{equation}
    \aleph_{m,R} \leq \sum_{l=3}^m \tau_{m,R}
\end{equation}
Jenson and Scheeres proposed a method to calculate TEMoN involving finding Z-eigenvalues of an $m+2$ order tensor derived from the method of Lagrange multipliers. However, this method requires finding more than just the eigenvector associated with the largest eigenvalue. Instead, the algorithm relies on methods other than symmetric higher-order power iteration to find as many eigenvectors as possible, checking each to find which yields the largest TEMoN. In general, finding all eigenvectors, especially those that are not associated with the largest eigenvalues, is a difficult task that has no guarantees of finding all eigenpairs and can be computationally intensive \cite{benson2019computing}. %Shortly, we will propose a method for computing TEMoN using the easier-to-compute largest D-eigenvalue of a related tensor.
We propose another nonlinearity index in a similar spirit called the D-Eigenvalue Measure of Nonlinearity (DEMoN). The development of this index also motivates a simpler method for computing TEMoN using the largest eigenvalue of some tensor computable by shifted symmetric higher order power iteration. DEMoN-m is given as %Being Scale-Free, also makes this measure of nonlinearity Super-Friendly, as the user does not need to specify or determine the size of the ball over which the optimization is performed. 
\begin{equation}
    \mu^{(m)}(t_f,t_0)=\sup_{\left\Vert\dxo\right\Vert_2=1}\left(\frac{\left\Vert\psm{m}{\xo}{f}{0}\delta\mathbf{x}^m(t_0)\right\Vert_2}{\left\Vert\ph{\xo}{f}{0}\dxo\right\Vert_2}\right)
    \label{eqn:demon}
\end{equation}
If the optimization were conducted over a ball of radius $R$ rather than the unit ball, the result would be scaled by a factor of $R^{m-1}$ given the homogeneity of multilinear operators. This makes the index scale-free in a sense unlike any index that involves a nonhomogeneous polynomial or full nonlinear function. 

All of the indices presented so far besides the original sampling-based index are scale-free in a similar fashion. The most significant order to calculate is the $m=2$ case since this accounts for the deviations from linearity of dynamics closest to the reference orbit. The index $\mu^{(2)}$ quantifies how large the norm of the quadratic term in the Taylor series can become relative to the norm of the linear term given the same input vector. In order to calculate DEMoN, a D-eigenvalue problem, and a rescaling of the resulting eigenvector can be substituted back into Eq.~\ref{eqn:demon}. This nonlinearity index can be calculated using only the largest D-eigenpair, avoiding the difficult task of computing all eigenpairs encountered by the original implementation of TEMoN. The vector $\dxo$ that maximizes the expression in Eq.~\ref{eqn:demon} is parallel to the maximizer of the related optimization
\begin{equation}\sup_{\delta\mathbf{x}_0^T\boldsymbol{\Phi}^T\boldsymbol{\Phi}\delta\mathbf{x}_0=1}\left\Vert\psm{m}{\xo}{f}{0}\delta\mathbf{x}^m(t_0)\right\Vert_2^2
    \label{eqn:demon1}
\end{equation}
where the arguments of the shorthand $\boldsymbol{\Phi}$ and $\delta\mathbf{x}_0$ are understood to match the arguments in the rest of the expression. Refer to Sec. \ref{sec:d-norm} to compute the vector $\delta\mathbf{x}_0^*$ maximizing the above expression using the theory of D-eigenvalues where $\mathbf{D}=\boldsymbol{\Phi}^T\boldsymbol{\Phi}$ and the square root of $\mathbf{D}$ need not be calculated using the Cholesky factorization, but instead is given by $\boldsymbol{\Phi}$. With this vector $\delta\mathbf{x}_0^*$ computed, DEMoN is the expression in Eq.~\ref{eqn:demon} evaluated at a normalized (to the unit ball) $\delta\mathbf{x}_0^*$. %It is important to note that only the largest D-eigenvalue needs to be calculated (which can be calculated from the largest Z-eigenvalue of a related tensor), and so this index is substantially easier to calculate than existing methods for calculating TEMoN which requires computation of eigenpairs other than just the relatively easy to find largest eigenpair. 
A disadvantage of this approach is that the sum of the squares of DEMoN for various orders plus unity (to account for the second-order CGT term) gives a more relaxed bound on the nonlinearity index $\aleph_{m,R}$ than the sum of TEMoN terms. On the other hand, a sum of DEMoN terms provides a bound on the nonlinearity index 
\begin{equation}
    \beth_{m,R}=\max_{\dxo\in\mathcal{B}_R} \frac{\left\Vert\frac{1}{2!}\mathbf{\Psi}^{(2)}\dxom{2}+...+\frac{1}{m!}\mathbf{\Psi}^{(m)}\dxom{m}\right\Vert_2}{\left\Vert\mathbf{C}^{(2)}\dxom{2}\right\Vert_2}
    \label{eqn:beth}
\end{equation}
\begin{equation}
    \beth_{m,R} \leq \sum_{l=2}^m \frac{R^{m-1}}{m!}\mu^{(m)}
    \end{equation}
While $\aleph_{m,R}$ gives the \textbf{error of the norm} of the $m$-th order approximation of the final deviation vector relative to the norm of the linear approximation, this index denoted $\beth_{m,R}$ (Beth from the Hebrew alphabet) gives the \textbf{norm of the error} of the $m$-th order approximation of the final deviation vector relative to the norm of the linear approximation. Both TEMoN and DEMoN as well as $\aleph$ and $\beth$ are valid indices to compute with differing meanings: an error of a norm or a norm of an error. %However, DEMoN and $\beth$ are substantially easier to compute and approximate respectively. Additionally, t
The nature of TEMoN and the index $\aleph$ implies that they may miss/understate strong nonlinearities with contributions in different directions than the linear contribution. Orthogonal contributions like these do not change the norm of the final deviation vector up to a linear approximation. For concreteness, we present the following extreme stylized example. Consider two dynamical systems in $\mathbb{R}^2$ and reference trajectories (one is the primed dynamical system and the other is the unprimed dynamical system) such that the state transition matrix and second-order state transition tensor for either system have all zero entries except for 
\begin{align}
    \Phi^1_1=1,& \quad\Psi^2_{1,1}=1\\
    (\Phi')^1_1=1,& \quad(\Psi')^1_{1,1}=1\\
\end{align}
In this situation, the values of DEMoN and $\beth$ are identical for the two systems (DEMoN has a value of 1). On the other hand, $\tau_{3,R}=0$ and $\tau'_{3,R}=\frac{1}{2}R$, while $\tau_{4,R}=\frac{1}{4}R^2=\tau'_{4,R}$, so that $\aleph_{m,R}=\frac{1}{4}R^2$ and $\aleph'_{m,R}=\frac{1}{2}R+\frac{1}{4}R^2$ for any $m$ higher than 4. While these two systems seem equally nonlinear, TEMoN and $\aleph$ index one of these systems as being substantially more nonlinear than the other. %This gives a potential advantage in interpretation in addition to the computational advantages for considering $\beth$ over $\aleph$ and DEMoN over TEMoN.
% A more realistic example is the scalar range measurement function. Consider the range from the origin in two dimensions
% \begin{equation}
%     h_r(\mathbf{r})=h_r(x,y)=\sqrt{x^2+y^2}
% \end{equation}
% where $\mathbf{r}=[x, y]^T$.
% Since this is a scalar measurement function, the first-order partial derivatives simply form a gradient which is (0,1)-tensor, and the second-order partial derivatives form the Hessian which is a (0,2)-tensor. That is, the partial derivative tensors are one order lower than in the case of a vector valued nonlinear function, and some of the formalism about tensor eigenvalues simplifies to that of regular matrix eigenvalues. The first order partial simplifies to
% \begin{equation}
%     \left(\frac{\partial h_r}{\partial \mathbf{r}}\right)^T=\hat{\mathbf{r}}
% \end{equation}
% which is the unit vector pointing from the origin to the position. On the other hand, the Hessian is
% \begin{equation}
%     \frac{\partial^2 h_r}{\partial \mathbf{r}^2}=\frac{\hat{\boldsymbol{\theta}}\hat{\boldsymbol{\theta}}^T}{\Vert \mathbf{r}\Vert_2}
% \end{equation}
% where $\hat{\boldsymbol{\theta}}$ is a direction orthogonal to $\hat{\boldmath{r}}$. Since the Hessian of the range function is an outer product of $\hat{\boldsymbol{\theta}}$ with itself, the vector $\hat{\boldmath{r}}$ given by the gradient is in the null space of the Hessian. Thus, TEMoN-3 will be zero for every point in the plane excluding the origin even though the range measurement is clearly nonlinear everywhere and the second-order term in the Taylor expansion is generically nonzero. 
The fact that TEMoN characterizes an error of a norm rather than a norm of an error and can underestimate the nonlinearity of some functions gives a potential advantage for considering $\beth$ over $\aleph$ and DEMoN over TEMoN. However, the error in the norm interpretation of TEMoN may be of interest still. Having derived DEMoN, we see that a slight modification of this D-eigenvalue methodology allows for computation of TEMoN. The main issue is that the higher-order CGT tensors are not necessarily even order or convex on the unit ball, so symmetric higher-order power iteration may not be convergent as is the case in calculating DEMoN which deals only in squares of tensors and not products of different tensors. As a result, TEMoN must be calculated using a shifted symmetric higher-order power iteration. For details, see Sec. \ref{appendix:temon}. 

Since each of the nonlinearity indices above results from some constrained optimization problem, each of them has some corresponding argument or input vector that leads to the nonlinearity index. These vectors which solve the constrained optimization problem can be used to identify possible direction for splitting of a Gaussian mixture model. $\nu^\Box$ was used for just this purpose in the context of splitting Gaussian mixtures in a coordinate axis aligned fashion \cite{losacco2024low}. Using vectors associated with other nonlinearity indices such as $\nu^2_2$ or $\mu^{(2)}$ can give a single optimal non-axis-aligned splitting direction. The former focuses on reducing absolute error from second-order nonlinearities, while the latter focuses on reducing percentage error from second-order nonlinearities. Similar constrained optimizations associated with D-eigenvalues can be posed where $\mathbf{D}=\mathbf{P}^{-1}$ is the inverse of the covariance matrix for the Gaussian distribution being split. These problems can be interpreted as finding the direction in which nonlinearity is greatest given inputs that are equally likely. Though not phrased in the tensor norm formalism we have presented, the vector that solves the constrained maximization problem associated with the (Frobenius, $\mathbf{P}^{-1}$)-norm of the second-order partial derivative tensor of a flow or measurement function has been employed for choosing a maximally nonlinear Gaussian mixture splitting direction \cite{tuggle2020model}.

\begin{figure}[hbt!]
    \centering
    \includegraphics[width=.45\textwidth]{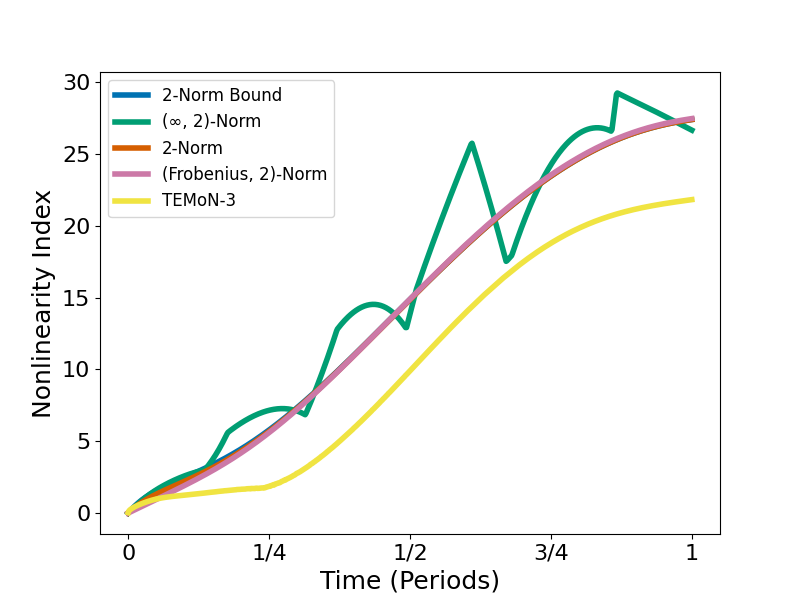}
    \includegraphics[width=.45\textwidth]{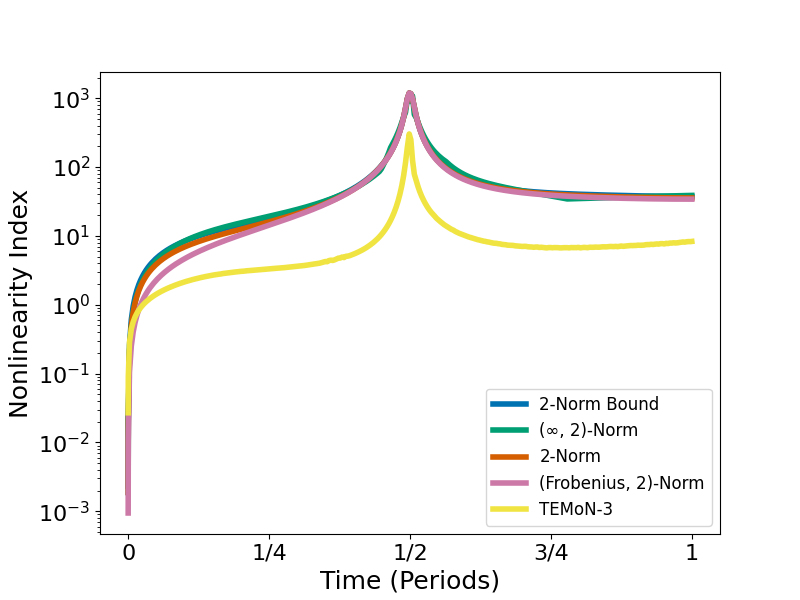}\\
    \includegraphics[width=.45\textwidth]{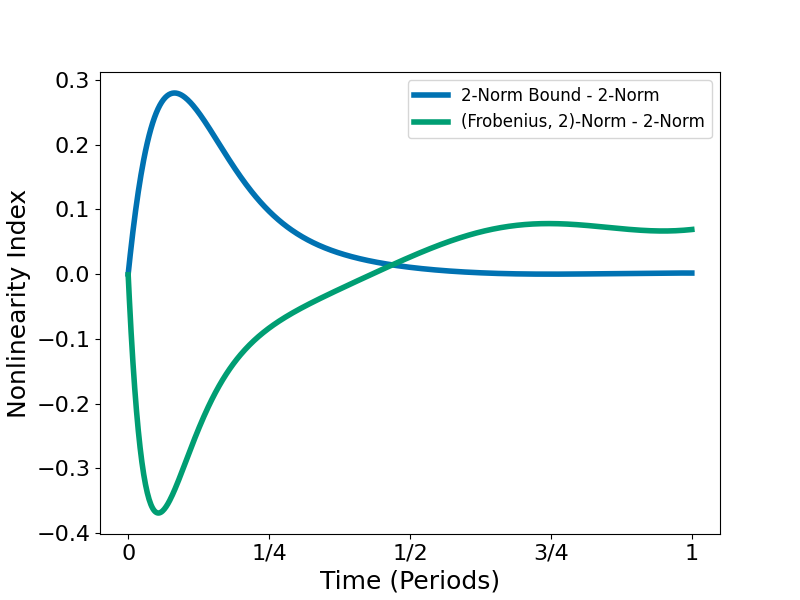}
    \includegraphics[width=.45\textwidth]{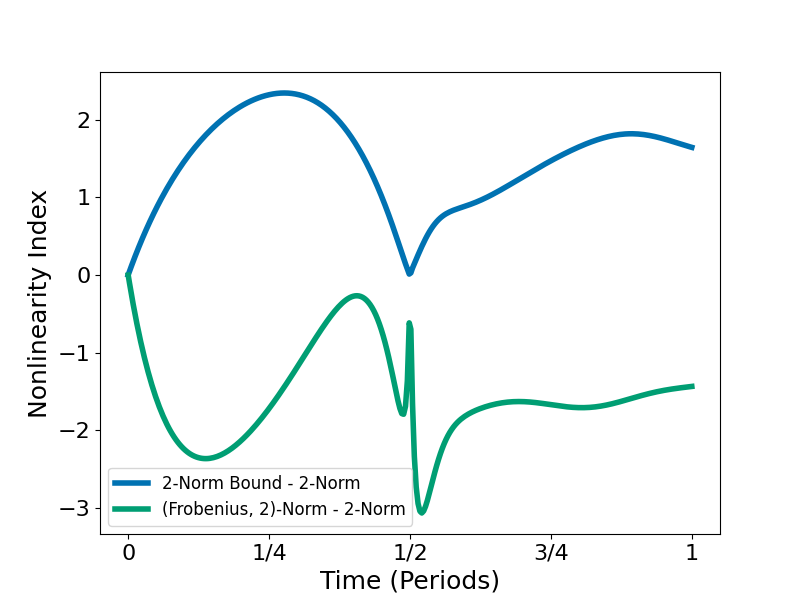}\\
    \includegraphics[width=.45\textwidth]{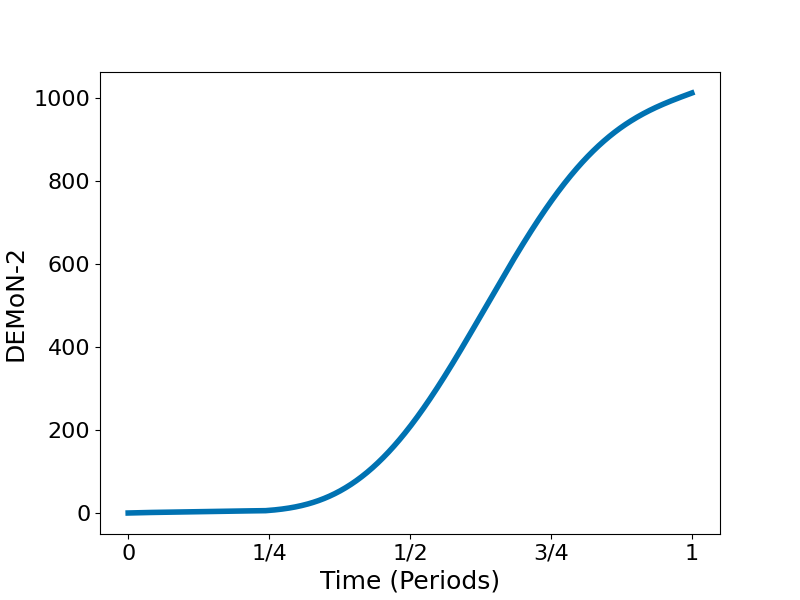}
    \includegraphics[width=.45\textwidth]{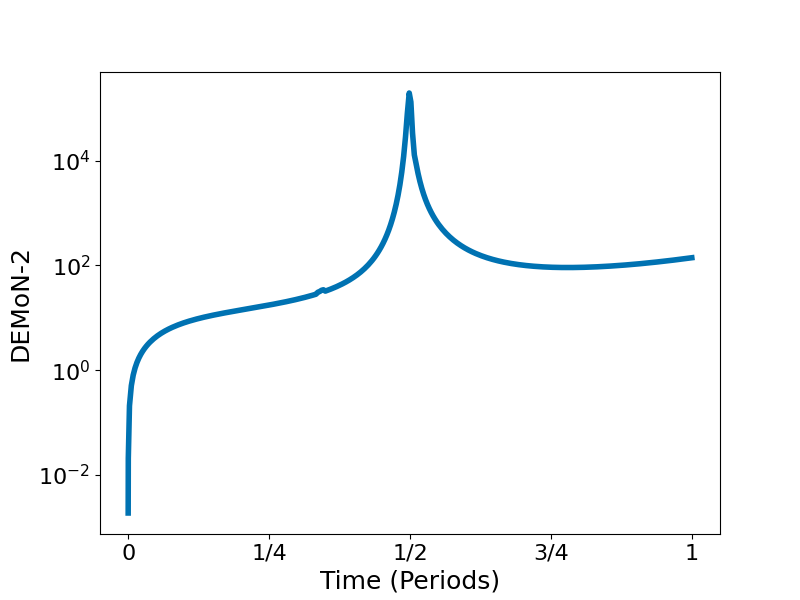}
    \caption{The nonlinearity indices associated with nondimensional circular two-body motion on the left and the Earth-Moon circular restricted three-body problem Gateway NRHO on the right.}
    \label{fig:nonlin}
\end{figure}

With the nonlinearity indices derived and presented above, we show two examples with a selection of the nonlinearity indices. First, we present the nonlinearity index associated with nondimensional circular two-body motion. Two-body dynamics were defined in Eq. \ref{eqn:two-body}. To further simplify the problem and make interpretation of the nonlinearity index simpler, we nondimensionalize the two-body problem using a reference circular orbit with semimajor axis $a$ to form the nondimensional state vector $\mathbf{x}'=[\mathbf{r}'^T, \mathbf{v}'^T]^T$.
\begin{align}
    \frac{d}{d\tau}\mathbf{x}'=\begin{bmatrix}
        \mathbf{v}'\\
        \dfrac{-\mathbf{r}'}{\Vert\mathbf{r}'\Vert_2^3}
    \end{bmatrix}
\end{align}
where $\Vert\mathbf{r}'\Vert_2^3$ denotes the cube of the 2-norm. The states and times are related by the following coordinate changes
\begin{align}
    \mathbf{r}'=\mathbf{r}/a\\
    \tau=\omega t\\
    \mathbf{v}'=\mathbf{v}/(a\omega)
\end{align}
where $\omega$ is the mean motion of the circular orbit with semimajor axis $a$. The initial conditions we employ for the generic planar nondimensional circular orbit are $\mathbf{x}_0=[1,0,0,0,1,0]^T$. A number of nonlinearity indices are presented on the left in Fig. \ref{fig:nonlin}. These include DEMoN-2 on its own due to scale considerations, TEMoN-3, and the nonlinearity indices associated with the 2-, (Frobenius, 2)-, and ($\infty$, 2)-norms as well as the flattening based 2-norm bound of the second-order state transition tensor. For the 2-norm bound based nonlinearity index, we take the quotient of the 2-norm bound in Sec. \ref{sec:bounds-2} applied to the second-order state transition tensor with the 2-norm of the state transition matrix. Since the values of the 2-norm, (Frobenius, 2)-norm, and 2-norm bound are very similar, we break out the difference between the 2-norm nonlinearity index and these two other indices on their own. In addition to the two-body example, we also use the Gateway NRHO beginning at apolune to demonstrate how these nonlinearity indices grow in the circular restricted three body problem on the right hand side of Fig. \ref{fig:nonlin}. The sensitivity to initial conditions exhibited by the three-body problem manifests itself in a much higher growth rate of the nonlinearity indices than in the two-body case. We also see a peak in nonlinearity near perilune which is consistent with previous analysis with TEMoN \cite{jenson2022semianalytical}.
Note that the 2-norm bound does in-fact give an upper bound on the 2-norm based nonlinearity index in both cases for all times, and quite a tight one at that. One of the main takeaways is that in the case of two- and three-body motion, these nonlinearity indices exhibit very similar patterns of growth over time and may be quite interchangeable-especially in the case of the 2-norm, (Frobenius, 2)-norm, and 2-norm upper bound based nonlinearity indices. The main exception is that the overall magnitudes differ between TEMoN, DEMoN, and the nonlinearity indices associated with various norms of the second-order state transition tensor. This implies that the simpler to compute nonlinearity indices such as the (Frobenius, 2)-norm or the ($\infty, 2$)-norm might be the best to use, though the ($\infty, 2$)-norm may be advised against due to its non-smoothness. 

%\FloatBarrier

%TEMoN can be calculated in a similar fashion. Consider the optimization problem (negativity problem?)

% An analog to the nonlinearity index $\aleph_{m,R}$ approximated/bounded by TEMoN can be computed exactly using the majorization minimization technique from Sec. \ref{sec:mm-opt}. While $\aleph_{m,R}$ gives the error of the norm of the $m$-th order approximation of the final deviation vector relative to the norm of the linear approximation, this index denoted $\beth_{m,R}$ (Beth from the Hebrew alphabet) will give the norm of the error of the $m$-th order approximation of the final deviation vector relative to the norm of the linear approximation.
% \begin{equation}
%     \beth_{m,R}=\max_{\dxo\in\mathcal{B}_R} \frac{\left\Vert\frac{1}{2!}\mathbf{\Psi}^{(2)}\dxom{2}+...+\frac{1}{m!}\mathbf{\Psi}^{(m)}\dxom{m}\right\Vert_2}{\left\Vert\mathbf{C}^{(2)}\dxom{2}\right\Vert_2}
%     \label{eqn:beth}
% \end{equation}
% Note that as opposed to this expression, Eq.~\ref{eqn:aleph} has terms in the numerator with contributions from the first order state transition matrix.

%\subsection{Non-Gaussian Summary Statistics}

\section{Conclusion}
In this paper, we have developed and synthesized a number of existing approaches to quantifying nonlinearity of a function or flow of a dynamical system into a single framework based on tensor operator norms. This led to an easy-to-compute and easy-to-interpret scale-free nonlinearity index that most closely parallels the original nonlinearity index presented over two decades ago, as well as two easier-to-compute versions of the modern tensor-eigenvalue-based nonlinearity index (TEMoN and DEMoN). One of these tensor eigenvalue based measures of nonlinearity (DEMoN) quantifies the norm of a deviation from linearity rather than a deviation of norm predicted by linear and nonlinear models as has been done previously. In addition to presenting nonlinearity indices as a means to quantify nonlinearity of a dynamical system and compare different formulations of a system, we have shown that the same tools for defining a nonlinearity index can be employed to estimate the error resulting from commonly used linear algorithms applied to a nonlinear system. In particular, we demonstrated the use of tensor norms to approximately bound the error in linearized guidance methods for stationkeeping or rendezvous. Additionally, we presented a method for analyzing the maximum error in the extended Kalman filter update step when a measurement is much more precise than the current state estimate. This error approximation methodology can be applied with some modifications to many linear algorithms in guidance, navigation, and control. The methodology is more efficient than sampling-based methods for quantifying error and gives an analytical grasp on error performance for varying distance scales from the point at which linearization was applied. Though we do not pursue it here, future work could involve using this same methodology in a statistical context for quantifying non-Gaussianity by examining norms of tensors related to higher-order moments of a distribution.

\section*{Acknowledgements}
Jackson Kulik was supported by the Department of Defense National Defense Science and Engineering Graduate Fellowship. M. Ruth was supported by the National Science Foundation Graduate Research Fellowship under Grant No. DGE-1650441. The authors would like to thank Robin Armstrong, Keith LeGrand, Spencer Boone, David Cunningham, and Erica Jenson for helpful discussion and comments. Code examples associated with this paper can be found at https://github.com/SIOSlab/STMInt.

\bibliography{sample}

\begin{thebibliography}{57}
\newcommand{\enquote}[1]{``#1''}
\providecommand{\natexlab}[1]{#1}
\providecommand{\url}[1]{\texttt{#1}}
\providecommand{\urlprefix}{URL }
\expandafter\ifx\csname urlstyle\endcsname\relax
  \providecommand{\doi}[1]{\discretionary{}{}{}https://doi.org/#1}\else
  \providecommand{\doi}[1]{\discretionary{}{}{}\urlstyle{rm}\url{https://doi.org/#1}}\fi

\bibitem[{Clohessy and Wiltshire(1960)}]{clohessy1960terminal}
Clohessy, W.~H., and Wiltshire, R.~S., \enquote{Terminal guidance system for satellite rendezvous,} \emph{Journal of the Aerospace Sciences}, Vol.~27, No.~9, 1960, pp. 653--658.
\newblock \doi{10.2514/8.8704}.

\bibitem[{Mullins(1992)}]{mullins1992initial}
Mullins, L.~D., \enquote{Initial value and two point boundary value solutions to the Clohessy-Wiltshire equations,} \emph{Journal of the Astronautical Sciences}, Vol.~40, No.~4, 1992, pp. 487--501.

\bibitem[{Schmidt(1981)}]{schmidt1981kalman}
Schmidt, S.~F., \enquote{The Kalman filter-Its recognition and development for aerospace applications,} \emph{Journal of Guidance and Control}, Vol.~4, No.~1, 1981, pp. 4--7.
\newblock \doi{10.2514/3.19713}.

\bibitem[{Junkins(1997)}]{junkins1997karman}
Junkins, J.~L., \enquote{von karman lecture: Adventures on the interface of dynamics and control,} \emph{Journal of Guidance, Control, and Dynamics}, Vol.~20, No.~6, 1997, pp. 1058--1071.
\newblock \doi{10.2514/2.4176}.

\bibitem[{Junkins and Singla(2004)}]{junkins2004nonlinear}
Junkins, J.~L., and Singla, P., \enquote{How nonlinear is it? A tutorial on nonlinearity of orbit and attitude dynamics,} \emph{The Journal of the Astronautical Sciences}, Vol.~52, No.~1, 2004, pp. 7--60.
\newblock \doi{10.1007/bf03546420}.

\bibitem[{Jenson and Scheeres(2023)}]{jenson2022semianalytical}
Jenson, E.~L., and Scheeres, D.~J., \enquote{Semianalytical Measures of Nonlinearity Based on Tensor Eigenpairs,} \emph{Journal of Guidance, Control, and Dynamics}, Vol.~46, No.~4, 2023, p. 638–653.
\newblock \doi{10.2514/1.g006760}, \urlprefix\url{http://dx.doi.org/10.2514/1.g006760}.

\bibitem[{Losacco et~al.(2024)Losacco, Fossà, and Armellin}]{losacco2024low}
Losacco, M., Fossà, A., and Armellin, R., \enquote{Low-Order Automatic Domain Splitting Approach for Nonlinear Uncertainty Mapping,} \emph{Journal of Guidance, Control, and Dynamics}, Vol.~47, No.~2, 2024, pp. 291--310.
\newblock \doi{10.2514/1.g007271}.

\bibitem[{Tuggle(2020)}]{tuggle2020model}
Tuggle, K.~E., \enquote{Model selection for Gaussian mixture model filtering and sensor scheduling,} Ph.D. thesis, 2020.

\bibitem[{Jenson and Scheeres(2024)}]{jenson2024bounding}
Jenson, E.~L., and Scheeres, D.~J., \enquote{Bounding nonlinear stretching about spacecraft trajectories using tensor eigenpairs,} \emph{Acta Astronautica}, Vol. 214, 2024, pp. 159--166.
\newblock \doi{10.1016/j.actaastro.2023.10.013}.

\bibitem[{Qi et~al.(2008)Qi, Wang, and Wu}]{qi2008d}
Qi, L., Wang, Y., and Wu, E.~X., \enquote{D-eigenvalues of diffusion kurtosis tensors,} \emph{Journal of Computational and Applied Mathematics}, Vol. 221, No.~1, 2008, pp. 150--157.
\newblock \doi{10.1016/j.cam.2007.10.012}.

\bibitem[{Boodram et~al.(2022{\natexlab{a}})Boodram, Boone, and McMahon}]{boone2022directional}
Boodram, O., Boone, S., and McMahon, J., \enquote{Directional State Transition Tensors for Capturing Dominant Nonlinear Effects in Orbital Dynamics,} \emph{Journal of Guidance, Control, and Dynamics}, 2022{\natexlab{a}}, pp. 1--12.
\newblock \doi{10.2514/1.G006910}.

\bibitem[{Kelly et~al.(2020)Kelly, Sinclair, and Majji}]{kellyhigher}
Kelly, P., Sinclair, A.~J., and Majji, M., \enquote{Higher Order State Transition Tensors for Keplerian Motion Using Universal Parameters,} \emph{Astrodynamics Specialist Conference}, 2020.

\bibitem[{Rein and Tamayo(2016)}]{rein2016second}
Rein, H., and Tamayo, D., \enquote{Second-order variational equations for N-body simulations,} \emph{Monthly Notices of the Royal Astronomical Society}, Vol. 459, No.~3, 2016, pp. 2275--2285.
\newblock \doi{10.1093/mnras/stw644}.

\bibitem[{Bani~Younes(2019)}]{bani2019exact}
Bani~Younes, A., \enquote{Exact Computation of High-Order State Transition Tensors for Perturbed Orbital Motion,} \emph{Journal of Guidance, Control, and Dynamics}, Vol.~42, No.~6, 2019, pp. 1365--1371.
\newblock \doi{10.2514/1.g003897}.

\bibitem[{Cunningham and Russell(2023)}]{cunningham2023interpolated}
Cunningham, D., and Russell, R.~P., \enquote{An Interpolated Second-Order Relative Motion Model for Gateway,} \emph{The Journal of the Astronautical Sciences}, Vol.~70, No.~4, 2023, p.~26.
\newblock \doi{10.1007/s40295-023-00393-9}.

\bibitem[{Boone and McMahon(2021)}]{boone2021orbital}
Boone, S., and McMahon, J., \enquote{Orbital Guidance Using Higher-Order State Transition Tensors,} \emph{Journal of Guidance, Control, and Dynamics}, Vol.~44, No.~3, 2021, pp. 493--504.
\newblock \doi{10.2514/1.g005493}.

\bibitem[{Kulik et~al.(2023)Kulik, Clark, and Savransky}]{kulik2023state}
Kulik, J., Clark, W., and Savransky, D., \enquote{State Transition Tensors for Continuous-Thrust Control of Three-Body Relative Motion,} \emph{Journal of Guidance, Control, and Dynamics}, 2023, pp. 1--10.
\newblock \doi{10.2514/1.g007311}.

\bibitem[{Park and Scheeres(2006)}]{park2006nonlinear}
Park, R.~S., and Scheeres, D.~J., \enquote{Nonlinear mapping of Gaussian statistics: theory and applications to spacecraft trajectory design,} \emph{Journal of guidance, Control, and Dynamics}, Vol.~29, No.~6, 2006, pp. 1367--1375.
\newblock \doi{10.2514/1.20177}.

\bibitem[{Majji et~al.(2008{\natexlab{a}})Majji, Junkins, and Turner}]{majji2008high}
Majji, M., Junkins, J.~L., and Turner, J.~D., \enquote{A high order method for estimation of dynamic systems,} \emph{The Journal of the Astronautical Sciences}, Vol.~56, No.~3, 2008{\natexlab{a}}, pp. 401--440.
\newblock \doi{10.1007/bf03256560}.

\bibitem[{Boodram et~al.(2022{\natexlab{b}})Boodram, Boone, and McMahon}]{boodramefficient}
Boodram, O., Boone, S., and McMahon, J., \enquote{Efficient Nonlinear Spacecraft Navigation Using Directional State Transition Tensors,} \emph{2022 AAS/AIAA Astrodynamics Specialist Conference AAS 22-670}, 2022{\natexlab{b}}.

\bibitem[{Kulik et~al.(2024)Kulik, Oertell, and Savransky}]{kulik2024overdetermined}
Kulik, J., Oertell, O., and Savransky, D., \enquote{Overdetermined Eigenvector Approach to Passive Angles-Only Relative Orbit Determination,} \emph{Journal of Guidance, Control, and Dynamics}, Vol.~47, 2024, pp. 1--9.
\newblock \doi{10.2514/1.g007744}.

\bibitem[{Rasotto et~al.(2016)Rasotto, Morselli, Wittig, Massari, Di~Lizia, Armellin, Valles, and Ortega}]{rasotto2016differential}
Rasotto, M., Morselli, A., Wittig, A., Massari, M., Di~Lizia, P., Armellin, R., Valles, C., and Ortega, G., \enquote{Differential algebra space toolbox for nonlinear uncertainty propagation in space dynamics,} 2016.

\bibitem[{Valli et~al.(2013)Valli, Armellin, Di~Lizia, and Lavagna}]{valli2013nonlinear}
Valli, M., Armellin, R., Di~Lizia, P., and Lavagna, M.~R., \enquote{Nonlinear mapping of uncertainties in celestial mechanics,} \emph{Journal of Guidance, Control, and Dynamics}, Vol.~36, No.~1, 2013, pp. 48--63.
\newblock \doi{10.2514/1.58068}.

\bibitem[{Wittig et~al.(2015)Wittig, Di~Lizia, Armellin, Makino, Bernelli-Zazzera, and Berz}]{wittig2015propagation}
Wittig, A., Di~Lizia, P., Armellin, R., Makino, K., Bernelli-Zazzera, F., and Berz, M., \enquote{Propagation of large uncertainty sets in orbital dynamics by automatic domain splitting,} \emph{Celestial Mechanics and Dynamical Astronomy}, Vol. 122, 2015, pp. 239--261.
\newblock \doi{10.1007/s10569-015-9618-3}.

\bibitem[{Servadio et~al.(2022)Servadio, Zanetti, and Armellin}]{servadio2022maximum}
Servadio, S., Zanetti, R., and Armellin, R., \enquote{Maximum a posteriori estimation of Hamiltonian systems with high order Taylor polynomials,} \emph{The Journal of the Astronautical Sciences}, Vol.~69, No.~2, 2022, pp. 511--536.
\newblock \doi{10.1007/s40295-022-00304-4}.

\bibitem[{Armellin and Di~Lizia(2018)}]{armellin2018probabilistic}
Armellin, R., and Di~Lizia, P., \enquote{Probabilistic optical and radar initial orbit determination,} \emph{Journal of Guidance, Control, and Dynamics}, Vol.~41, No.~1, 2018, pp. 101--118.
\newblock \doi{10.2514/1.g002217}.

\bibitem[{Di~Lizia et~al.(2014)Di~Lizia, Armellin, Bernelli-Zazzera, and Berz}]{di2014high}
Di~Lizia, P., Armellin, R., Bernelli-Zazzera, F., and Berz, M., \enquote{High order optimal control of space trajectories with uncertain boundary conditions,} \emph{Acta Astronautica}, Vol.~93, 2014, pp. 217--229.
\newblock \doi{10.1016/j.actaastro.2013.07.007}.

\bibitem[{Lizia et~al.(2008)Lizia, Armellin, and Lavagna}]{lizia2008application}
Lizia, P.~D., Armellin, R., and Lavagna, M., \enquote{Application of high order expansions of two-point boundary value problems to astrodynamics,} \emph{Celestial Mechanics and Dynamical Astronomy}, Vol. 102, 2008, pp. 355--375.
\newblock \doi{10.1007/s10569-008-9170-5}.

\bibitem[{Majji et~al.(2008{\natexlab{b}})Majji, Junkins, and Turner}]{majji2008jth}
Majji, M., Junkins, J., and Turner, J., \enquote{Jth moment extended Kalman filtering for estimation of nonlinear dynamic systems,} \emph{AIAA Guidance, Navigation and Control Conference and Exhibit}, American Institute of Aeronautics and Astronautics, 2008{\natexlab{b}}, p. 7386.
\newblock \doi{10.2514/6.2008-7386}.

\bibitem[{Zanetti et~al.(2010)Zanetti, DeMars, and Bishop}]{zanetti2010underweighting}
Zanetti, R., DeMars, K.~J., and Bishop, R.~H., \enquote{Underweighting nonlinear measurements,} \emph{Journal of guidance, control, and dynamics}, Vol.~33, No.~5, 2010, pp. 1670--1675.
\newblock \doi{10.2514/1.50596}.

\bibitem[{Golub(2013)}]{golub2013matrix}
Golub, G., \emph{Matrix computations}, JHU press, 2013.
\newblock \doi{10.56021/9781421407944}.

\bibitem[{Kolda and Bader(2006)}]{kolda2006matlab}
Kolda, T.~G., and Bader, B.~W., \enquote{MATLAB tensor toolbox,} Tech. rep., Sandia National Laboratories (SNL), Albuquerque, NM, and Livermore, CA~…, 2006.

\bibitem[{Lim(2005)}]{lim2005singular}
Lim, L.-H., \enquote{Singular values and eigenvalues of tensors: a variational approach,} \emph{1st IEEE International Workshop on Computational Advances in Multi-Sensor Adaptive Processing, 2005.}, IEEE, 2005, pp. 129--132.

\bibitem[{De~Lathauwer et~al.(2000)De~Lathauwer, De~Moor, and Vandewalle}]{de2000best}
De~Lathauwer, L., De~Moor, B., and Vandewalle, J., \enquote{On the best rank-1 and rank-(r 1, r 2,..., rn) approximation of higher-order tensors,} \emph{SIAM journal on Matrix Analysis and Applications}, Vol.~21, No.~4, 2000, pp. 1324--1342.

\bibitem[{Kolda and Mayo(2011)}]{kolda2011shifted}
Kolda, T.~G., and Mayo, J.~R., \enquote{Shifted power method for computing tensor eigenpairs,} \emph{SIAM Journal on Matrix Analysis and Applications}, Vol.~32, No.~4, 2011, pp. 1095--1124.
\newblock \doi{10.1137/100801482}.

\bibitem[{Kofidis and Regalia(2002)}]{kofidis2002best}
Kofidis, E., and Regalia, P.~A., \enquote{On the best rank-1 approximation of higher-order supersymmetric tensors,} \emph{SIAM Journal on Matrix Analysis and Applications}, Vol.~23, No.~3, 2002, pp. 863--884.
\newblock \doi{10.1137/s0895479801387413}.

\bibitem[{Boone and McMahon(2023)}]{boone2023directional}
Boone, S., and McMahon, J., \enquote{Directional State Transition Tensors for Capturing Dominant Nonlinear Effects in Orbital Dynamics,} \emph{Journal of Guidance, Control, and Dynamics}, Vol.~46, No.~3, 2023, pp. 431--442.
\newblock \doi{10.2514/1.g006910}.

\bibitem[{Kolda and Mayo(2014)}]{kolda2014adaptive}
Kolda, T.~G., and Mayo, J.~R., \enquote{An adaptive shifted power method for computing generalized tensor eigenpairs,} \emph{SIAM Journal on Matrix Analysis and Applications}, Vol.~35, No.~4, 2014, pp. 1563--1581.
\newblock \doi{10.1137/140951758}.

\bibitem[{Drakakis and Pearlmutter(2009)}]{drakakis2009calculation}
Drakakis, K., and Pearlmutter, B.~A., \enquote{On the calculation of the l2→ l1 induced matrix norm,} \emph{International Journal of Algebra}, Vol.~3, No.~5, 2009, pp. 231--240.

\bibitem[{Lewis(2010)}]{lewis2010top}
Lewis, A.~D., \enquote{A top nine list: Most popular induced matrix norms,} \emph{Queen’s University, Kingston, Ontario, Tech. Rep}, 2010, pp. 1--13.

\bibitem[{Li et~al.(2020)Li, Chen, Li, and Zhao}]{li2020eigenvalue}
Li, S., Chen, Z., Li, C., and Zhao, J., \enquote{Eigenvalue bounds of third-order tensors via the minimax eigenvalue of symmetric matrices,} \emph{Computational and Applied Mathematics}, Vol.~39, 2020, pp. 1--14.
\newblock \doi{10.1007/s40314-020-01245-0}.

\bibitem[{Vallado(2001)}]{vallado2001fundamentals}
Vallado, D.~A., \emph{Fundamentals of astrodynamics and applications}, Springer Science \& Business Media, 2001, Vol.~4, Chap.~1.

\bibitem[{Koon et~al.(2000)Koon, Lo, Marsden, and Ross}]{koon2000dynamical}
Koon, W.~S., Lo, M.~W., Marsden, J.~E., and Ross, S.~D., \enquote{Dynamical systems, the three-body problem and space mission design,} \emph{Equadiff 99: (In 2 Volumes)}, World Scientific, 2000, pp. 1167--1181.
\newblock \doi{10.1142/9789812792617_0222}.

\bibitem[{Nat(2019)}]{NationalAA2019}
\enquote{National Aeronautics and Space Administration (NASA) White Paper: Gateway Destination Orbit Model: A Continuous 15 Year NRHO Reference Trajectory,} 2019.

\bibitem[{Virtanen et~al.(2020)Virtanen, Gommers, Oliphant, Haberland, Reddy, Cournapeau, Burovski, Peterson, Weckesser, Bright, {van der Walt}, Brett, Wilson, Millman, Mayorov, Nelson, Jones, Kern, Larson, Carey, Polat, Feng, Moore, {VanderPlas}, Laxalde, Perktold, Cimrman, Henriksen, Quintero, Harris, Archibald, Ribeiro, Pedregosa, {van Mulbregt}, and {SciPy 1.0 Contributors}}]{2020SciPy-NMeth}
Virtanen, P., Gommers, R., Oliphant, T.~E., Haberland, M., Reddy, T., Cournapeau, D., Burovski, E., Peterson, P., Weckesser, W., Bright, J., {van der Walt}, S.~J., Brett, M., Wilson, J., Millman, K.~J., Mayorov, N., Nelson, A. R.~J., Jones, E., Kern, R., Larson, E., Carey, C.~J., Polat, {\.I}., Feng, Y., Moore, E.~W., {VanderPlas}, J., Laxalde, D., Perktold, J., Cimrman, R., Henriksen, I., Quintero, E.~A., Harris, C.~R., Archibald, A.~M., Ribeiro, A.~H., Pedregosa, F., {van Mulbregt}, P., and {SciPy 1.0 Contributors}, \enquote{{{SciPy} 1.0: Fundamental Algorithms for Scientific Computing in Python},} \emph{Nature Methods}, Vol.~17, 2020, pp. 261--272.
\newblock \doi{10.1038/s41592-019-0686-2}.

\bibitem[{Howell and Marchand(2005)}]{howell2005natural}
Howell, K.~C., and Marchand, B.~G., \enquote{Natural and non-natural spacecraft formations near the L1 and L2 libration points in the Sun--Earth/Moon ephemeris system,} \emph{Dynamical Systems}, Vol.~20, No.~1, 2005, pp. 149--173.
\newblock \doi{10.1080/1468936042000298224}.

\bibitem[{Griffith et~al.(2004)Griffith, Turner, Vadali, and Junkins}]{griffith2004higher}
Griffith, D., Turner, J., Vadali, S., and Junkins, J., \enquote{Higher order sensitivities for solving nonlinear two-point boundary-value problems,} \emph{AIAA/AAS Astrodynamics Specialist Conference and Exhibit}, American Institute of Aeronautics and Astronautics, 2004, p. 5404.
\newblock \doi{10.2514/6.2004-5404}.

\bibitem[{Junkins et~al.(2009)Junkins, Turner, and Majji}]{junkins2009generalizations}
Junkins, J.~L., Turner, J.~D., and Majji, M., \enquote{Generalizations and applications of the lagrange implicit function theorem,} \emph{The Journal of the Astronautical Sciences}, Vol.~57, No. 1-2, 2009, pp. 313--345.
\newblock \doi{10.1007/bf03321507}.

\bibitem[{Kulik and Savransky(2022)}]{kulik2022relative}
Kulik, J., and Savransky, D., \enquote{Relative Transfer Singularities and Multi-Revolution Lambert Uniqueness,} \emph{AIAA SCITECH 2022 Forum}, American Institute of Aeronautics and Astronautics, 2022, p. 0958.
\newblock \doi{10.2514/6.2022-0958}.

\bibitem[{Fitzgerald(1998)}]{fitzgerald1998pinch}
Fitzgerald, R.~J., \enquote{Pinch Points of Debris from a Satellite Breakup,} \emph{Journal of guidance, control, and dynamics}, Vol.~21, No.~5, 1998, pp. 813--815.
\newblock \doi{10.2514/2.7617}.

\bibitem[{Zhang and Zhou(2010)}]{zhang2010second}
Zhang, G., and Zhou, D., \enquote{A second-order solution to the two-point boundary value problem for rendezvous in eccentric orbits,} \emph{Celestial Mechanics and Dynamical Astronomy}, Vol. 107, 2010, pp. 319--336.
\newblock \doi{10.1007/s10569-010-9269-3}.

\bibitem[{Sinclair et~al.(2014)Sinclair, Sherrill, and Lovell}]{sinclair2014calibration}
Sinclair, A.~J., Sherrill, R.~E., and Lovell, T.~A., \enquote{Calibration of linearized solutions for satellite relative motion,} \emph{Journal of Guidance, Control, and Dynamics}, Vol.~37, No.~4, 2014, pp. 1362--1367.
\newblock \doi{10.2514/1.g000037}.

\bibitem[{Kelly et~al.(2022)Kelly, Arya, Junkins, and Majji}]{kellyindex}
Kelly, P., Arya, V., Junkins, J., and Majji, M., \enquote{Nonlinearity Index for State-Costate Dynamics of Optimal Control Problems,} \emph{Astrodynamics Specialist Conference}, 2022.

\bibitem[{Alfriend and Yan(2005)}]{alfriend2005evaluation}
Alfriend, K.~T., and Yan, H., \enquote{Evaluation and comparison of relative motion theories,} \emph{Journal of Guidance, Control, and Dynamics}, Vol.~28, No.~2, 2005, pp. 254--261.
\newblock \doi{10.2514/1.6691}.

\bibitem[{Abdallah et~al.(2013)Abdallah, Newman, and Omran}]{abdallah2013measuring}
Abdallah, A.~M., Newman, B.~A., and Omran, A.~M., \enquote{Measuring aircraft nonlinearity across aerodynamic attitude flight envelope,} \emph{AIAA Atmospheric Flight Mechanics (AFM) Conference}, American Institute of Aeronautics and Astronautics, 2013, p. 4985.
\newblock \doi{10.2514/6.2013-4985}.

\bibitem[{Tapolcai et~al.(2017)Tapolcai, Omran, and Newman}]{tapolcai2017aircraft}
Tapolcai, D.~P., Omran, A., and Newman, B., \enquote{Aircraft stall phenomenon analysis using nonlinearity index theory,} \emph{Aerospace Science and Technology}, Vol.~68, 2017, pp. 288--298.
\newblock \doi{10.1016/j.ast.2017.04.008}.

\bibitem[{Benson and Gleich(2019)}]{benson2019computing}
Benson, A.~R., and Gleich, D.~F., \enquote{Computing tensor Z-eigenvectors with dynamical systems,} \emph{SIAM Journal on Matrix Analysis and Applications}, Vol.~40, No.~4, 2019, pp. 1311--1324.
\newblock \doi{10.1137/18m1229584}.

\end{thebibliography}

\section{Appendix}

\subsection{Sufficiency of Partial Symmetry for Higher-Order Power Iteration}
\label{appendix:partial_sym}
We begin by showing an alternative expression for Eq.~\ref{eqn:power_iter} that references $\tilde{\mathbf{B}}$ rather than $\hat{\mathbf{B}}$. The contraction of the symmetrized tensor $\hat{\mathbf{B}}$ with $2m-1$ copies of the vector $\mathbf{x}$ can be computed instead as the average of all $2m$ different ways we could perform the $2m-1$ contractions with the non-symmetric tensor $\tilde{\mathbf{B}}$. Each of the $2m$ ways of performing the contraction can be described by which dimension is not contracted with:
\begin{equation}
    \left(\hat{\mathbf{B}}\mathbf{x}^{2m-1}\right)_l=\frac{1}{2m}\sum_{j=1}^{2m}\tilde{\mathbf{B}}_{i_1,...,i_{2m}}\delta^{i_j}_l\prod_{k\neq j} x^{i_k}
    \label{eqn:symmetrizing_contraction}
\end{equation}
where $\delta^{i_j}_l$ is the Kronecker delta given by contracting an index from then Euclidean metric tensor with the Euclidean inverse metric tensor.
The expression in Eq.~\ref{eqn:symmetrizing_contraction} is true regardless of any partial symmetry of $\tilde{\mathbf{B}}$. However, we note a few special properties of the tensors under consideration: $\mathbf{B}$ is partially symmetric under any permutation of the $m$ covariant indices, leading $\tilde{\mathbf{B}}$ to be symmetric under permutations of $i_1..i_m$ amongst themselves and $i_{m+1}..i_{2m}$ amongst themselves. Additionally, $\tilde{\mathbf{B}}$ is also symmetric under the permutation $(i_1...i_m,i_{m+1}..i_{2m})\rightarrow (i_{m+1}...i_{2m},i_1...i_m)$ as well as any permutations generated by the two mentioned above. As a result, for all $k\in\{1...2m\}$, there exists a permutation $\sigma$ such that the first element of $\sigma((i_1...i_m,i_{m+1}..i_{2m}))$ is $i_k$. 

An example with a $4$-dimensional tensor $\tilde{\mathbf{B}}$ is most instructive:
\begin{equation}
    \tilde{\mathbf{B}}_{i_1,i_2,i_3,i_4}=\tilde{\mathbf{B}}_{i_2,i_1,i_3,i_4}=\tilde{\mathbf{B}}_{i_3,i_4,i_1,i_2}=
    \tilde{\mathbf{B}}_{i_4,i_3,i_1,i_2}
\end{equation}
The first equality comes from partial symmetry of the covariant indices of $\mathbf{B}$, the second equality comes from the partial symmetry of the square of a tensor, and the fourth equality is generated by both types of permutations. As a result, for all $j$, the vectors below are equal:
\begin{equation}
\tilde{\mathbf{B}}_{i_1,...,i_{2m}}\prod_{k\neq j} x^{i_k}
\end{equation}
and thus
\begin{equation}
    \hat{\mathbf{B}}\mathbf{x}^{2m-1}=\tilde{\mathbf{B}}\mathbf{x}^{2m-1}
\end{equation}
This means we could perform symmetric higher-order power iteration successfully, only using our partially symmetric $\tilde{\mathbf{B}}$ rather than the fully symmetric $\hat{\mathbf{B}}$ and yield the exact same steps at each iteration. In fact, this opens up the avenue to an even simpler and more efficient version of power iteration that directly relies on $\mathbf{B}$:
\begin{equation}
    (\hat{\mathbf{B}}\mathbf{x}^{2m-1})_j=(\tilde{\mathbf{B}}\mathbf{x}^{2m-1})_j=(\mathbf{B}\mathbf{x}^{m-1})^{i_1}_j\boldsymbol{\delta}_{i_1,i_2}(\mathbf{B}\mathbf{x}^{m})^{i_2}
    \label{eqn:simplified_power}
\end{equation}
The rightmost expression can be obtained by first calculating the matrix $\mathbf{B}\mathbf{x}^{m-1}$, then calculating the vector $\mathbf{B}\mathbf{x}^{m}=(\mathbf{B}\mathbf{x}^{m-1})\mathbf{x}$ by a single matrix vector product. Then, a final matrix vector product between the two gives the result from Eq.~\ref{eqn:simplified_power}.

\subsection{Implementation of the $(2,\mathbf{D})$-Norm Computation}
\label{appendix:partial_sym_d}
The tensor $\hat{\mathbf{B}}_\mathbf{D}$ is a supersymmetric tensor since $\hat{\mathbf{B}}$ is a supersymmetric tensor. Thus, symmetric higher-order power iteration could be applied directly to calculate the Z-eigenpair of $\hat{\mathbf{B}}_\mathbf{D}$. In order to save on operations, if the Cholesky decomposition is performed to compute a lower triangular $\mathbf{D}^{1/2}$, then one iteration of symmetric higher-order power iteration can be taken more efficiently with the following steps. First, given the current iterate $\mathbf{x}$, compute 
\begin{equation}
    \mathbf{z}=\mathbf{D}^{-1/2}\mathbf{x}
\end{equation}
by forward substitution. Next, compute the vector
\begin{equation}
    \hat{\mathbf{B}}(\mathbf{D}^{-1/2}\mathbf{x})^{2m-1}=\hat{\mathbf{B}}\mathbf{z}^{2m-1}
\end{equation}
This can be performed efficiently as in Eq.~\ref{eqn:simplified_power}. Finally, the new iterate is given by back substitution to solve an upper triangular linear system to find 
\begin{equation}
    \hat{\mathbf{B}}_\mathbf{D}\mathbf{x}^{2m-1}=(\mathbf{D}^{1/2})^{-T}(\hat{\mathbf{B}}\mathbf{z}^{2m-1})
\end{equation}
which is the numerator of the next iterate of symmetric higher-order power iteration. Above, $\hat{\mathbf{B}}\mathbf{z}^{2m-1}$ is reinterpreted as a column vector to make the linear algebra operations to be performed more clear. The above sequence of operations can be derived in index notation:
\begin{align}
    (\hat{\mathbf{B}}_\mathbf{D}\mathbf{x}^{2m-1})_k&=\hat{B}_{j_1...j_{2m}}(\mathbf{D}^{-1/2})^{j_1}_{i_1}...(\mathbf{D}^{-1/2})^{j_{2m-1}}_{i_{2m-1}}(\mathbf{D}^{-1/2})^{j_{2m}}_{k}x^{i_1}...x^{i_{2m-1}}\\
    &=(\hat{\mathbf{B}}\mathbf{z}^{2m-1})_{j_{2m}}(\mathbf{D}^{-1/2})^{j_{2m}}_{k}
\end{align}
In this argument, the result is a row vector (covariant vector) while the input $\mathbf{x}$ is assumed to be a column vector (contravariant vector).

\subsection{Calculating TEMoN Using Shifted Symmetric Higher-Order Power Iteration}
\label{appendix:temon}
In order to find TEMoN, we solve two optimization problems over the unit ball, one for a positive case and one for a negative case:
\begin{align}
        \tau_{m}^\pm=\max_{\dxo\in\mathcal{B}_1} \frac{\pm\mathbf{C}^{(m)}\dxom{m}}{\mathbf{C}^{(2)}\dxom{2}}\\
        \tau_{m,R}=R^{m-2}\max(\tau_{m}^+,\tau_{m}^-)
    \label{eqn:TEMONpm}
\end{align}
In either case, we reformulate the problem with a change of coordinates
\begin{equation}
    \mathbf{y}=\boldsymbol{\Phi}\mathbf{x}
\end{equation}
and calculate the maximum Z-eigenpair of each of the two tensors given by
\begin{equation}
    \pm\left(\mathbf{C}^{(m)}_{\boldsymbol{\Phi}}\right)_{i_1...i_{m}}=\pm\mathbf{C}^{(m)}_{j_1...j_{m}}\left(\boldsymbol{\Phi}^{-1}\right)^{j_1}_{i_1}...\left(\boldsymbol{\Phi}^{-1}\right)^{j_{m}}_{i_{m}}
    \label{eqn:d-eig-tensor-temon}
\end{equation}
We can do so by performing shifted symmetric higher-order power iteration on the symmetrized tensor \cite{kolda2011shifted}. Let $\hat{\mathbf{C}}^{(m)}_{\boldsymbol{\Phi}}$ be the symmetrization of $\mathbf{C}^{(m)}_{\boldsymbol{\Phi}}$, and define
\begin{equation}
    \alpha=(m-1)\sum_{i_1,...,i_m}\left\vert \left(\hat{\mathbf{C}}^{(m)}_{\boldsymbol{\Phi}}\right)_{i_1...i_{m}}\right\vert
\end{equation}
Then, symmetric shifted higher-order power iteration is guaranteed to converge to (a typically large) eigenvector of $\mathbf{C}^{(m)}_{\boldsymbol{\Phi}}$. The iteration is given:
\begin{equation}
    \mathbf{y}_{n+1}^\pm=\frac{\pm\hat{\mathbf{C}}^{(m)}_{\boldsymbol{\Phi}}(\mathbf{y}^\pm_n)^{m-1}+\alpha \mathbf{y}_n^\pm}{\left\Vert\pm\hat{\mathbf{C}}^{(m)}_{\boldsymbol{\Phi}}(\mathbf{y}^\pm_n)^{m-1}+\alpha \mathbf{y}_n^\pm\right\Vert}_2
    \label{eqn:power_iter_temon}
\end{equation}
Using the inverse coordinate change on the converged to Z-eigenvector $\mathbf{y}^*$ we arrive at
\begin{equation}
    \mathbf{x}^\pm_*=\boldsymbol{\Phi}^{-1}\mathbf{y
}^\pm_*\end{equation}
Finally, our optimization is given by 
\begin{equation}
    \tau_{m}^\pm=\frac{\pm\mathbf{C}^{(m)}(\mathbf{x}^\pm_*)^m}{\mathbf{C}^{(2)}(\mathbf{x}^\pm_*)^2}
\end{equation}
The CGTs have fewer symmetries than the STTs and their squares. As a result, we do not present a method leveraging partial symmetry for computing TEMoN without performing symmetrization as we do in the case of DEMoN.

\end{document}